\documentclass[12pt]{article}
\title{Random walks and electric networks}
\author{Peter G. Doyle \and J. Laurie Snell}
\date{
\footnotesize
Version 3.02,  5 January 2000 \\
Copyright (C) 1999, 2000 Peter G. Doyle and J. Laurie Snell\\
Derived from work(s)\\
Copyright (C) 1984 The Mathematical Association of America\\
This work is freely redistributable under the terms of\\
the GNU General Public License\\
as published by the Free Software Foundation.\\
This work comes with ABSOLUTELY NO WARRANTY.
}
\usepackage{theorem}
\theorembodyfont{\rmfamily}
\newcommand{\prt}[2]{\section{#2}\label{part:#1}}
\newcommand{\chap}[2]{\subsection{#2}\label{chap:#1}}
\newcommand{\sect}[2]{\subsubsection{#2}\label{sect:#1}}
\newcommand{\partref}[1]{Section \ref{part:#1}}
\newcommand{\sectref}[1]{Section \ref{sect:#1}}
\newcommand{\book}{work}
\usepackage[dvips]{graphics}
\newcommand{\barefig}[1]{\includegraphics{figures/#1.ps}}
\newcommand{\fig}[2]{
\begin{figure}[htbp]
   \begin{center}
      \barefig{#2}
   \end{center}
   \caption{$\clubsuit$}
   \label{fig:#1}
\end{figure}
}
\newcommand{\prefig}[1]{\fig{#1}{#1}}
\newcommand{\figref}[1]{\ref{fig:#1}}
\newenvironment{touted}{}{}
\newcommand{\tout}[1]{\begin{touted}\textbf{#1}\ }
\newcommand{\toutend}{\end{touted}}
\newcommand{\proofstart}{\textbf{Proof.}\ }
\newcommand{\proofend}{$\diamondsuit$}
\newcommand{\whystart}{\textbf{Why.}\ }
\newcommand{\whyend}{$\heartsuit$}
\newcommand{\dt}[1]{\emph{#1}}
\newenvironment{exercises}{}{}
\newtheorem{exercise}{Exercise}[subsection]
\newcommand{\exstart}{\begin{exercises}}
\newcommand{\exend}{\end{exercises}}
\newcommand{\ex}[2]{\begin{exercise} \label{ex:#1} #2 \end{exercise}}
\newcommand{\exref}[1]{\ref{ex:#1}}
\newcommand{\propstart}{\begin{description}}
\newcommand{\prop}{\item}
\newcommand{\propend}{\end{description}}

\newcommand{\given}{\mbox{\ given\ }}
\newcommand{\goto}{\rightarrow}
\newcommand{\union}{\cup}
\newcommand{\gapprox}{\stackrel{>}{\approx}}
\newcommand{\recip}[1]{\frac{1}{#1}}
\newcommand{\binom}[2]{{#1 \choose #2}}
\newcommand{\multinom}{\binom}
\newcommand{\braces}[1]{\left ( #1 \right )}
\newcommand{\floor}[1]{\left \lfloor #1 \right \rfloor}
\newcommand{\R}{{\mathbf{R}}}
\newcommand{\Z}{{\mathbf{Z}}}
\newcommand{\Prob}{{\mathbf{P}}}
\newcommand{\Expect}{{\mathbf{E}}}

\newcommand{\Gbar}{\overline{G}}

\newcommand{\half}{{\frac{1}{2}}}

\newcommand{\fourth}{{\frac{1}{4}}}
\newcommand{\Rinfty}{{R_\infty}}
\newcommand{\eff}{{\mbox{eff}}}
\newcommand{\esc}{{\mbox{esc}}}
\newcommand{\eps}[1]{#1^{(\epsilon)}}
\newcommand{\Peps}{\eps{P}}
\newcommand{\peps}{\eps{p}}
\newcommand{\Reff}{R_\eff}
\newcommand{\Ceff}{C_\eff}
\newcommand{\pesc}{p_\esc}
\newcommand{\NTtwo}{{\mbox{NT}_2}}
\newcommand{\NTthree}{{\mbox{NT}_3}}
\newcommand{\NTsemi}{{\mbox{NT}_{2.5849\ldots}}}
\newcommand{\pbar}{\bar{p}}
\newcommand{\Pbar}{\bar{P}}
\newcommand{\Cbar}{\bar{C}}
\newcommand{\mat}[1]{{\mathbf{#1}}}
\newcommand{\matP}{\mat{P}}
\newcommand{\matPbar}{\mat{\Pbar}}
\newcommand{\orig}{\mat{0}}
\newcommand{\ball}[2]{#1^{(#2)}}
\newcommand{\sphere}[2]{\partial \ball{#1}{#2}}
\newcommand{\fuzz}[2]{#1_#2}
\newcommand{\matrixfive}[5]{\pmatrix{#1\cr#2\cr#3\cr#4\cr#5}}
\newdimen\snellbaselineskip
\newdimen\snellskip
\snellbaselineskip=\baselineskip
\def\srule{\omit\kern.5em\vrule\kern-.5em}
\newbox\bigstrutbox
\setbox\bigstrutbox=\hbox{\vrule height14.5pt depth9.5pt width0pt}
\def\bigstrut{\relax\ifmmode\copy\bigstrutbox\else\unhcopy\bigstrutbox\fi}
\def\middlehrule#1#2{\noalign{\kern-\snellbaselineskip\kern\snellskip}
&\multispan#1\strut\hrulefill
&\omit\hbox to.5em{\hrulefill}\vrule 
height \snellskip\kern-.5em&\multispan#2\hrulefill\cr}

\makeatletter
\def\bordermatrix#1{\begingroup \m@th
  \@tempdima 8.75\p@
  \setbox\z@\vbox{%
    \def\cr{\crcr\noalign{\kern2\p@\global\let\cr\endline}}%
    \ialign{$##$\hfil\kern2\p@\kern\@tempdima&\thinspace\hfil$##$\hfil
      &&\quad\hfil$##$\hfil\crcr
      \omit\strut\hfil\crcr\noalign{\kern-\snellbaselineskip}%
      #1\crcr\omit\strut\cr}}%
  \setbox\tw@\vbox{\unvcopy\z@\global\setbox\@ne\lastbox}%
  \setbox\tw@\hbox{\unhbox\@ne\unskip\global\setbox\@ne\lastbox}%
  \setbox\tw@\hbox{$\kern\wd\@ne\kern-\@tempdima\left(\kern-\wd\@ne
    \global\setbox\@ne\vbox{\box\@ne\kern2\p@}%
    \vcenter{\kern-\ht\@ne\unvbox\z@\kern-\snellbaselineskip}\,\right)$}%
  \null\;\vbox{\kern\ht\@ne\box\tw@}\endgroup}

\makeatother
\begin{document}
\maketitle
\section*{Preface}
Probability theory, like much of mathematics, is indebted to 
physics as a source of problems and intuition for solving these problems. 
Unfortunately, the level of abstraction of current mathematics often makes
it difficult for anyone but an expert to appreciate this fact.
In this \book\ we will look at the interplay of physics
and mathematics in terms of an example where the mathematics 
involved is at the college level.
The example is the relation between elementary electric network theory
and random walks.

Central to the \book\ will be Polya's beautiful theorem that a 
random walker on an infinite street network in $d$-dimensional space
is bound to return to the starting point when $d=2$,
but has a positive probability of escaping to infinity
without returning to the starting point when $d \geq 3$.
Our goal will be to interpret this theorem as a statement
about electric networks,
and then to prove the theorem
using techniques from classical electrical theory.
The techniques referred to go back to Lord Rayleigh,
who introduced them in connection with
an investigation of musical instruments.
The analog of Polya's theorem in this connection
is that wind instruments are possible in our three-dimensional world,
but are not possible in Flatland (Abbott \cite{abbott}).

The connection between random walks and electric networks has been
recognized for some time
(see
Kakutani \cite{kakutani},
Kemeny, Snell, and Knapp \cite{kemenysnellknapp},
and Kelly \cite{kelly}).
As for Rayleigh's method,
the authors first learned it from Peter's father Bill Doyle,
who used it to explain a mysterious comment in
Feller (\cite{feller1}, p. 425, Problem 14).
This comment suggested that a random walk in two dimensions
remains recurrent when some of the streets are blocked,
and while this is ticklish to prove probabilistically, 
it is an easy consequence of Rayleigh's method.
The first person to apply Rayleigh's method to random walks
seems to have been
Nash-Williams \cite{nashwilliams}.
Earlier, 
Royden \cite{royden}
had applied Rayleigh's method to an equivalent problem.
However, the true importance of Rayleigh's method for probability theory
is only now becoming appreciated.
See, for example,
Griffeath and Liggett \cite{griffeathliggett},
Lyons \cite{lyons},
and Kesten \cite{kesten}.

Here's the plan of the \book:
In \partref{I}\ we will restrict ourselves to the study of
random walks on finite networks.
Here we will establish the connection between
the electrical concepts of current and voltage
and corresponding descriptive quantities of random walks
regarded as finite state Markov chains.
In \partref{II} we will consider random walks on infinite networks.
Polya's theorem will be proved using Rayleigh's method,
and the proof will be compared with the classical proof
using probabilistic methods.
We will then discuss walks on more general infinite graphs,
and use Rayleigh's method to derive certain extensions of
Polya's theorem.
Certain of the results in \partref{II} were obtained
by Peter Doyle in work on his Ph.D. thesis.

To read this \book,
you should have a knowledge of the basic concepts of probability theory
as well as a little electric network theory and linear algebra.
An elementary introduction to finite Markov chains as presented by
Kemeny, Snell, and Thompson \cite{kemenysnellthompson}
would be helpful.

The work of Snell was carried out while enjoying the hospitality 
of Churchill College and the Cambridge Statistical Laboratory supported
by an NSF Faculty Development Fellowship.
He thanks Professors Kendall and Whittle
for making this such an enjoyable and rewarding visit.
Peter Doyle thanks his father
for teaching him how to think like a physicist.
We both thank Peter Ney for assigning
the problem in Feller that started all this,
David Griffeath for suggesting the example
to be used in our first proof that 3-dimensional random walk is recurrent
(\sectref{6.9}),
and Reese Prosser for keeping us going by his friendly and helpful hectoring.
Special thanks are due Marie Slack,
our secretary extraordinaire,
for typing the original and the excessive number of revisions one is led to
by computer formatting.

\prt{I}{Random walks on finite networks}

\chap{1}{Random walks in one dimension}

\sect{1.1}{A random walk along Madison Avenue}

A \dt{random walk}, or \dt{drunkard's walk},
was one of the first chance processes studied in probability;
this chance process continues to play an important role
in probability theory and its applications.
An example of a random walk may be described as follows:

A man walks along a 5-block stretch of Madison Avenue.
He starts at corner $x$ and,
with probability 1/2, walks one block to the right and,
with probability 1/2, walks one block to the left;
when he comes to the next corner
he again randomly chooses his direction along Madison Avenue.
He continues until he reaches corner 5, which is home,
or corner 0, which is a bar.
If he reaches either home or the bar, he stays there.
(See Figure \figref{1.1}.)
\prefig{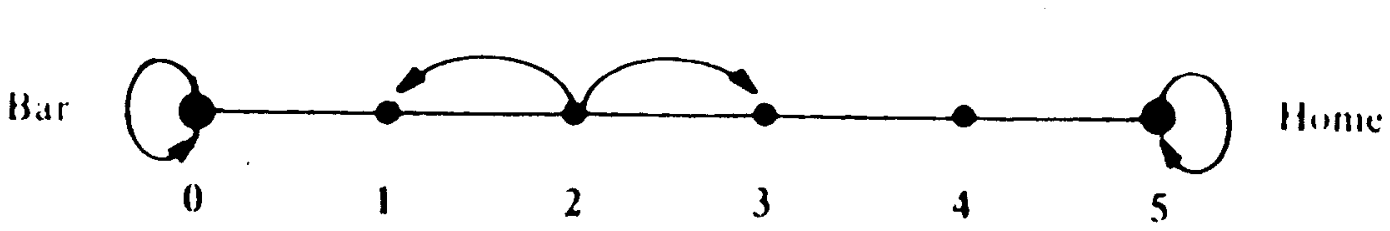}

The problem we pose is to find the probability $p(x)$ that the man,
starting at corner $x$, will reach home
before reaching the bar.
In looking at this problem, we will not be so much
concerned with the particular form of the solution,
which turns out to be $p(x)=x/5$, as with its general properties,
which we will eventually describe by saying
``$p(x)$ is the unique solution to a certain Dirichlet problem.''

\sect{1.2}{The same problem as a penny matching game}

In another form, the problem is posed in terms of the following game:
Peter and Paul match pennies; they have a total of 5 pennies;
on each match, Peter wins one penny from Paul with probability 1/2 and 
loses one with probability 1/2;
they play until Peter's fortune reaches 0 (he is ruined)
or reaches 5 (he wins all Paul's money).
Now $p(x)$ is the probability that
Peter wins if he starts with $x$ pennies.

\sect{1.3}{The probability of winning: basic properties}

Consider a random walk on the integers
$0, 1, 2, \ldots, N$.
Let $p(x)$ be the probability, starting at $x$,
of reaching $N$ before 0.
We regard $p(x)$ as a function defined on the points
$x=0, 1, 2, \ldots,N$.
The function $p(x)$ has the following properties:
\propstart
\prop{(a)} $p(0)=0$.
\prop{(b)} $p(N)=1$.
\prop{(c)} $p(x)=\half p(x - 1) + \half p(x + 1)$ for $x=1,2,\ldots,N- 1$.
\propend

Properties (a) and (b) follow from our convention that 0 and N are traps;
if the walker reaches one of these positions, he stops there;
in the game interpretation,
the game ends when one player has all of the pennies.
Property (c) states that, for an interior point, the 
probability $p(x)$ of reaching home from $x$
is the average of the probabilities $p(x - 1)$ and $p(x + 1)$
of reaching home from the points that the walker may go to from $x$.
We can derive (c) from the following basic fact about probability:

\tout{Basic Fact.}
Let $E$ be any event, and $F$ and $G$ be events
such that one and only one of the events $F$ or $G$ will occur.
Then
\[
\Prob(E) = \Prob(F) \cdot \Prob(E \given F) + \Prob(G)  \cdot \Prob(E \given G)
.
\]
\toutend

In this case, let $E$ be the event
``the walker ends at the bar'',
$F$ the event
``the first step is to the left'',
and $G$ the event
``the first step is to the right''.
Then, if the walker starts at $x$,
$\Prob(E) = p(x)$,
$\Prob(F) = \Prob(G)  = \half$,
$\Prob(E \given F) = p(x - 1)$,
$\Prob(E \given G) = p(x + 1)$,
and (c) follows.

\sect{1.4}{An electric network problem: the same problem?}

Let's consider a second apparently very different problem.
We connect equal resistors in series and put a unit voltage
across the ends as in Figure \figref{1.2}.
\prefig{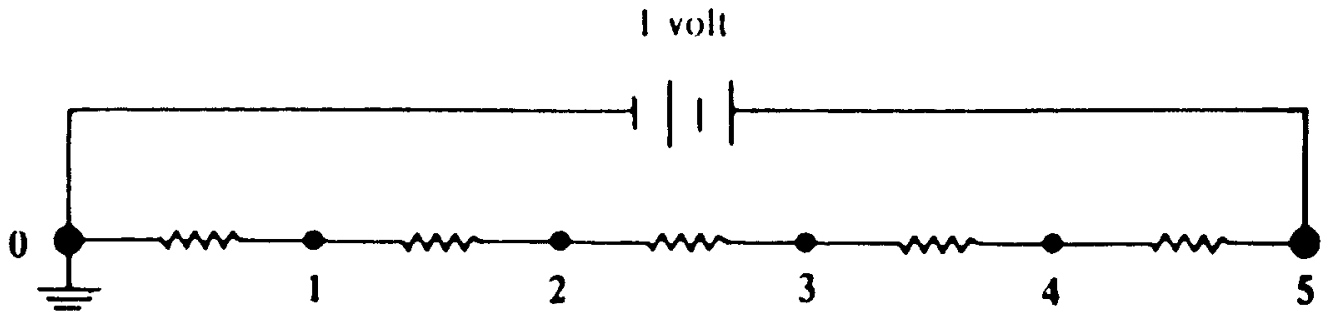}

Voltages $v(x)$ will be established at the points
$x = 0,1,2,3,4,5$.
We have grounded the point $x=0$ so that $v(0) = 0$.
We ask for the voltage $v(x)$
at the points $x$ between the resistors.
If we have $N$ resistors,
we make $v(0) = 0$ and $v(N) = 1$,
so $v(x)$ satisfies properties (a) and (b) of \sectref{1.3}.
We now show that $v(x)$ also satisfies (c).

By Kirchhoff's Laws, the current flowing into $x$ must be equal to 
the current flowing out.
By Ohm's Law, if points $x$ and $y$ are connected by
a resistance of magnitude $R$,
then the current $i_{xy}$ that flows from $x$ to $y$ is equal to
\[
i_{xy} = \frac{v(x)-v(y)}{R}
.
\]
Thus for $x=1,2,\ldots,N-1$,
\[
\frac{v(x-1)-v(x)}{R} + \frac{v(x+1)-v(x)}{R} = 0
.
\]
Multiplying through by $R$ and solving for $v(x)$ gives
\[
v(x) = \frac{v(x+1)+v(x-1)}{2}
\]
for $x=1,2,\ldots,N-1$.
Therefore, $v(x)$ also satisfies property (c).

We have seen that $p(x)$ and $v(x)$ both satisfy properties (a), (b), 
and (c) of \sectref{1.3}.
This raises the question: are $p(x)$ and $v(x)$ equal?
For this simple example, we can easily find $v(x)$ using Ohm's Law,
find $p(x)$ using elementary probability,
and see that they are the same.
However,
we want to illustrate a principle that will work for very general circuits.
So instead
we shall prove that these two functions are the same by showing 
that there is only one function that satisfies these properties,
and we shall prove this
by a method that will apply to more general situations than
points connected together in a straight line.

\exstart
\ex{1.4.1}{
Referring to the random walk along Madison Avenue,
let $X=p(1)$, $Y=p(2)$, $Z=p(3)$, and $W=p(4)$.
Show that properties (a), (b), and (c) of \sectref{1.3}\ determine
a set of four linear equations
with variables $X$, $Y$, $Z$ and $W$.
Show that these equations have a unique solution.
What does this say about $p(x)$ and $v(x)$ for this special case?
}

\ex{1.4.2}{
Assume that our walker has a tendency to drift in one direction:
more specifically, assume that each step is to the right with 
probability $p$ or to the left with probability $q=1-p$.
Show that properties (a), (b), and (c) of \sectref{1.3}\ should be
replaced by
\propstart
\prop{(a)} $p(0)=0$.
\prop{(b)} $p(N)=1$.
\prop{(c)} $p(x)= q \cdot p(x - 1) + p \cdot p(x + 1)$.
\propend
}

\ex{1.4.3}{
In our electric network problem,
assume that the resistors are not necessarily equal.
Let $R_x$ be the resistance between $x$ and $x + 1$. 
Show that
\[
v(x) =
\frac
{\frac{1}{R_{x-1}}}{\frac{1}{R_{x-1}}+\frac{1}{R_x}}
v(x-1)
+
\frac
{\frac{1}{R_x}}{\frac{1}{R_{x-1}}+\frac{1}{R_x}}
v(x+1)
.
\]
How should the resistors be chosen to correspond to the random 
walk of
Exercise \exref{1.4.2}?
}
\exend

\sect{1.5}{Harmonic functions in one dimension; the Uniqueness Principle}

Let $S$ be the set of points
$S = \{0,1,2,\ldots,N\}$.
We call the points of the set
$D = \{1,2,\ldots,N-1\}$ 
the \dt{interior points} of $S$
and those of
$B = \{0,N\}$
the \dt{boundary points} of $S$.
A function $f(x)$ defined on $S$ is
\dt{harmonic} if, 
at points of $D$,
it satisfies the averaging property
\[
f(x) = \frac{f(x - 1) + f(x + 1)}{2}
.
\]

As we have seen, $p(x)$ and $v(x)$ are harmonic functions on $S$ having 
the same values on the boundary:
$p(0) = v(0) = 0$; $p(N) = v(N) = 1$. 
Thus both $p(x)$ and $v(x)$ solve the problem of finding a harmonic function 
having these boundary values.
Now the problem of finding a harmonic function
given its boundary values is called
the \dt{Dirichlet problem}, and the
\dt{Uniqueness Principle} for the Dirichlet problem
asserts that there cannot be two different harmonic functions
having the same boundary values.
In particular, it follows that $p(x)$ and $v(x)$ are really the same function,
and this is what we have been hoping to show.
Thus the fact that $p(x) = v(x)$
is an aspect of a general fact about harmonic functions.

We will approach the Uniqueness Principle by way of
the \dt{Maximum Principle} for harmonic functions,
which bears the same relation to the
Uniqueness Principle as Rolle's Theorem does to the Mean Value
Theorem of Calculus.

\tout{Maximum Principle}.
A harmonic function $f(x)$ defined on $S$
takes on its maximum value $M$ and its minimum value $m$ on the boundary.
\toutend

\proofstart
Let $M$ be the largest value of $f$.
Then if $f(x) = M$ for $x$ in $D$,
the same must be true for $f(x - 1)$ and $f(x + 1)$
since $f(x)$ is the average of these two values.
If $x-1$  is still an interior point, the same
argument implies that $f(x - 2) = M$;
continuing in this way, we eventually conclude that $f(0) = M$.
That same argument works for the minimum value $m$.
\proofend

\tout{Uniqueness Principle.}
If $f(x)$ and $g(x)$ are harmonic functions on $S$
such that $f(x) = g(x)$ on $B$, then $f(x) = g(x)$ for all $x$.
\toutend

\proofstart
Let $h(x) = f(x) - g(x)$.
Then if $x$ is any interior point,
\[
\frac{h(x - 1) + h(x + 1)}{2}
=
\frac{f(x - 1) + f(x + 1)}{2}
-
\frac{g(x - 1) + g(x + 1)}{2}
,
\]
and $h$ is harmonic.
But $h(x) = 0$ for $x$ in $B$,
and hence, by the Maximum Principle,
the maximum and mininium values of $h$ are 0.
Thus $h(x)=0$ for all $x$,
and $f(x) = g(x)$ for all $x$.
\proofend

Thus we finally prove that $p(x)=v(x)$; but what does $v(x)$ equal?
The Uniqueness Principle shows us a way to find a concrete
answer: just guess.
For if we can find any harmonic function
$f(x)$ having the right boundary values,
the Uniqueness Principle guarantees that
\[
p(x) = v(x) = f(x)
.
\]

The simplest function to try for $f(x)$ would be a linear function; 
this leads to the solution $f(x) = x/N$.
Note that $f(0) = 0$ and $f(N) = 1$ and
\[
\frac{f(x-1) + f(x+1)}{2}
=
\frac{x-1+x+1}{2N}
=
\frac{x}{N}
=
f(x)
.
\]
Therefore $f(x) = p(x) = v(x) = x/N$.

As another application of the Uniqueness Principle,
we prove that our walker will eventually reach 0 or $N$.
Choose a starting point $x$ with $0 < x < N$.
Let $h(x)$ be the probability that the walker never reaches the 
boundary $B = \{0,N\}$.
Then
\[
h(x) =
\half h(x+1) + \half h(x-1)
\]
and h is harmonic.
Also $h(0) = h(N) = 0$;
thus, by the Maximum Principle,
$h(x) = 0$ for all $x$.

\exstart
\ex{1.5.1}{
 Show that you can choose $A$ and $B$
so that the function
$f(x) = A(q/p)^x + B$
satisfies the modified properties (a), (b) and (c) of
Exercise \exref{1.4.2}.
Does this show that $f(x) = p(x)$?
}

\ex{1.5.2}{
Let $m(x)$ be the expected number of steps, starting at $x$, 
required to reach 0 or $N$ for the first time.
It can be proven that $m(x)$ is finite.
Show that $m(x)$ satisfies the conditions
\propstart
\prop{(a)} $m(0)=0$.
\prop{(b)} $m(N)=0$.
\prop{(c)} $m(x)=\half m(x + 1) + \half m(x - 1) + 1$.
\propend
}

\ex{1.5.3}{
Show that the conditions in
Exercise \exref{1.5.2}\ have a unique solution. 
Hint: show
that if $m$ and $\bar{m}$ are two solutions,
then $f = m - \bar{m}$ is harmonic 
with $f (0) = f (N) = 0$
and hence $f(x) = 0$ for all $x$.
}

\ex{1.5.4}{
Show that you can choose $A$, $B$, and $C$
such that $f(x) = A + B x + C x^2$
satisfies all the conditions of
Exercise \exref{1.5.2}.
Does this show that $f(x) = m(x)$ for
this choice of $A$, $B$, and $C$?
}

\ex{1.5.5}{
Find the expected duration of the walk down Madison Avenue as a
function of the walker's starting point (1, 2, 3, or 4).
}
\exend

\sect{1.6}{The solution as a fair game (martingale)}

Let us return to our interpretation of a random walk as Peter's fortune
in a game of penny matching with Paul.
On each match,
Peter wins one penny with probability 1/2
and loses one penny with probability 1/2. 
Thus, when Peter has $k$ pennies his expected fortune after the next play
is
\[
\half (k-1) + \half (k+1) = k
,
\]
so his expected fortune after the next play
is equal to his present fortune.
This says that he is playing a
\dt{fair game};
a chance process that can be interpreted as a player's fortune
in a fair game is called a 
\dt{martingale}.

Now assume that Peter and Paul have a total of $N$ pennies.
Let $p(x)$ be the probability that,
when Peter has $x$ pennies,
he will end up with all $N$ pennies.
Then Peter's expected final fortune in this game is
\[
(1 - p(x)) \cdot 0 + p(x) \cdot N = p(x) \cdot N
.
\]

If we could be sure that a fair game remains fair to the end of the game,
then we could conclude that Peter's expected final fortune is 
equal to his starting fortune $x$, i.e., $x = p(x) \cdot N$.
This would give $p(x) = x/N$
and we would have found the probability that Peter wins
using the fact that a fair game remains fair to the end.
Note that the time the game ends is a random time, namely,
the time that the walk first reaches 0 or $N$ for the first time.
Thus the question is, is the fairness of
a game preserved when we stop at a random time?

Unfortunately, this is not always the case.
To begin with, if Peter somehow has knowledge of what the future holds
in store for him, 
he can decide to quit when he gets to the end of a winning streak.
But even if we restrict ourselves to stopping rules
where the decision to stop or continue
is independent of future events, fairness may not be preserved. 
For example, assume that Peter is allowed to go into debt and can 
play as long as he wants to.
He starts with 0 pennies and decides to play
until his fortune is 1 and then quit.
We shall see that a random walk on the set of all integers,
starting at 0,
will reach the point 1 if we wait long enough.
Hence, Peter will end up one penny ahead by this system of stopping.

However, there are certain conditions under which we can guarantee
that a fair game remains fair when stopped at a random time.
For our purposes,
the following standard result of martingale theory will do:

\tout{Martingale Stopping Theorem.}
A fair game that is stopped at a random time
will remain fair to the end of the game
if it is assumed that there is a finite amount of money in the world
and a player must stop if he wins all this money
or goes into debt by this amount.
\toutend

This theorem would justify the above argument to obtain
$p(x) = x/N$.

Let's step back and see how this martingale argument worked.
We began with a harmonic function,
the function $f(x) = x$,
and interpreted it as the player's fortune in a fair game.
We then considered the player's expected final fortune in this game.
This was another harmonic function
having the same boundary values
and we appealed to the Martingale Stopping Theorem
to argue that this function must be the same as the original function.
This allowed us to write down an expression 
for the probability of winning,
which was what we were looking for.

Lurking behind this argument is a general principle:
If we are given boundary values of a function,
we can come up with a harmonic function
having these boundary values by assigning to each point the 
player's expected final fortune
in a game where the player starts from the given point
and carries out a random walk until he reaches a boundary point,
where he receives the specified payoff.
Furthermore,
the Martingale Stopping Theorern allows us to conclude
that there can be no other harmonic function with these boundary values.
Thus martingale theory allows us to establish existence and uniqueness
of solutions to a Dirichlet problem.
All this isn't very exciting for the cases we've been considering,
but the nice thing is that the same arguments carry through
to the more general situations that we will be considering later on.

The study of martingales was originated by
Levy \cite{levy} and
Ville \cite{ville}.
Kakutani \cite{kakutani}
showed the connection between random walks and harmonic functions.
Doob \cite{doob}
developed martingale stopping theorems
and showed how to exploit the preservation of fairness to solve a 
wide variety of problems in probability theory.
An informal discussion of martingales may be found in
Snell \cite{snell}.

\exstart
\ex{1.6.1}{
Consider a random walk with a drift;
that is, there is a probability $p \neq \half$ of
going one step to the right and a probability $q = 1 - p$ of going 
one step to the left.
(See Exercise \exref{1.4.2}.)
Let $w(x) = (q/p)^x$;
show that,
if you interpret $w(x)$ as your fortune when you are at $x$,
the resulting game is fair.
Then use the Martingale Stopping Theorem to argue that
\[
w(x) = p(x) w(N) + (1 - p(x)) w(0)
.
\]
Solve for $p(x)$ to obtain
\[
p(x)
=
\frac
{\braces{\frac{q}{p}}^x -1}
{\braces{\frac{q}{p}}^N -1}
.
\]
}

\ex{1.6.2}{
You are gambling against a professional gambler;
you start with $A$ dollars and the gambler with $B$ dollars;
you play a game in which
you win one dollar with probability $p < \half$
and lose one dollar with probability $q = 1 - p$;
play continues until you or the gambler runs out of money.
Let $R_A$ be the probability that you are ruined.
Use the result of
Exercise \exref{1.6.1}\ to show that
\[
R_A
=
\frac
{1-\braces{\frac{p}{q}}^B}
{1-\braces{\frac{p}{q}}^N}
\]
with $N = A + B$.
If you start with 20 dollars and the gambler with 50 dollars
and $p = .45$,
find the probability of being ruined.
}

\ex{1.6.3}{
The gambler realizes that the probability of ruining you is at 
least $1 - (p/q)^B$
(Why?).
The gambler wants to make the probability at least .999. 
For this, $(p/q)^B$ should be at most .001.
If the gambler offers you a game with $p = .499$,
how large a stake should she have?
}
\exend

\chap{2}{Random walks in two dimensions}

\sect{2.1}{An example}

We turn now to the more complicated problem of a random walk on a 
two-dimensional array.
In Figure \figref{2.1} we illustrate such a walk.
\prefig{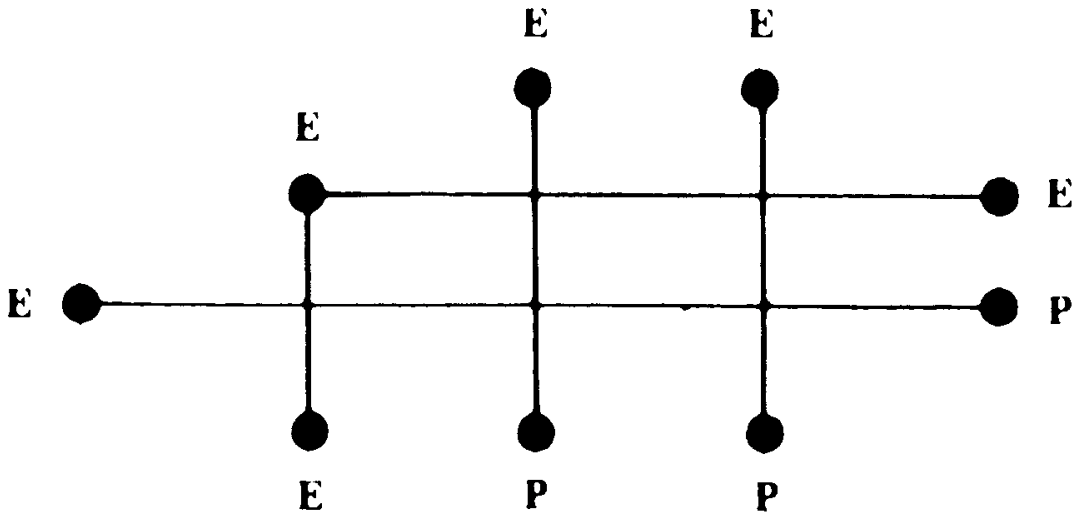}
The large dots represent boundary points;
those marked $E$ indicate escape routes and those 
marked $P$ are police.
We wish to find the probability $p(x)$ that our walker,
starting at an interior point $x$,
will reach an escape route before he reaches a policeman.
The walker moves from $x = (a, b)$ to each of the four
neighboring points
$(a + 1, b)$, $(a - 1, b)$, $(a, b + 1)$, $(a, b - 1)$
with probability $\fourth$.
If he reaches a boundary point, he remains at this point.

The corresponding voltage problem is shown in
Figure \figref{2.2}.
\prefig{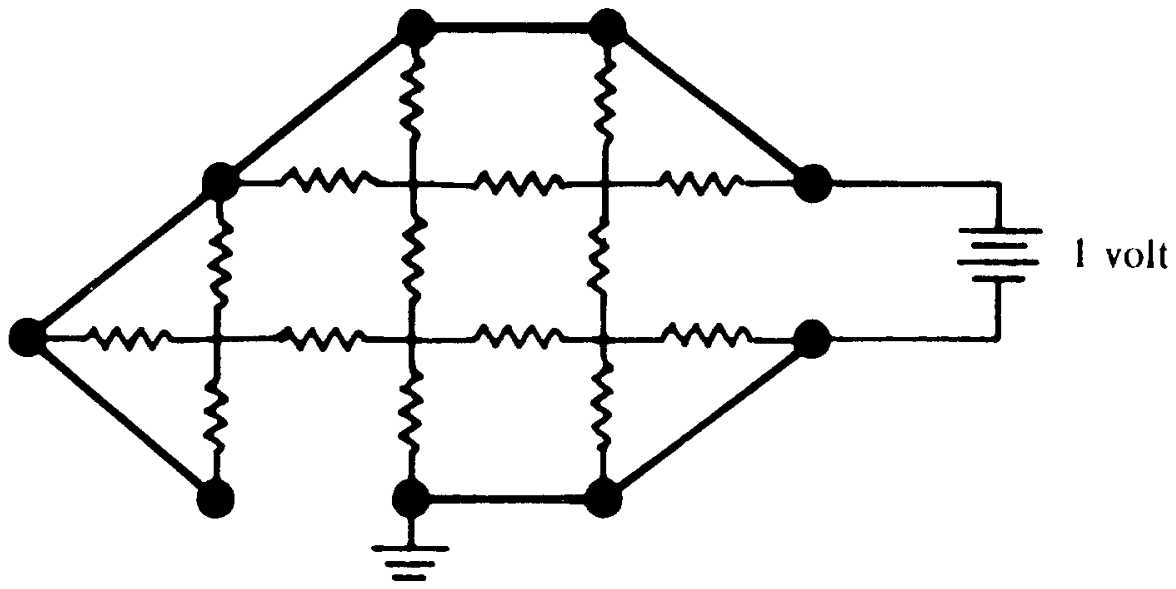}
The boundary points $P$ are grounded and points $E$ are connected and 
fixed at one volt by a one-volt battery.
We ask for the voltage $v(x)$ at the interior points.

\sect{2.2}{Harmonic functions in two dimensions}

We now define harmonic functions for sets of lattice points in 
the plane
(a lattice point is a point with integer coordinates).
Let $S = D \union B$
be a finite set of lattice points such that
(a) $D$ and $B$ have no points in common,
(b) every point of $D$ has its four neighboring points in $S$,
and
(c) every point of $B$
has at least one of its four neighboring points in $D$.
We assume further that $S$ hangs together in a nice way,
namely,
that for any two points $P$ and $Q$ in $S$,
there is a sequence of points $P_j$ in $D$ such that
$P, P_1,P_2, \ldots, P_n, Q$
forms a path from $P$ to $A$.
We call the points of $D$ the \dt{interior points}
of $S$ and the points of $B$ the \dt{boundary points} of $S$.

A function $f$ defined on $S$ is \dt{harmonic} if,
for points $(a, b)$ in $D$, 
it has the averaging property
\[
f (a, b) =
\frac{f(a + 1, b) +f(a - 1, b) +f(a, b + 1) +f(a, b - 1)}{4}
.
\]
Note that there is no restriction on the values of $f$
at the boundary points.

We would like to prove that $p(x) = v(x)$
as we did in the one-dimensional case.
That $p(x)$ is harmonic follows again by 
considering all four possible first steps;
that $v(x)$ is harmonic follows again by 
Kirchhoff's Laws
since the current coming into $x = (a, b)$ is
\[
\frac{v(a + 1, b) - v(a, b)}{R}
+
\frac{v(a - 1, b) - v(a, b)}{R}
+
\frac{v(a, b + 1) - v(a, b)}{R}
+
\frac{v(a, b - 1) - v(a, b)}{R}
=
0
.\]
Multiplying through by $R$ and solving for $v(a, b)$ gives
\[
v(a, b) =
\frac{v(a + 1, b) + v(a - 1, b) + v(a, b + 1) + v(a, b - 1)}{4}
.
\]

Thus $p(x)$ and $v(x)$ are harmonic functions with the same boundary
values.
To show from this that they are the same,
we must extend the Uniqueness Principle to two dimensions.

We first prove the Maximum Principle.
If $M$ is the maximum value of $f$
and if $f (P) = M$ for $P$ an interior point,
then since $f(P)$ is the average of the values of $f$ at its neighbors,
these values must all equal $M$ also. 
By working our way due south, say, repeating this argument at every 
step, we eventually reach a boundary point $Q$
for which we can conclude that $f(Q) = M$.
Thus a harmonic function always attains its maximum (or minimum)
on the boundary;
this is the Maximum Principle.
The proof of the Uniqueness Principle goes through as before
since again the difference of two harmonic functions is harmonic.

The fair game argument, using the Martingale Stopping Theorem, 
holds equally well and again gives an alternative proof of the 
existence and uniqueness to the solution of the Dirichlet problem.

\exstart
\ex{2.2.1}{
Show that if $f$ and $g$ are harmonic functions
so is $h = a \cdot f + b \cdot g$
for constants $a$ and $b$.
This is called the \dt{superposition principle}.
}

\ex{2.2.2}{
Let $B_1,B_2, \ldots, B_n$ be the boundary points for a region $S$.
Let $e_j(a, b)$ be a function that is harmonic in $S$
and has boundary value 1 at $B_j$
and 0 at the other boundary points.
Show that if arbitrary boundary values $v_1, v_2, \ldots, v_n$
are assigned,
we can find the harmonic function $v$ with these values
from the solutions $e_1,e_2,\ldots,e_n$.
}
\exend

\sect{2.3}{The Monte Carlo solution}

Finding the exact solution to a Dirichlet problem in two 
dimensions is not always a simple matter,
so before taking on this problem,
we will consider two methods for generating approximate solutions.
In this section we will present a method using random walks.
This method is known as a \dt{Monte Carlo method},
since random walks are random, and gambling involves
randomness, and there is a famous gambling casino in Monte Carlo.
In \sectref{2.4}, we will describe a much more
effective method for finding approximate solutions, called the
\dt{method of relaxations}.

We have seen that the solution to the Dirichlet problem can be
found by finding the value of a player's final winning in the
following game:
Starting at $x$ the player carries out a random walk
until reaching a boundary point.
He is then paid an amount $f(y)$
if $y$ is the boundary point first reached.
Thus to find $f(x)$,
we can start many random walks at $x$
and find the average final winnings for these walks.
By the law of averages (the law of large numbers in probability theory),
the estimate that we obtain this way will
approach the true expected final winning $f(x)$.

Here are some estimates obtained this way by starting 10,000
random walks from each of the interior points and, for each $x$,
estimating $f(x)$ by the average winning of the random walkers who
started at this point.
\[
\matrix{
&&1&1\cr
&1.824&.785&1\cr
1&.876&.503&.317&0\cr
&1&0&0
}
\]

This method is a colorful way to solve the problem, but quite inefficient.
We can use probability theory to estimate how inefficient it is.
We consider the case with boundary values I or 0 as in our example.
In this case, the expected final winning is just 
the probability that the walk ends up at a boundary point with value 1.
For each point $x$, assume that we carry out $n$ random walks;
we regard each random walk to be an experiment
and interpret the outcome of the $i$th experiment to be a ``success''
if the walker ends at a boundary point with a 1 and a ``failure'' otherwise.
Let $p =  p(x)$ be the unknown probability for success
for a walker starting at $x$ and $q = 1 - p$.
How many walks should we carry out to get a reasonable
estimate for $p$?
We estimate $p$ to be the fraction $\pbar$ of the walkers
that end at a 1.

We are in the position of a pollster who wishes to estimate the
proportion $p$ of people in the country
who favor candidate $A$ over $B$.
The pollster chooses a random sample of $n$ people
and estimates $p$ as the proportion $\pbar$ of voters in his sample
who favor $A$.
(This is a gross oversimplification of what a pollster does, of course.)
To estimate the number $n$ required,
we can use the central limit theorem.
This theorem states that,
if $S_n$, is the number of successes in $n$ independent experiments,
each having probability $p$ for success, then for any $k > 0$
\[
\Prob \braces{ -k < \frac{S_n-np}{\sqrt{npq}} < k }
\approx A(k),
\]
where $A (k)$ is the area under the normal curve between $-k$ and $k$.
For $k = 2$ this area is approximately .95;
what does this say about $\pbar = S_n/n$?
Doing a little rearranging, we see that
\[
\Prob \braces{ -2 < \frac{\pbar - p}{\sqrt{\frac{pq}{n}}} < 2 }
\approx
.95
\]
or
\[
\Prob \braces{ -2 \frac{\sqrt{pq}}{n} < \pbar - p < 2 \frac{\sqrt{pq}}{n} }
\approx .95
.
\]
Since $\sqrt{pq} \leq \half$,
\[
\Prob \braces{ - \frac{1}{\sqrt{n}} < \pbar - p < \frac{1}{\sqrt{n}} }
\gapprox .95
.
\]
Thus, if we choose $\frac{1}{\sqrt{n}} = .01$,
or $n =10,000$,
there is a 95 percent chance that our estimate $\pbar = S_n/n$
will not be off by more than .01.
This is a large number for rather modest accuracy;
in our example we carried out
10,000 walks from each point and this required about 5 seconds on 
the Dartmouth computer.
We shall see later, when we obtain an exact solution,
that we did obtain the accuracy predicted.

\exstart
\ex{2.3.1}{
You play a game in which you start a random walk at the center 
in the grid shown in
Figure \figref{2.3}.
\prefig{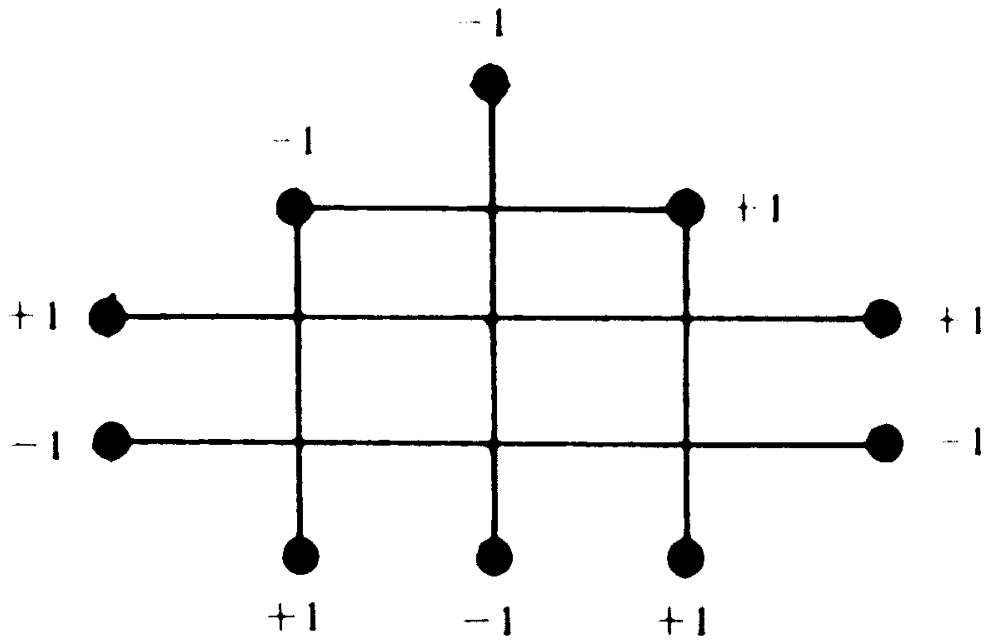}
When the walk reaches the boundary,
you receive a payment of $+1$ or $-1$
as indicated at the boundary points.
You wish to simulate this game
to see if it is a favorable game to play;
how many simulations would you need to be reasonably 
certain of the value of this game to an accuracy of .01?
Carry out such a simulation
and see if you feel that it is a favorable game.
}
\exend

\sect{2.4}{The original Dirichlet problem; the method of relaxations}

The Dirichlet problem we have been studying
is not the original Dirichlet problem,
but a discrete version of it.
The original Dirichlet problem
concerns the distribution of temperature, say,
in a continuous medium;
the following is a representative example.

Suppose we have a thin sheet of metal
gotten by cutting out a small square
from the center of a large square.
The inner boundary is kept at temperature 0
and the outer boundary is kept at temperature 1
as indicated in Figure \figref{2.4}.
\prefig{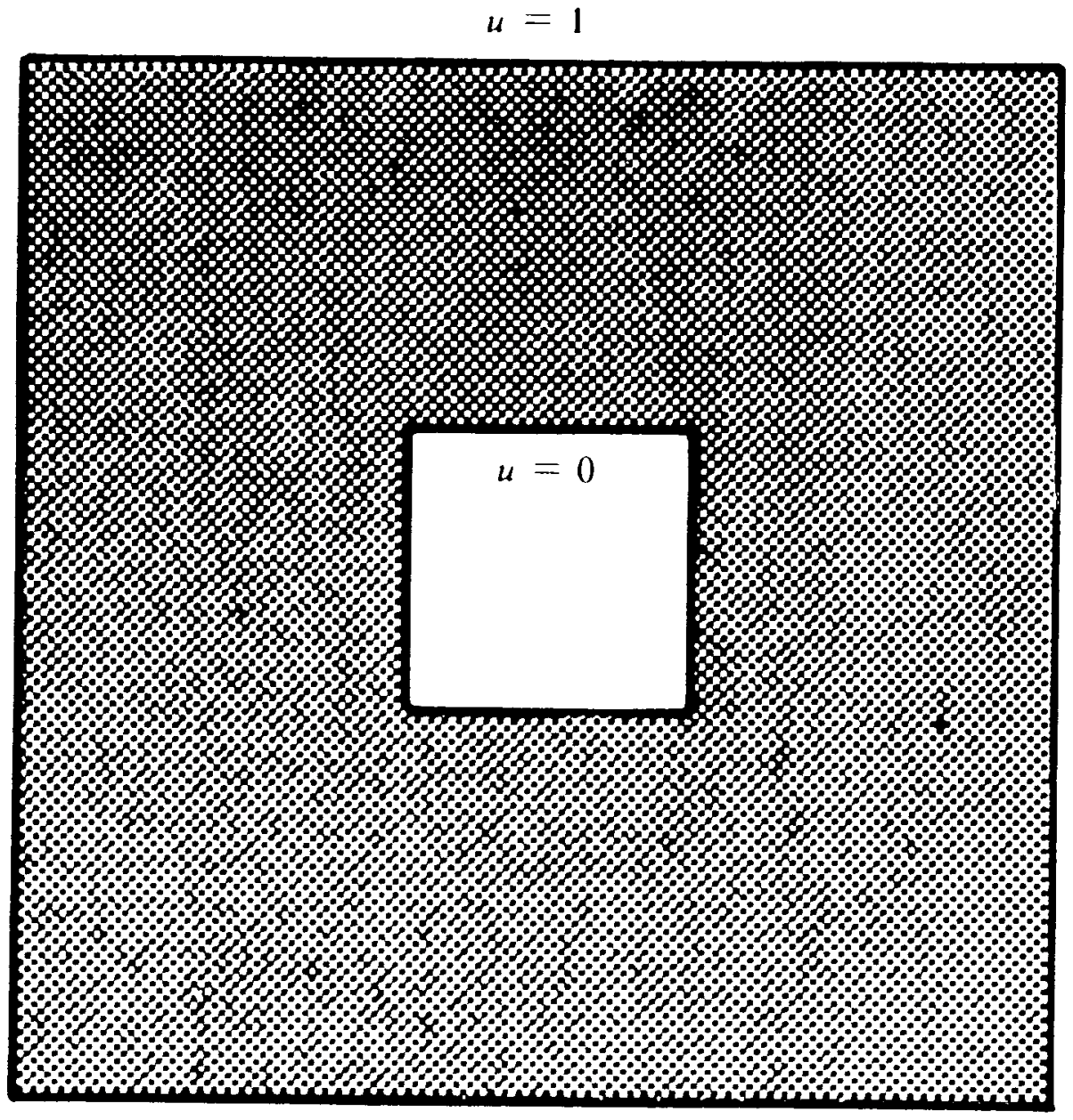}
The problem is to find the temperature at points 
in the rectangle's interior.
If $u(x, y)$ is the temperature at $(x, y)$, 
then $u$ satisfies Laplace's differential equation
\[
u_{xx} + u_{yy} = 0
.
\]
A function that satisfies this differential equation is called 
\dt{harmonic}.
It has the property that the value $u(x, y)$
is equal to the average of the values over
any circle with center $(x, y)$ lying inside the region.
Thus to determine the temperature $u(x, y)$,
we must find a harmonic function defined in 
the rectangle that takes on the prescribed boundary values.
We have a problem entirely analogous to our discrete Dirichlet problem,
but with continuous domain.

The \dt{method of relaxations} was
introduced as a way to get approximate solutions to the original
Dirichlet problem.
This method is actually more closely connected to
the discrete Dirichlet problem than to
the continuous problem.
Why? Because, faced with the continuous problem just described,
no physicist will hesitate to replace it
with an analogous discrete problem,
approximating the continuous medium 
by an array of lattice points such as that depicted in
Figure \figref{2.5},
and searching for a function that is harmonic
in our discrete sense and that takes on the appropriate
boundary values.
It is this approximating discrete problem
to which the method of relaxations applies.
\prefig{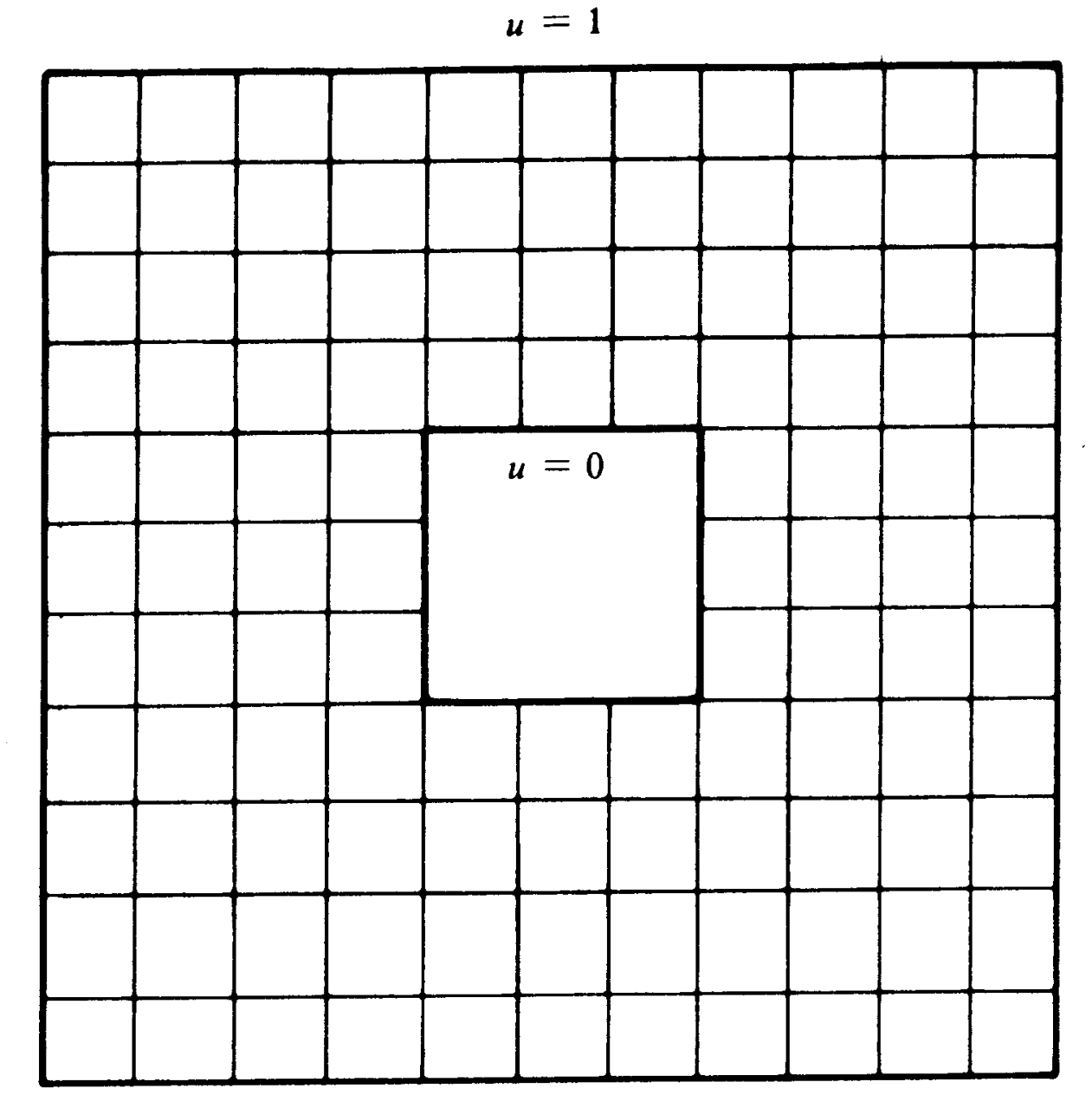}

Here's how the method goes.
Recall that we are looking for a function
that has specified boundary values,
for which the value at any interior point
is the average of the values at its neighbors.
Begin with any function
having the specified boundary values,
pick an interior point,
and see what is happening there.
In general,
the value of the function at the point we are looking at
will not be equal to
the average of the values at its neighbors.
So adjust the value of the function to
be equal to the average of the values at its neighbors.
Now run through the rest of the interior points,
repeating this process.
When you have adjusted the values at all of the interior points,
the function that results will not be harmonic,
because most of the time after adjusting the value at a 
point to be the average value at its neighbors,
we afterwards came along
and adjusted the values at one or more of those neighbors,
thus destroying the harmony.
However, the function that results
after running through all the interior points,
if not harmonic,
is more nearly harmonic than the function we started with;
if we keep repeating this averaging process, 
running through all of the interior points again and again,
the function will approximate more and more closely
the solution to our Dirichlet problem.

We do not yet have the tools to prove that this method works for a
general initial guess; this will have to wait until later
(see Exercise \exref{3.5.2}).
We will start with a special choice of initial values
for which we can prove that the method works
(see Exercise \exref{2.4.2}).

We start with all interior points 0 and keep the boundary points fixed.
\[
\matrix{
&&1&1\cr
&1&0&0&1\cr
1&0&0&0&0\cr
&1&0&0
}
\]
After one iteration we have:
\[
\matrix{
&&1&1\cr
&1&.547&.648&1\cr
1&.75&.188&.047&0\cr
&1&0&0
}
\]

Note that we go from left to right moving up each column
replacing each value by the average of the four neighboring values.
The computations for this first iteration are
\[
.75 = (1/4)(1 + 1 + 1 + 0)
\]
\[
.1875 = (1/4)(.75 + 0 + 0 + 0)
\]
\[
.5469 = (1/4)(.1875 + 1 + 1 + 0)
\]
\[
.0469 = (1/4)(.1875 + 0 + 0 + 0)
\]
\[
.64844 = (1/4)(.0469 + .5769 + 1 + 1)
\]
We have printed the results to three decimal places.
We continue the iteration
until we obtain the same results to three decimal places.
This occurs at iterations 8 and 9.
Here's what we get:
\[
\matrix{
&&1&1\cr
&1&.823&.787&1\cr
1&.876&.506&.323&0\cr
&1&0&0
}
\]

We see that we obtain the same result to three places
after only nine iterations
and this took only a fraction of a second of computing time.
We shall see that these results are correct to three place accuracy. 
Our Monte Carlo method took several seconds of computing time
and did not even give three place accuracy.

The classical reference for the method of relaxations as a means 
of finding approximate solutions to continuous problems is
Courant, Friedrichs, and Lewy \cite{courantfriedrichslewy}.
For more information on the relationship
between the original Dirichlet problem and the discrete analog, 
see Hersh and Griego \cite{hershgriego}.

\exstart
\ex{2.4.1}{
Apply the method of relaxations to the discrete problem 
illustrated in Figure \figref{2.5}.
}

\ex{2.4.2}{
Consider the method of relaxations started with an initial guess
with the property that the value at each point is $\leq$
the average of the values at the neighbors of this point.
Show that the successive values at a point $u$ are monotone
increasing with a limit $f(u)$
and that these limits provide a solution to the Dirichlet problem.
}
\exend

\sect{2.5}{Solution by solving linear equations}

In this section we will show how to find an exact solution
to a two-dimensional Dirichlet problem
by solving a system of linear equations.
As usual, we will illustrate the method in the case of the example 
introduced in \sectref{2.1}.
This example is shown again in Figure \figref{2.6};
the interior points have been labelled $a$, $b$, $c$, $d$, and $e$.
\prefig{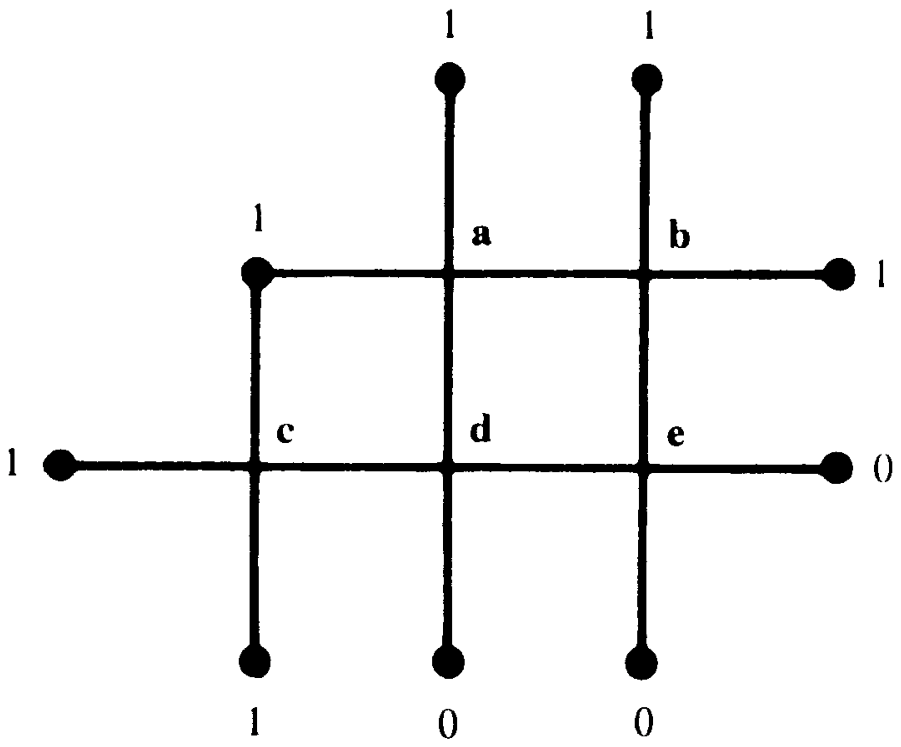}
By our averaging property, we have
\[
x_a = \frac{x_b + x_d + 2}{4}
\]
\[
x_b = \frac{x_a + x_c + 2}{4}
\]
\[
x_c = \frac{x_d + 3}{4}
\]
\[
x_d = \frac{x_a + x_c + x_e}{4}
\]
\[
x_e = \frac{x_b + x_d}{4}
.
\]
We can rewrite these equations in matrix form as
\[
\matrixfive
{1 & -1/4 & 0 & -1/4 & 0}
{ -1/4 & 1 & 0 & 0 & -1/4}
{ 0 & 0 & 1 & -1/4 & 0}
{ -1/4 & 0 & -1/4 & 1 & -1/4}
{ 0 & -1/4 & 0 & -1/4 & 1}
\matrixfive{x_a}{x_b}{x_c}{x_d}{x_e}
=
\matrixfive{1/2}{1/2}{3/4}{0}{0}
.
\]
We can write this in symbols as
\[
\mat{A} \mat{x} = \mat{u}
.
\]
Since we know there is a unique solution,
$\mat{A}$ must have an inverse and
\[
\mat{x} = \mat{A}^{-1} \mat{u}
.
\]
Carrying out this calculation we find
\[
\mbox{Calculated $\mat{x}$}
=
\matrixfive{.823}{.787}{.876}{.506}{.323}
.
\]

Here, for comparison, are the approximate solutions found earlier:
\[
\mbox{Monte Carlo $\mat{x}$}
=
\matrixfive{.824}{.785}{.876}{.503}{.317}
.
\]
\[
\mbox{Relaxed $\mat{x}$}
=
\matrixfive{.823}{.787}{.876}{.506}{.323}
.
\]
We see that our Monte Carlo approximations were fairly good
in that no error of the simulation is greater than .01,
and our relaxed approximations were very good indeed,
in that the error does not show up at all.

\exstart
\ex{2.5.1}{
Consider a random walker on the graph of Figure \figref{2.7}.
\prefig{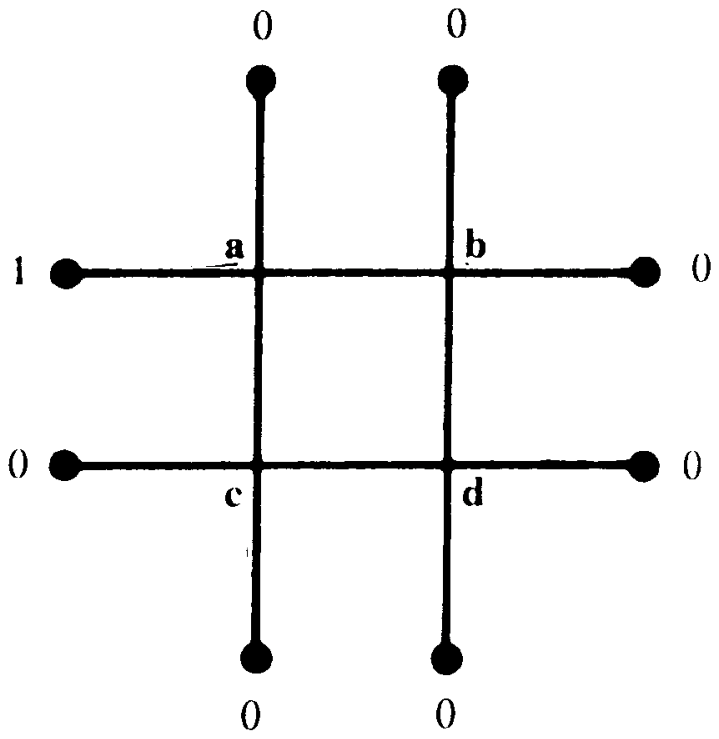}
Find the probability of reaching the point with a 1 before any of the 
points with 0's for each starting point $a, b, c, d$.
}

\ex{2.5.2}{
Solve the discrete Dirichlet problem for the graph shown in 
Figure \figref{2.8}.
\prefig{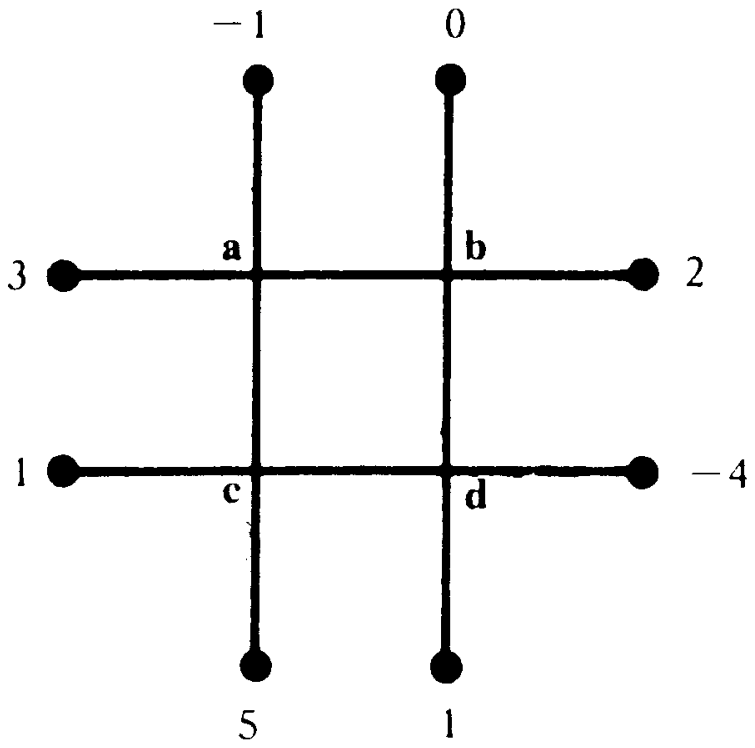}
The interior points are $a, b, c, d$.
(Hint: See Exercise \exref{2.2.2}.)
}

\ex{2.5.3}{
Find the exact value, for each possible starting point,
for the game described in Exercise \exref{2.3.1}.
Is the game favorable starting in the center?
}
\exend

\sect{2.6}{Solution by the method of Markov chains}

In this section,
we describe how the Dirichlet problem can be solved
by the method of Markov chains.
This method may be viewed as a more sophisticated version
of the method of linear equations.

A \dt{finite Markov chain}
is a special type of chance process that may be
described informally as follows:
we have a set $S = \{s_1, s_2, \ldots, s_r\}$ of \dt{states}
and a chance process that moves around through these states.
When the process is in state $s_i$,
it moves with probability $P_{ij}$ to the state $s_j$.
The transition probabilities $P_{ij}$ are represented by
an $r$-by-$r$ matrix $\matP$
called the \dt{transition matrix}.
To specify the chance process completely 
we must give, in addition to the transition matrix,
a method for starting the process.
We do this by specifying a specific state in which the process 
starts.

According to Kemeny, Snell, and Thompson \cite{kemenysnellthompson},
in the Land of Oz,
there are three kinds of weather: rain, nice, and snow.
There are never two nice days in a row.
When it rains or snows,
half the time it is the same the next day.
If the weather changes,
the chances are equal for a change
to each of the other two types of weather.
We regard the weather in the Land of Oz
as a Markov chain with transition matrix:
\[
\mat {P} = \bordermatrix{
        & \mbox {R} & \mbox {N} & \mbox {S} \cr
\mbox {R} &     1/2 &     1/4 &     1/4 \cr
\mbox {N} &     1/2 &       0 &     1/2 \cr
\mbox {S} &     1/4 &     1/4 &     1/2}
.
\]

When we start in a particular state,
it is natural to ask for the probability
that the process is in each of the possible states
after a specific number of steps.
In the study of Markov chains,
it is shown that this information is provided by
the powers of the transition matrix.
Specifically, if $\matP^n$
is the matrix $\matP$ raised to the $n$th power,
the entries $P^n_{ij}$ represent the probability
that the chain, started in state $s_i$, will, after $n$ steps,
be in state $s_j$.
For example, the fourth power of the transition matrix $\matP$ for the 
weather in the Land of Oz is
\[
\mat {P}^4  = \bordermatrix{
            &\mbox{R}&\mbox{N}&\mbox{S} \cr
\mbox{R} & .402      & .199      & .398 \cr
\mbox{N} & .398      & .203      & .398 \cr
\mbox{S} & .398      & .199      & .402 \cr}
.
\]

Thus, if it is raining today in the Land of Oz,
the probability that the weather will be nice four days from now is .199.
Note that the probability of a particular type of weather
four days from today
is essentially independent of the type of weather today.
This Markov chain is an example of a type of chain
called a regular chain.
A Markov chain is a \dt{regular} chain
if some power of the transition matrix has no zeros.
In the study of regular Markov chains,
it is shown that the probability of being in a state
after a large number of steps
is independent of the starting state.

As a second example, we consider a random walk in one dimension.
Let us assume that the walk is stopped when it reaches either 
state 0 or 4.
(We could use 5 instead of 4, as before, but we want to keep the 
matrices small.)
We can regard this random walk as a Markov chain
with states 0, 1, 2, 3, 4
and transition matrix given by
\[
\matP  = \bordermatrix{
  &0&1&2&3&4 \cr
0& 1&0&0&0&0 \cr
1& 1/2&0&1/2&0&0 \cr
2& 0&1/2&0&1/2&0 \cr
3& 0&0&1/2&0&1/2 \cr
4& 0&0&0&0&1}
.
\]

The states 0 and 4 are \dt{traps} or \dt{absorbing states}.
These are states that, once entered, cannot be left.
A Markov chain is called \dt{absorbing}
if it has at least one absorbing state
and if, from any state, it is possible (not necessarily in one step)
to reach at least one absorbing state. 
Our Markov chain has this property
and so is an absorbing Markov chain.
The states of an absorbing chain that are not traps
are called \dt{non-absorbing}.

When an absorbing Markov chain is started in a non-absorbing state,
it will eventually end up in an absorbing state.
For non-absorbing state $s_i$
and absorbing state $s_j$,
we denote by $B_{ij}$
the probability that the chain
starting in $s_i$ will end up in state $s_j$.
We denote by $\mat{B}$
the matrix with entries $B_{ij}$.
This matrix will have as many rows as non-absorbing states 
and as many columns as there are absorbing states.
For our random walk example,
the entries $B_{x,4}$
will give the probability that our random walker,
starting at $x$, will reach 4 before reaching 0.
Thus, if we can find the matrix $\mat{B}$ by Markov chain techniques,
we will have a way to solve the Dirichlet problem.

We shall show, in fact, that the Dirichlet problem has a natural
generalization in the context of absorbing Markov chains
and can be solved by Markov chain methods.

Assume now that $\matP$ is an absorbing Markov chain
and that there are $u$ absorbing states and $v$ non-absorbing states.
We reorder the states so that the absorbing states come first
and the non-absorbing states come last.
Then our transition matrix has the canonical form:
\[
\matP = \pmatrix{
\mat{I} & \mat{0} \cr
\mat{R} & \mat{Q} }
.
\]
Here $\mat{I}$ is a $u$-by-$u$ identity matrix;
$\mat0$ is a matrix of dimension $u$-by-$v$ with all entries 0.

For our random walk example this canonical form is:
\[
\bordermatrix{
 &0&4& 1&2&3 \cr
0& 1&0& 0&0&0 \cr
4& 0&1& 0&0&0 \cr
1& 1/2&0& 0&1/2&0 \cr
2& 0&0& 1/2&0&1/2 \cr
3& 0&1/2& 0&1/2&0 }
.
\]
The matrix $\mat{N} = (\mat{I} - \mat{Q})^{-1}$
is called the \dt{fundamental matrix}
for the absorbing chain $\matP$.
The entries $N_{ij}$ of this matrix
have the following probabilistic interpretation:
$N_{ij}$ is the expected number of times that the chain
will be in state $s_j$ before absorption
when it is started in $s_i$. 
(To see why this is true,
think of how $(\mat{I} - \mat{Q})^{-1}$ would look if it were written as a
geometric series.)
Let $\mat{1}$ be a column vector of all 1's.
Then the vector $\mat{t} = \mat{N} \mat{I}$
gives the expected number of steps before absorption for each
starting state.

The absorption probabilities $\mat{B}$
are obtained from $\mat{N}$ by the matrix formula
\[
\mat{B} = (\mat{I} - \mat{Q})^{-1} \mat{R}
.
\]
This simply says that to get
the probability of ending up at a given absorbing state,
we add up the probabilities
of going there from all the non-absorbing states,
weighted by the number of times
we expect to be in those (non-absorbing) states.

For our random walk example
\[
\mat{Q} = \pmatrix{
0&\half&0 \cr
\half&0&\half \cr
0&\half&0 \cr}
\]
\[
\mat{I}-\mat{Q} = \pmatrix{
1&-\half&0 \cr
-\half&1&-\half \cr
0&-\half&1}
\]
\[
\mat{N} = (\mat{I}-\mat{Q})^{-1} = \bordermatrix{
 &1&2&3 \cr
1 & \frac{3}{2}&1&\half \cr
2& 1&2&1 \cr
3& \half&1&\frac{3}{2} } 
\]
\[
\mat{t} = \mat{N} \mat{1} =
\pmatrix{
\frac{3}{2}&1&\half \cr
1&2&1 \cr
\half&1&\frac{3}{2} } 
\pmatrix{1 \cr 1 \cr 1}
=
\pmatrix{3 \cr 4 \cr 3}
\]
\[
\mat{B} = \mat{N} \mat{R} =
\pmatrix{
\frac{3}{2}&1&\half \cr
1&2&1 \cr
\half&1&\frac{3}{2} }
\pmatrix{
\half&0 \cr
0&0 \cr
0&\half}
=
\bordermatrix{
 &0&4 \cr
1& \frac{3}{4}&\frac{1}{4} \cr
2& \half&\half \cr
3& \frac{1}{4}&\frac{3}{4} }
.
\]

Thus, starting in state 3, the probability is $3/4$ of reaching 4 before 0;
this is in agreement with our previous results.
From $\mat{t}$ we see that the expected duration of the game,
when we start in state 2, is 4.

For an absorbing chain $\matP$,
the $n$th power $\matP^n$ of the transition probabilities
will approach a matrix $\matP^\infty$ of the form
\[
\matP^\infty = \pmatrix{
\mat{I} & \mat{0} \cr
\mat{B} & \mat{Q} }
.
\]

We now give our Markov chain version of the Dirichlet problem.
We interpret the absorbing states as boundary states and the 
non-absorbing states as interior states.
Let $B$ be the set of boundary states
and $D$ the set of interior states.
Let $f$ be a function with domain the state space of a Markov chain $\matP$
such that for $i$ in $D$
\[
f(i) = \sum_j P_{ij} f(j)
.
\]
Then $f$ is a \dt{harmonic function} for $\matP$.
Now $f$ again has an averaging property
and extends our previous definition.
If we represent $f$ as a column vector $\mat{f}$,
$f$ is harmonic if and only if
\[
\matP \mat{f} = \mat{f}
.
\]
This implies that
\[
\matP^2 \mat{f} = \matP \cdot \matP \mat{f} = \matP \mat{f} = \mat{f}
\]
and in general
\[
\matP^n \mat{f} = \mat{f}
.
\]

Let us write the vector $\mat{f}$ as
\[
\mat{f} =
\pmatrix{\mat{f}_B \cr \mat{f}_D}
\]
where $\mat{f}_B$ represents the values of $f$ on the boundary
and $\mat{f}_D$ values on the interior.
Then we have
\[
\pmatrix{\mat{f}_B \cr \mat{f}_D}
=
\pmatrix{
\mat{I} & \mat{0} \cr
\mat{B} & \mat{Q} }
\pmatrix{\mat{f}_B \cr \mat{f}_D} 
\]
and
\[
\mat{f}_D = \mat{B} \mat{f}_B
.
\]

We again see that the values of a harmonic function are determined by
the values of the function at the boundary points.

Since the entries $B_{ij}$ of $\mat{B}$
represent the probability, starting in $i$,
that the process ends at $j$,
our last equation states that if you play a game in which
your fortune is $f_j$ when you are in state $j$,
then your expected final fortune is equal to your initial fortune;
that is, fairness is preserved. 
As remarked above,
from Markov chain theory
$\mat{B} = \mat{N} \mat{R}$ where $\mat{N} = (\mat{I} - \mat{Q})^{-1}$.
Thus
\[
\mat{f}_D = (\mat{I} - \mat{Q})^{-1} \mat{R} \mat{f}_B
.
\]
(To make the correspondence between this solution and the
solution of \sectref{2.5}, put
$\mat{A}  = \mat{I} - \mat{Q}$ and $\mat{u} = \mat{R} \mat{f}_B$.)

A general discussion of absorbing Markov chains may be found
in Kemeny, Snell, and Thompson \cite{kemenysnellthompson}.

\exstart
\ex{2.6.1}{
Consider the game played on the grid in Figure \figref{2.9}.
\prefig{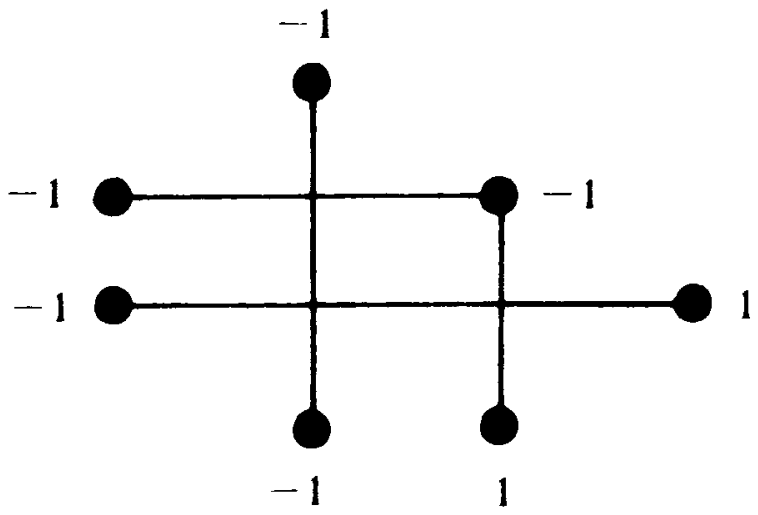}
You start at an interior point
and move randomly until a boundary point is reached
and obtain the payment indicated at this point.
Using Markov chain methods find, 
for each starting state, the expected value of the game.
Find also the expected duration of the game.
}
\exend

\chap{3}{Random walks on more general networks}

\sect{3.1}{General resistor networks and reversible Markov chains}

Our networks so far have been very special networks with unit resistors.
We will now introduce general resistor networks,
and consider what it means to carry out a random walk on such a network.

A \dt{graph}
is a finite collection of \dt{points}
(also called \dt{vertices} or \dt{nodes})
with certain pairs of points connected by \dt{edges}
(also called \dt{branches}).
The graph is \dt{connected}
if it is possible to go between any two points by moving along the edges.
(See Figure \figref{3.1}.)
\prefig{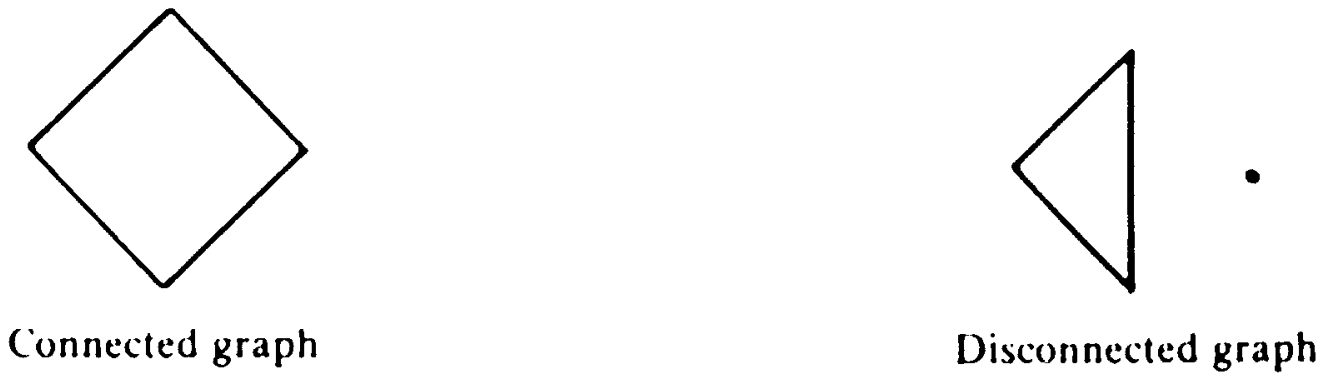}

We assume that $G$ is a connected graph and assign to each edge $xy$ a
resistance $R_{xy}$;
an example is shown in Figure \figref{3.2}.
\prefig{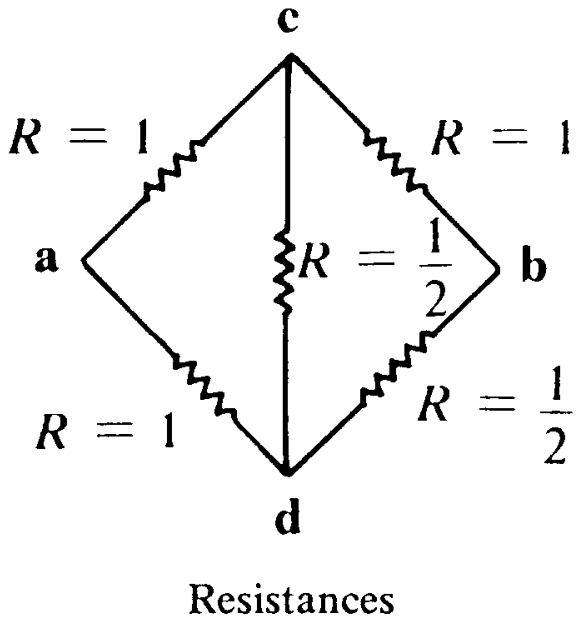}
The \dt{conductance} of an edge $xy$ is
$C_{xy} = 1/R_{xy}$;
conductances for our example are shown in Figure \figref{3.3}.
\prefig{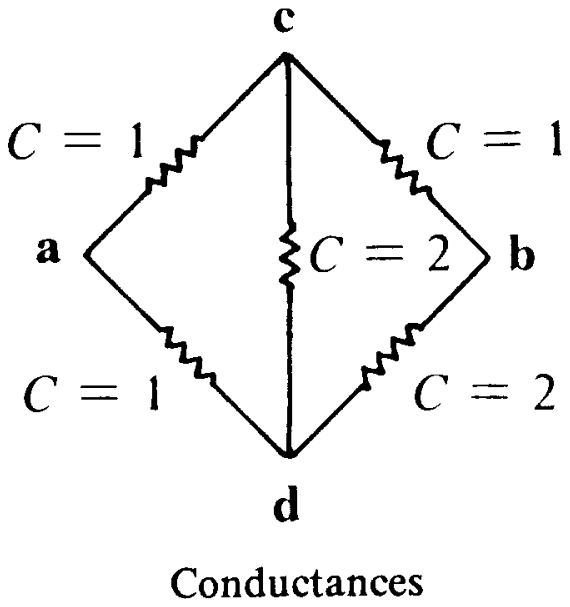}

We define a \dt{random walk} on $G$
to be a Markov chain with transition matrix $\matP$ given by
\[
P_{xy} = \frac{C_{xy}}{C_x}
\]
with $C_x = \sum_y C_{xy}$.
For our example,
$C_a = 2$, $C_b = 3$, $C_c = 4$, and $C_d = 5$,
and the transition matrix $\matP$ for the associated random walk is
\[
\bordermatrix{
&a&b&c&d \cr
a& 0&0&\half&\half \cr
b& 0&0&\frac{1}{3}&\frac{2}{3} \cr
c& \frac{1}{4}&\frac{1}{4}&0&\half \cr
d& \frac{1}{5}&\frac{2}{5}&\frac{2}{5}&0 }
\]
Its graphical representation is shown in Figure \figref{3.4}.
\prefig{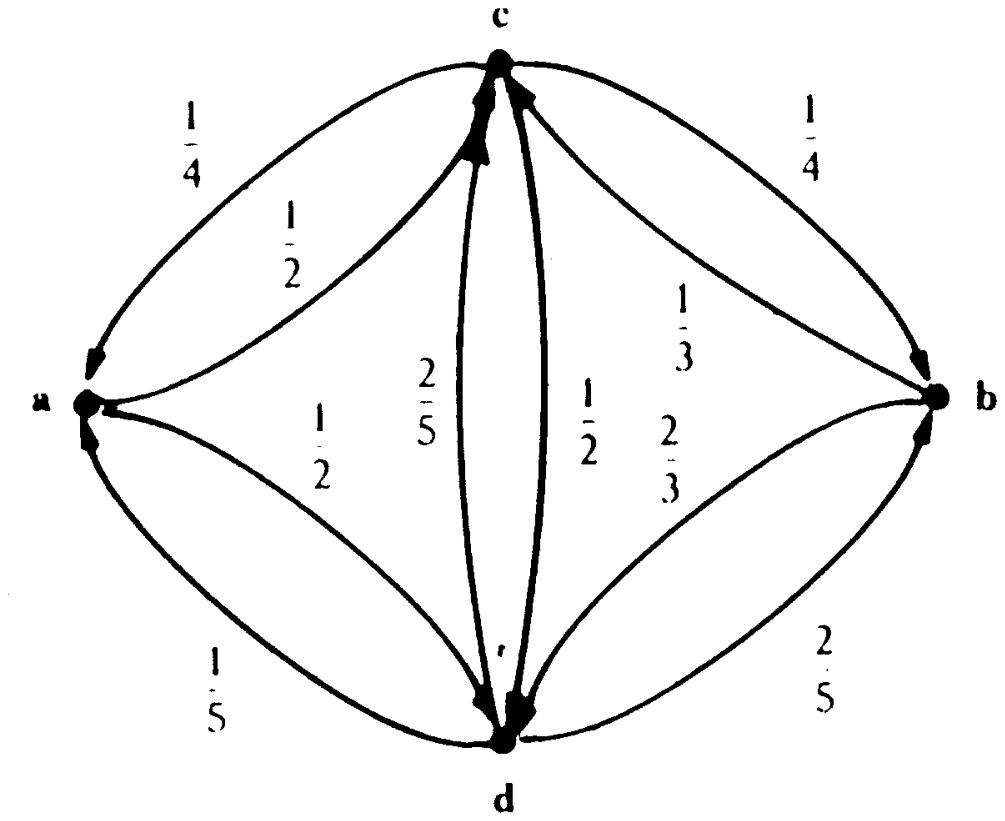}

Since the graph is connected,
it is possible for the walker to go between any two states.
A Markov chain with this property is called an
\dt{ergodic} Markov chain.
Regular chains, which were introduced in \sectref{2.6},
are always ergodic,
but ergodic chains are not always regular
(see Exercise \exref{3.1.1}).

For an ergodic chain, there is a unique probability vector $\mat{w}$ that
is a fixed vector for $\matP$, i.e., $\mat{w} \matP = \mat{w}$.
The component $w_j$ of $\mat{w}$
represents the proportion of times, in the long run,
that the walker will be in state $j$.
For random walks determined by electric networks,
the fixed vector is given by
$w_j = C_j/C$, where $C = \sum_x C_x$.
(You are asked to prove this in Exercise \exref{3.1.2}.)
For our example 
$C_a = 2$, $C_b = 3$, $C_c = 4$, $C_d = 5$, and $C=14$.
Thus
$\mat{w} = (2/14,3/14,4/14,5/14)$.
We can check that $\mat{w}$ is a fixed vector by noting that
\[
\pmatrix{
\frac{2}{14}&
\frac{3}{14}&
\frac{4}{14}&
\frac{5}{14}}
\pmatrix{
0&0&\half&\half \cr
0&0&\frac{1}{3}&\frac{2}{3} \cr
\frac{1}{4}&\frac{1}{4}&0&\half \cr
\frac{1}{5}&\frac{2}{5}&\frac{2}{5}&0 }
=
\pmatrix{
\frac{2}{14}&
\frac{3}{14}&
\frac{4}{14}&
\frac{5}{14}}
.
\]

In addition to being ergodic,
Markov chains associated with networks
have another property called \dt{reversibility}.
An ergodic chain is said to be \dt{reversible}
if $w_x P_{xy} = w_y P_{yx}$ for all $x, y$. 
That this is true for our network chains
follows from the fact that
\[
C_x P_{xy}
= C_x \frac{C_{xy}}{C_x}
= C_{xy}
= C_{yx}
= C_y \frac{C_{yx}}{C_y}
= C_y P_{yx}
.
\]
Thus, dividing the first and last term by $C$,
we have
$w_x P_{xy} = w_y P_{yx}$.

To see the meaning of reversibility,
we start our Markov chain with initial probabilities $\mat{w}$
(in equilibrium) and observe a few states,
for example
\[
a \;\; c \;\; b \;\; d
.
\]
The probability that this sequence occurs is
\[
w_a P_{ac} P_{cb} P_{bd}
=
\frac{2}{14} \cdot \frac{1}{2} \cdot \frac{1}{4} \cdot \frac{2}{3}
=
\frac{1}{84}
.
\]
The probability that the reversed sequence
\[
d \;\; b \;\; c \;\; a
\]
occurs is
\[
w_d P_{db} P_{bc} P_{ca}
=
\frac{5}{14} \cdot \frac{2}{5} \cdot \frac{1}{3} \cdot \frac{1}{4}
=
\frac{1}{84}
.
\]
Thus the two sequences have the same probability of occurring.

In general,
when a reversible Markov chain is started in equilibrium,
probabilities for sequences in the correct order of time
are the same as those with time reversed.
Thus, from data, we would never be able to tell the direction of time.

If $\matP$ is any reversible ergodic chain,
then $\matP$ is the transition matrix
for a random walk on an electric network;
we have only to define
$C_{xy} = w_x P_{xy}$.
Note, however,
if $P_{xx} \neq 0$ the resulting network will need a
conductance from $x$ to $x$ (see Exercise \exref{3.1.4}).
Thus reversibility characterizes
those ergodic chains that arise from electrical networks.
This has to do with the fact that
the physical laws that govern the behavior of steady electric currents
are invariant under time-reversal
(see Onsager \cite{onsager}).

When all the conductances of a network are equal,
the associated random walk on the graph $G$ of the network
has the property that,
from each point,
there is an equal probability of moving to each of the points
connected to this point by an edge.
We shall refer to this random walk as
\dt{simple random walk} on $G$.
Most of the examples we have considered so far are simple random walks.
Our first example of a random walk on Madison Avenue
corresponds to simple random walk on the graph with
points $0,1,2,\ldots,N$ and edges the streets connecting these points.
Our walks on two dimensional graphs were also simple random walks.

\exstart
\ex{3.1.1}{
Give an example of an ergodic Markov chain that is not regular.
(Hint: a chain with two states will do.)
}

\ex{3.1.2}{
Show that, if $\matP$ is the transition matrix for a random walk 
determined by an electric network,
then the fixed vector $\mat{w}$
is given by $w_x = \frac{C_x}{C}$
where $C_x = \sum_y C_{xy}$ and $C= \sum_x C_x$.
}

\ex{3.1.3}{
Show that, if $\matP$ is a reversible Markov chain
and $a, b, c$ are any three states,
then the probability, starting at $a$,
of the cycle $abca$ is the same as the probability
of the reversed cycle $acba$.
That is $P_{ab}P_{bc}P_{ca} = P_{ac}P_{cb}P_{ba}$.
Show, more generally, that the probability of going around any 
cycle in the two different directions is the same.
(Conversely, if this cyclic condition is satisfied,
the process is reversible.
For a proof, see Kelly \cite{kelly}.)
}

\ex{3.1.4}{
Assume that $\matP$ is a reversible Markov chain
with $P_{xx}=0$ for all $x$.
Define an electric network by
$C_{xy} = w_x P_{xy}$.
Show that the Markov chain associated with this circuit is $\matP$.
Show that we can allow $P_{xx} > 0$
by allowing a conductance from $x$ to $x$.
}

\ex{3.1.5}{
For the \dt{Ehrenfest urn model},
there are two urns that together contain $N$ balls.
Each second, one of the $N$ balls is chosen at random
and moved to the other urn.
We form a Markov chain with states the number of balls in
one of the urns.
For $N = 4$,
the resulting transition matrix is
\[
\matP = \bordermatrix{
 &0&1&2&3&4 \cr
0& 0&1&0&0&0 \cr
1& \frac{1}{4}&0&\frac{3}{4}&0&0 \cr
2& 0&\half&0&\half&0 \cr
3& 0&0&\frac{3}{4}&0&\frac{1}{4} \cr
4& 0&0&0&1&0 }
.
\]
Show that the fixed vector $\mat{w}$ is the binomial distribution
$\mat{w} =
(\frac{1}{16}, \frac{4}{16}, \frac{6}{16}, \frac{4}{16}, \frac{1}{16})$.
Determine the electric network associated with this chain.
}
\exend

\sect{3.2}{Voltages for general networks; probabilistic interpretation}

We assume that we have a network of resistors assigned to the edges of a
connected graph.
We choose two points $a$ and $b$
and put a one-volt battery
across these points establishing a voltage
$v_a = 1$ and $v_b = 0$,
as illustrated in Figure \figref{3.5}.
\prefig{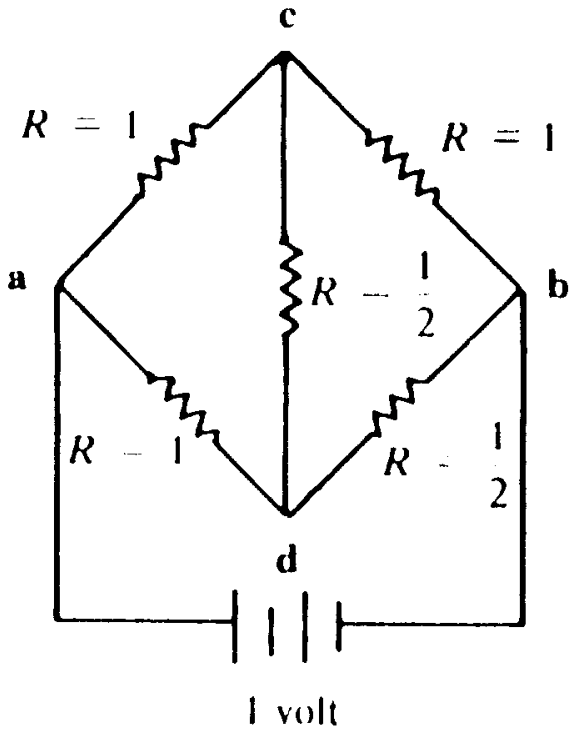}
We are interested in finding the voltages $v_x$
and the currents $i_{xy}$
in the circuit
and in giving a probabilistic interpretation to these quantities.

We begin with the probabilistic interpretation of voltage.
It will come as no surprise that we will interpret the voltage
as a hitting probability,
observing that both functions are harmonic and that they have the 
same boundary values.

By Ohm's Law, the currents through the resistors are determined 
by the voltages by
\[
i_{xy} = \frac{v_x-v_y}{R_{xy}}
=
(V_x - v_y) C_{xy}
.
\]
Note that $i_{xy}=-i_{yx}$.
Kirchhoff's Current Law requires that the total current
flowing into any point other than $a$ or $b$ is 0.
That is, for $x \neq a,b$
\[
\sum_y i_{xy}
= 0
.
\]
This will be true if
\[
\sum_y (v_x - v_y) C_{xy}
= 0
\]
or
\[
v_x \sum_y C_{xy} = \sum_y C_{xy} v_y
.
\]
Thus Kirchhoff's Current Law requires that our voltages have the 
property that
\[
v_x = \sum_y \frac{C_{xy}}{C_x} v_y = \sum_y P_{xy} v_y
\]
for $x \neq a,b$.
This means that the voltage $v_x$ is harmonic at all points
$x \neq a,b$.

Let $h_x$ be the probability, starting at $x$,
that state $a$ is reached before $b$.
Then $h_x$ is also harmonic at all points $x \neq a,b$.
Furthermore
\[
v_a = h_a = 1
\]
and
\[
v_b = h_b = 0
.
\]
Thus if we modify $\matP$ by making $a$ and $b$ absorbing states,
we obtain an absorbing Markov chain $\matPbar$
and $v$ and $h$ are both solutions to the Dirichlet problem for the Markov 
chain with the same boundary values. Hence $v = h$.

For our example,
the transition probabilities $\Pbar_{xy}$ are shown in Figure \figref{3.6}.
\prefig{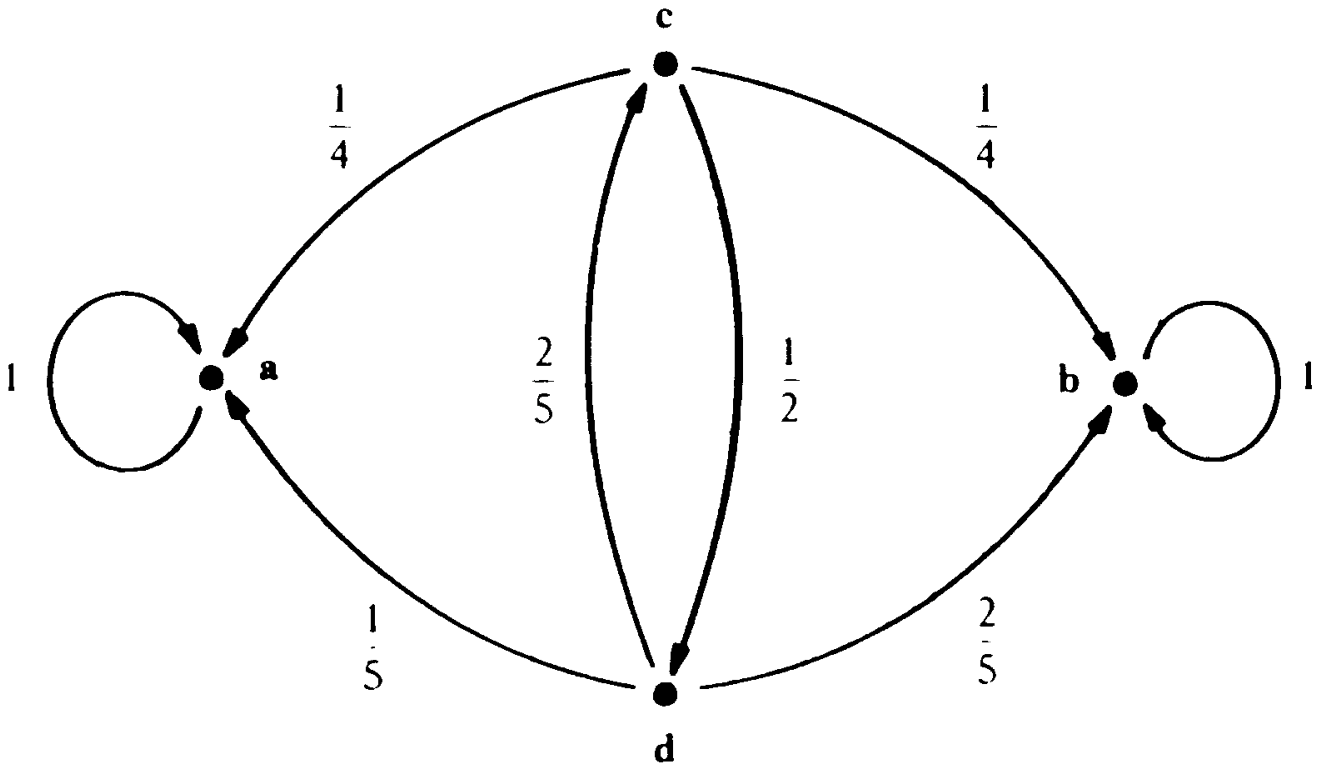}
The function $v_x$ is harmonic for $\Pbar$
with boundary values
$v_a=1, v_b=0$.

To sum up, we have the following:

\tout{Intrepretation of Voltage.}
When a unit voltage is applied between $a$ and $b$,
making
$v_a=1$ and $v_b=0$,
the voltage $v_x$ at any point $x$
represents the probability that a walker starting from $x$ will 
return to $a$ before reaching $b$.
\toutend

In this probabilistic interpretation of voltage,
we have assumed a unit voltage,
but we could have assumed an arbitrary voltage $v_a$
between $a$ and $b$.
Then the hitting probability $h_x$
would be replaced by an expected value
in a game where the player starts at $x$
and is paid $v_a$
if $a$ is reached before $b$
and 0 otherwise.

Let's use this interpretation of voltage
to find the voltages for our example.
Referring back to Figure \figref{3.6},
we see that
\[
v_a=1
\]
\[
v_b = 0
\]
\[
v_c = \frac{1}{4} + \frac{1}{2} v_d
\]
\[
v_d = \frac{1}{5} + \frac{2}{5} v_c
.
\]

Solving these equations yields
$v_c = \frac{7}{16}$ and $v_d = \frac{3}{8}$.
From these voltages we obtain the current $i_{xy}$.
For example
$i_{cd} = (\frac{7}{16} - \frac{3}{8}) \cdot 2 = \frac{1}{8}$.
The resulting voltages and currents are shown in Figure \figref{3.7}. 
\prefig{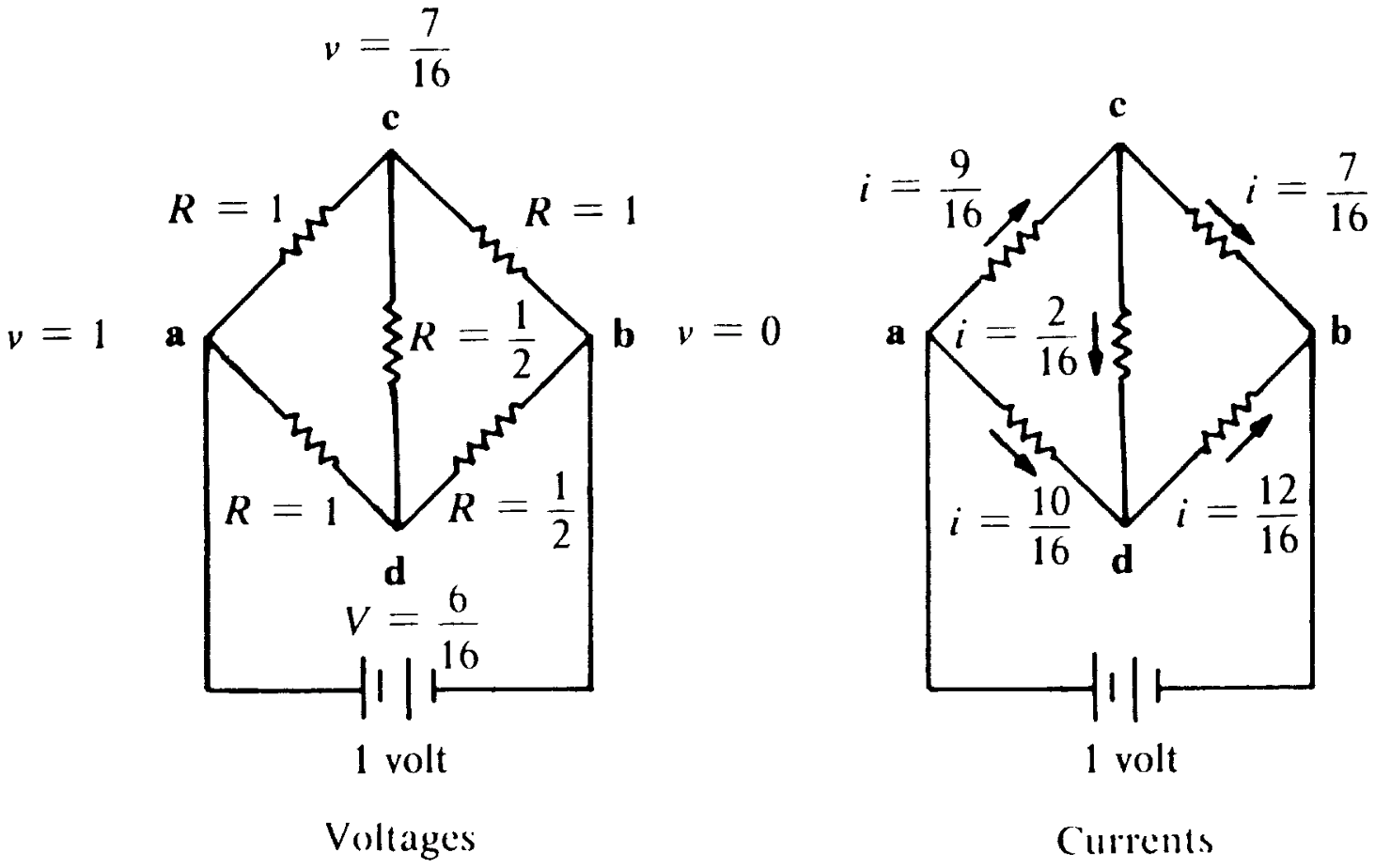}
The voltage at $c$ is $\frac{7}{16}$ and so this is also the probability, 
starting at $c$, of reaching $a$ before $b$.

\sect{3.3}{Probabilistic interpretation of current}

We turn now to the probabilistic interpretation of current.
This interpretation is found by taking a naive view
of the process of electrical conduction:
We imagine that positively charged particles
enter the network at point $a$
and wander around from point to point
until they finally arrive at point $b$,
where they leave the network.
(It would be more realistic to imagine negatively
charged particles entering at $b$ and leaving at $a$,
but realism is not what we're after.)
To determine the current $i_{xy}$
along the branch from $x$ to $y$,
we consider that in the course of its peregrinations
the point may pass once or several times
along the branch from $x$ to $y$,
and in the opposite direction from $y$ to $x$.
We may now hypothesize that the current $i_{xy}$
is proportional to
the expected net number of movements
along the edge from $x$ to $y$,
where movements from $y$ back to $x$ are counted as negative.
This hypothesis is correct,
as we will now show.

The walker begins at $a$ and walks until he reaches $b$;
note that if he returns to $a$ before reaching $b$,
he keeps on going.
Let $u_x$ be the expected number of visits to state $x$
before reaching $b$.
Then $u_b = 0$ and, 
for $x \neq a,b$,
\[
u_x = \sum_y u_y P_{yx}
.
\]
This last equation is true because,
for $x \neq a,b$,
every entrance to $x$ must come from some $y$.

We have seen that $C_x P_{xy} = C_y P_{yx}$;
thus
\[
u_x =
\sum_y u_y \frac{P_{xy} C_x}{C_y}
\]
or
\[
\frac{u_x}{C_x}
=
\sum_y P_{xy} \frac{u_y}{C_y}
.
\]
This means that
$v_x = u_x / C_x$
is harmonic for $x \neq a,b$.
We have also
$v_b =0$ and $v_a = u_a / C_a$.
This implies that $v_x$ is the voltage at $x$ when 
we put a battery from $a$ to $b$
that establishes a voltage
$u_a/C_a$ at $a$
and voltage 0 at $b$.
(We remark that the expression $v_x = u_x / C_x$
may be understood physically by viewing $u_x$ as charge
and $C_x$ as capacitance;
see Kelly \cite{kelly}
for more about this.)

We are interested in the current that flows from $x$ to $y$.
This is
\[
i_{xy}
=
(v_x -v_y) C_{xy}
=
\braces{ \frac{u_x}{C_x} - \frac{u_y}{C_y} } C_{xy}
=
\frac{u_x C_{xy}}{C_x} - \frac{u_y C_{yx}}{C_y}
=
u_x P_{xy} - u_y P_{yx}
.
\]
Now $u_x P_{xy}$ is
the expected number of times our walker will go from $x$ to $y$ 
and $u_y P_{yx}$ is the expected number of times he will go from $y$ to $x$. 
Thus the current $i_{xy}$
is the expected value for the net number of times the 
walker passes along the edge from $x$ to $y$.
Note that for any particular walk this net value will be an integer,
but the expected value will not.

As we have already noted,
the currents $i_{xy}$ here are not those of our
original electrical problem,
where we apply a 1-volt battery,
but they are proportional to those original currents.
To determine the constant of proportionality,
we note the following characteristic property of 
the new currents $i_{xy}$:
The total current flowing into the network at $a$ (and out at $b$)
is 1.
In symbols,
\[
\sum_y i_{ay} = 1
.
\]
Indeed,
from our probabilistic interpretation of $i_{xy}$
this sum represents the expected value
of the difference between the number of times our walker
leaves $a$ and enters $a$.
This number is necessarily one and so the current flowing into $a$ is 1.

This unit current flow from $a$ to $b$
can be obtained from the currents in the original circuit,
corresponding to a 1-volt battery,
by dividing through by the total amount of current
$\sum_y i_{ay}$
flowing into $a$;
doing this to the currents in our example yields
the unit current flow shown in Figure \figref{3.8}.
\prefig{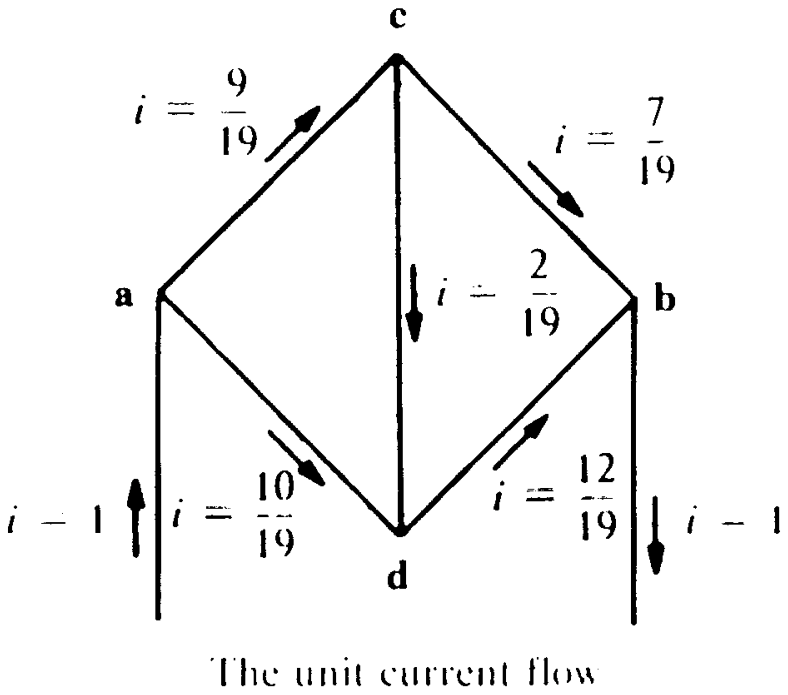}

This shows that the constant of proportionality
we were seeking to determine
is the reciprocal of the amount of current
that flows through the circuit
when a I-volt battery is applied between $a$ and $b$.
This quantity,
called the effective resistance between $a$ and $b$,
is discussed in detail in \sectref{3.4}.

To sum up, we have the following:

\tout{Interpretation of Current.}
When a unit current flows into $a$ and out of $b$,
the current $i_{xy}$
flowing through the branch connecting $x$ to $y$
is equal to
the expected net number of times that a walker,
starting at $a$ and walking until he reaches $b$,
will move along the branch from $x$ to $y$.
These currents are proportional to
the currents that arise when a unit voltage is applied between $a$ and $b$,
the constant of proportionality being the 
effective resistance of the network.
\toutend

We have seen that we can estimate the voltages by simulation.
We can now do the same for the currents.
We have to estimate the expected value
for the net number of crossings of $xy$.
To do this,
we start a large number of walks at $a$
and, for each one, record the net number of crossings of each edge
and then average these as an estimate for the expected value.
Carrying out 10,000 such walks yielded the results shown in 
Figure \figref{3.9}.
\prefig{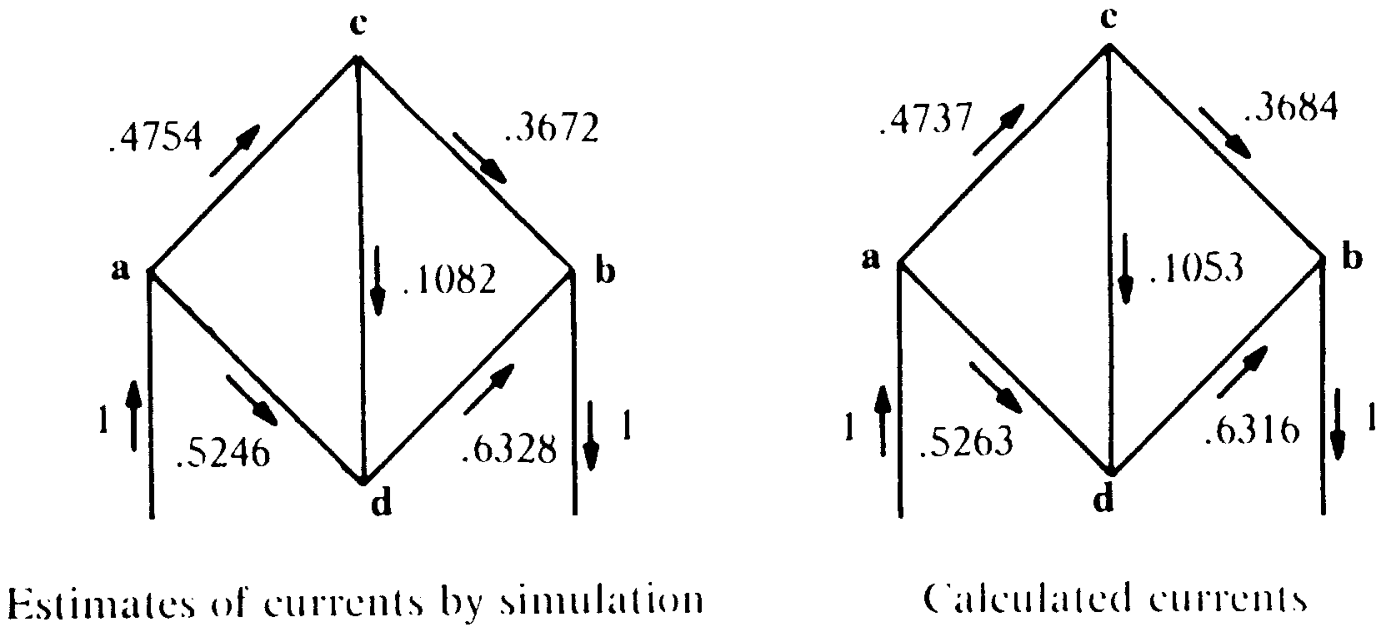}

The results of simulation are in good agreement with
the theoretical values of current.
As was the case for estimates of the voltages by simulation,
we have statistical errors.
Our estimates have the property that the total current flowing into $a$ is 1,
out of $b$ is 1,
and into any other point it is 0.
This is no accident;
the explanation is that the history of each walk
would have these properties,
and these properties are not destroyed by averaging.

\exstart
\ex{3.3.1}{
Kingman \cite{kingman}
introduced a different model for current flow.
Kelly \cite{kelly}
gave a new interpretation of this model.
Both authors use continuous time.
A discrete time version of Kelly's interpretation
would be the following:
At each point of the graph there is a black or a white button.
Each second an edge is chosen;
edge $xy$ is chosen with probability $C_{xy}/C$
where $C$ is the sum of the conductances.
The buttons on the edge chosen are then interchanged.
When a button reaches $a$ it is painted black,
and when a button reaches $b$ it is painted white.
Show that there is a limiting probability $p_x$
that site $x$ has a black button
and that $p_x$ is the voltage $v_x$ at $x$
when a unit voltage is imposed between $a$ and $b$.
Show that the current $i_{xy}$ is proportional to
the net flow of black buttons along the edge $xy$.
Does this suggest a hypothesis about
the behavior of conduction electrons in metals?
}
\exend

\sect{3.4}{Effective resistance and the escape probability}

When we impose a voltage $v$ between points $a$ and $b$,
a voltage $v_a = v$ is established at $a$
and $v_b=0$,
and a current $i_a = \sum_x i_{ax}$
will flow into the circuit from the outside source.
The amount of current that flows depends upon
the overall resistance in the circuit.
We define the \dt{effective resistance}
$\Reff$ between $a$ and $b$ by
$\Reff=v_a/i_a$.
The reciprocal quantity $\Ceff = 1/\Reff = i_a/v_a$
is the \dt{effective conductance}.
If the voltage between $a$ and $b$ is 
multiplied by a constant,
then the currents are multiplied by the same constant,
so $\Reff$ depends only on the ratio of $v_a$ to $i_a$.

Let us calculate $\Reff$ for our example.
When a unit voltage was imposed,
we obtained the currents shown in Figure \figref{3.7}.
The total current flowing into the circuit is
$i_a = 9/16 + 10/16 = 19/16$.
Thus the effective resistance is
\[
\Reff
=
\frac{v_a}{i_a}
=
\frac{1}{\frac{19}{16}}
=
\frac{16}{19}
.
\]

We can interpret the effective conductance probabilistically
as an escape probability.
When $v_a=1$,
the effective conductance equals the total current $i_a$ flowing into $a$.
This current is
\[
i_a
=
\sum_y (v_a - v_y) C_{ay}
=
\sum_y (v_a - v_y) \frac{C_{ay}}{C_a} C_a
=
C_a (1-\sum_y P_{ay} v_y)
=
C_a p_\esc
\]
where $p_\esc$ is the probability, starting at $a$,
that the walk reaches $b$ before returning to $a$.
Thus
\[
\Ceff = C_a p_\esc
\]
and
\[
p_\esc = \frac{\Ceff}{C_a}
.
\]

In our example $C_a = 2$
and we found that $i_a = 19/16$.
Thus
\[
p_\esc = \frac{19}{32}
.
\]

In calculating effective resistances,
we shall use two important facts about electric networks.
First, if two resistors are connected in series,
they may be replaced by a single resistor whose resistance
is the sum of the two resistances.
(See Figure \figref{3.10}.)
\prefig{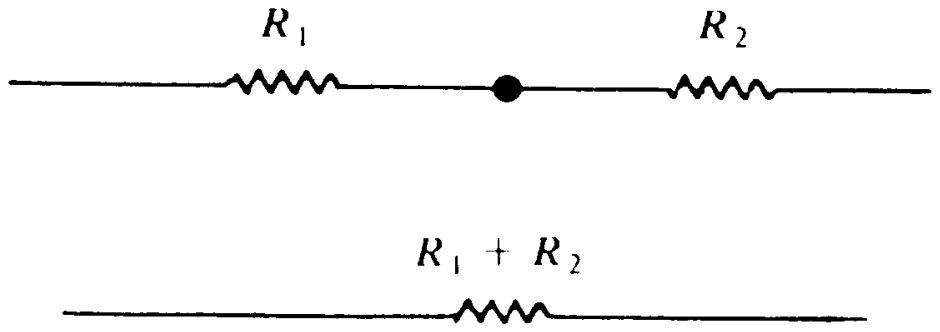}
Secondly, two resistors in parallel may be replaced by a single resistor
with resistance $R$ such that
\[
\recip{R} = \recip{R_1} + \recip{R_2}
=
\frac{R_1 R_2}{R_1 + R_2}
.
\]
(See Figure \figref{3.11}.)
\prefig{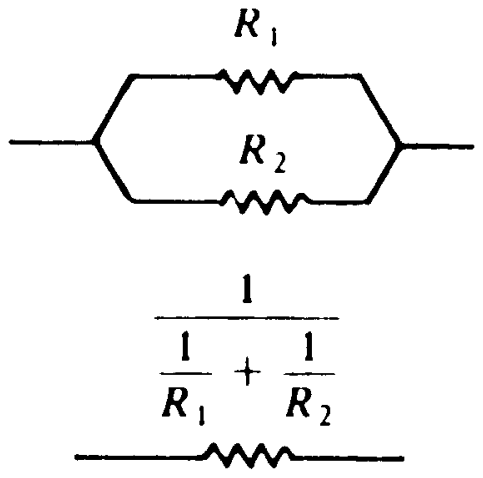}

The second rule can be stated more simply in terms of conductances:
If two resistors are connected in parallel, they may be replaced by 
a single resistor whose conductance is the sum of the two conductances.

We illustrate the use of these ideas to compute the effective 
resistance between two adjacent points of a unit cube of unit resistors,
as shown in Figure \figref{3.12}.
\prefig{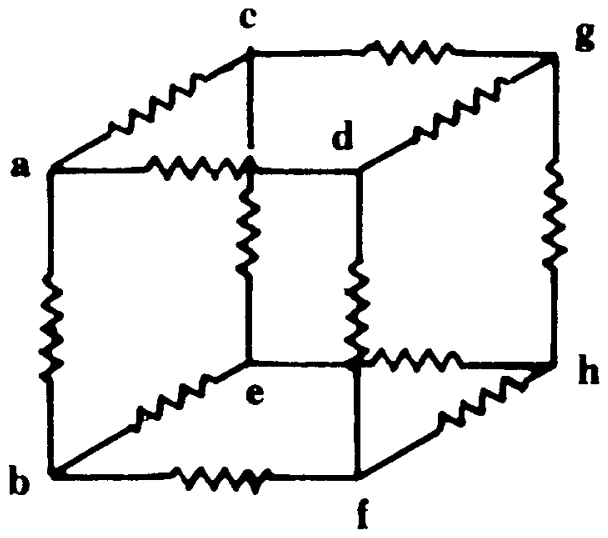}
We put a unit battery between $a$ and $b$.
Then, by symmetry, the voltages at $c$ and $d$
will be the same as will those at $e$ and $f$.
Thus our circuit is equivalent to
the circuit shown in Figure \figref{3.13}.
\prefig{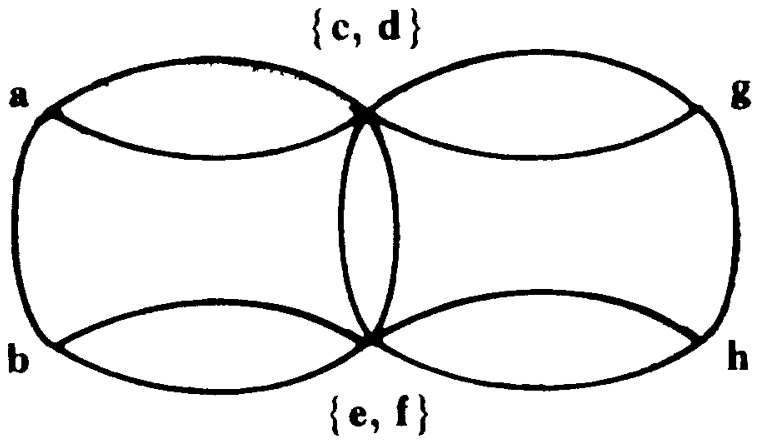}

Using the laws for the effective resistance
of resistors in series and parallel,
this network can be successively reduced to a single 
resistor of resistance $7/12$ ohms, as shown in Figure \figref{3.14}.
\prefig{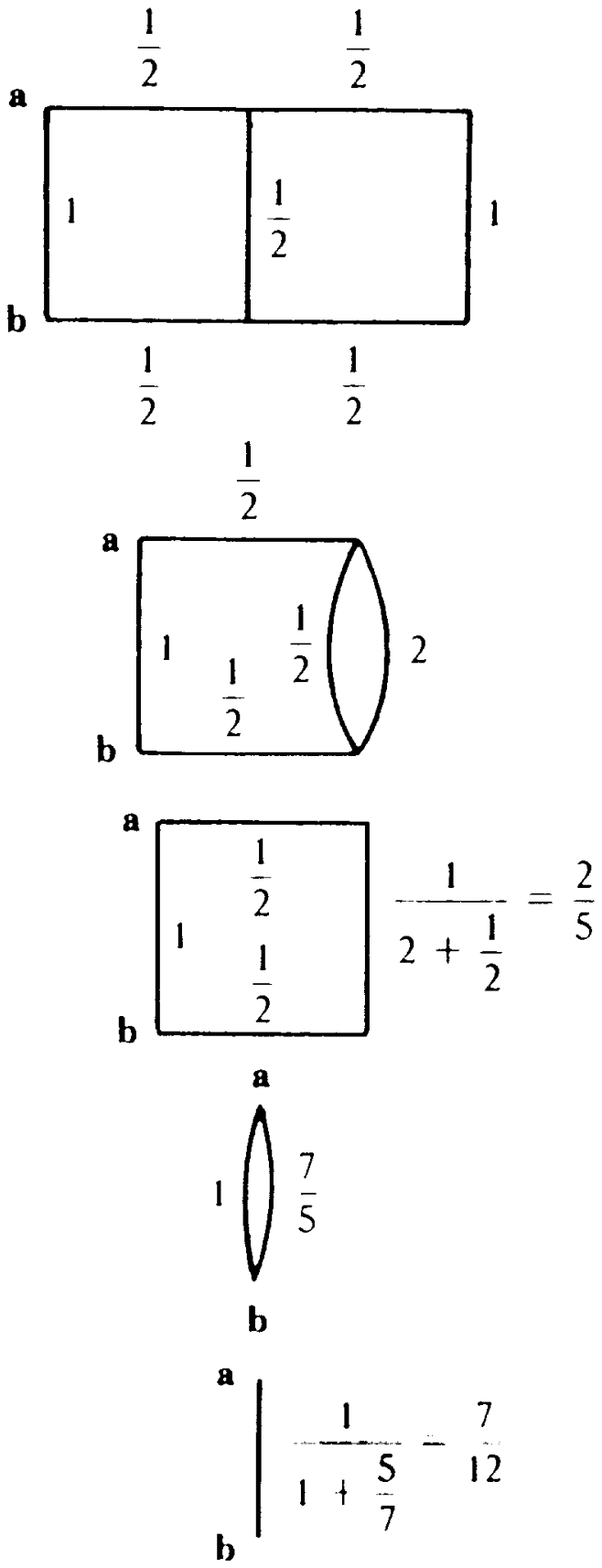}
Thus the effective resistance is $7/12$.
The current flowing into $a$ from the battery will
be $i_a = \recip{\Reff} = 12/7$.
The probability that a walk starting at $a$ will reach $b$
before returning to $a$ is
\[
p_\esc = \frac{i_a}{C_a} = \frac{\frac{12}{7}}{3} = \frac{4}{7}
.
\]

This example and many other interesting connections between 
electric networks and graph theory
may be found in Bollobas \cite{bollobas}.

\exstart
\ex{3.4.1}{
A bug walks randomly on the unit cube
(see Figure \figref{3.15}).
\prefig{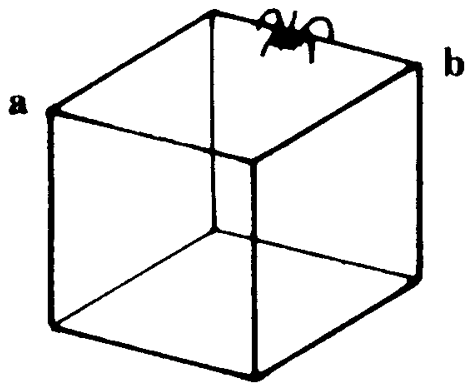}
If the bug starts at $a$,
what is the probability that it reaches food at $b$
before returning to $a$?
}

\ex{3.4.2}{
Consider the Ehrenfest urn model with $N = 4$
(see Exercise \exref{3.1.5}).
Find the probability, starting at 0, that state 4 is 
reached before returning to 0.
}

\ex{3.4.3}{
Consider the ladder network shown in
Figure \figref{3.16}.
\prefig{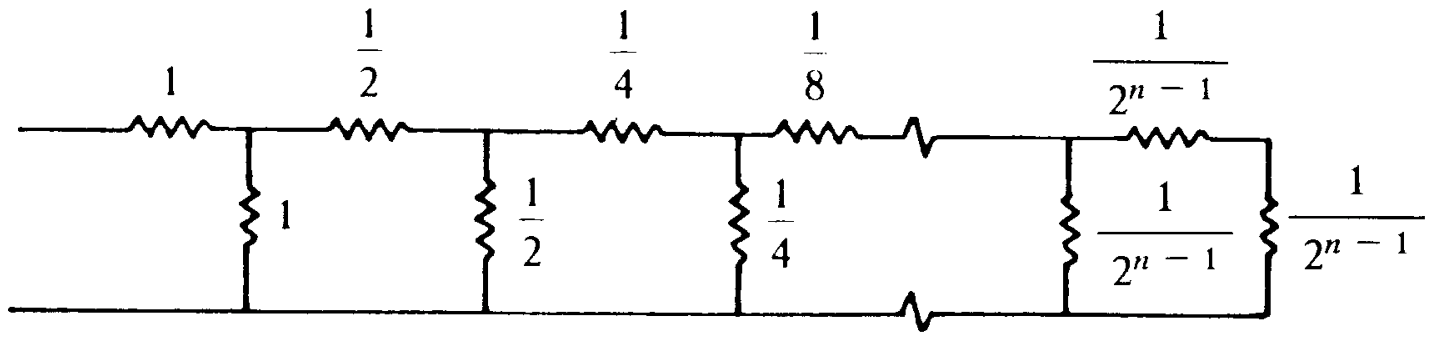}
Show that if $R_n$ is the effective resistance
of a ladder with $n$ rungs then $R_1 = 2$ and
\[
R_{n+1} = \frac{2+2R_n}{2+R_n}
.
\]
Use this to show that
$\lim_{n \goto \infty} R_n = \sqrt{2}$.
}

\ex{3.4.4}{
A drunken tourist starts at her hotel and walks at random 
through the streets of the idealized Paris shown in
Figure \figref{3.17}.
\prefig{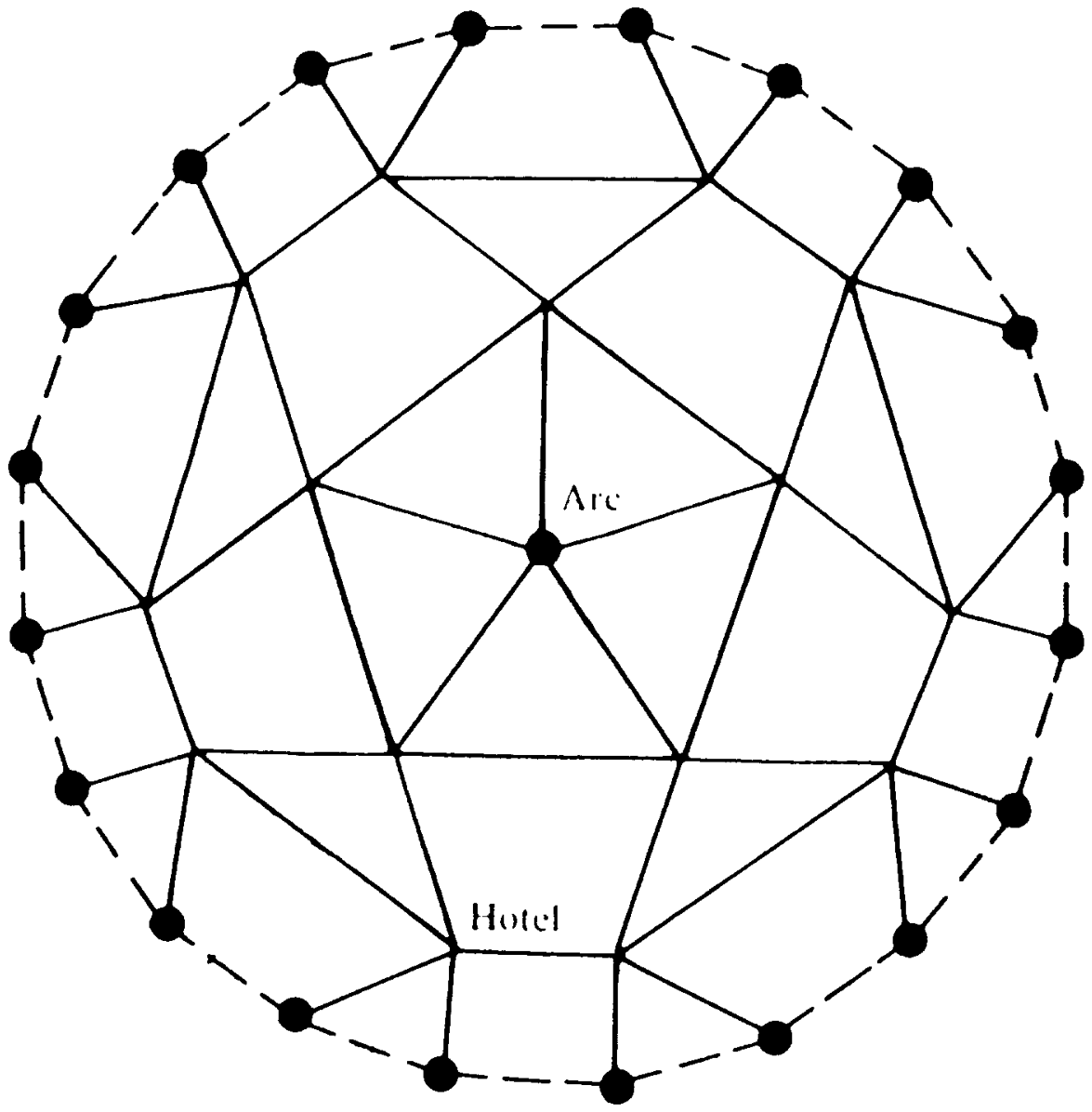}
Find the probability that she reaches the Arc de Triomphe
before she reaches the outskirts of town.
}
\exend

\sect{3.5}{Currents minimize energy dissipation}

We have seen that when we impose a voltage between $a$ and $b$
voltages $v_x$ are established at the points and currents $i_{xy}$
flow through the resistors.
In this section we shall give a characterization of the currents in 
terms of a quantity called \dt{energy dissipation}.
When a current $i_{xy}$ flows through a resistor,
the energy dissipated is
\[
i_{xy}^2 R_{xy}
;
\]
this is the product of the current $i_{xy}$
and the voltage $v_{xy} = i_{xy} R_{xy}$.
The
\dt{total energy dissipation} in the circuit is
\[
E = \half \sum_{x,y} i_{xy}^2 R_{xy}
.
\]
Since $i_{xy} R_{xy} = v_x - v_y$,
we can also write the energy dissipation as
\[
E = \half \sum_{x,y} i_{xy} (v_x - v_y)
.
\]
The factor 1/2 is necessary in this formulation since
each edge is counted twice in this sum.
For our example, we see from Figure \figref{3.7} that
\[
E
=
\braces{\frac{9}{16}}^2 \cdot 1
+
\braces{\frac{10}{16}}^2 \cdot 1
+
\braces{\frac{7}{16}}^2 \cdot 1
+
\braces{\frac{2}{16}}^2 \cdot \half
+
\braces{\frac{12}{16}}^2 \cdot \half
=
\frac{19}{16}
.
\]

If a source (battery) establishes voltages
$v_a$ and $v_b$ at $a$ and $b$,
then the energy supplied is $(v_a - v_b) i_a$
where $i_a = \sum_x i_{ax}$.
By conservation of energy,
we would expect this to be equal to the energy dissipated.
In our example $v_a - v_b = 1$ and
$i_a = \frac{19}{16}$,
so this is the case.
We shall show that this is true in a somewhat more general setting.

Define a \dt{flow} $\mat{j}$ from $a$ to $b$
to be an assignment of numbers $j_{xy}$ to pairs $xy$
such that
\propstart
\prop{(a)}
$j_{xy} = - j_{yx}$
\prop{(b)}
$\sum_y j_{xy} = 0$ if $x \neq a,b$
\prop{(c)}
$j_{xy} = 0$ if $x$ and $y$ are not adjacent.
\propend

We denote by $j_x = \sum_y j_{xy}$ the flow into $x$ from the outside.
By (b)
$j_x = 0$ for $x \neq a,b$.
Of course $j_b = - j_a$.
To verify this, note that
\[
j_a + j_b
=
\sum_x j_x
=
\sum_x \sum_y j_{xy}
=
\half
\sum_{x,y} (j_{xy} + j_{yx})
=
0
,
\]
since $j_{xy} = - j_{yx}$.

With this terminology,
we can now formulate the following version of
the principle of conservation of energy:

\tout{Conservation of Energy.}
Let $w$ be any function defined on the points of the graph
and $\mat{j}$ a flow from $a$ to $b$.
Then
\[
(w_a - w_b) j_a
=
\half
\sum_{x,y}
(w_x - w_y) j_{xy}
.
\]
\toutend

\proofstart
\begin{eqnarray*}
\sum_{x,y} (w_x - w_y) j_{xy}
&=&
\sum_x (w_x \sum_y j_{xy}) - \sum_y (w_y \sum_x j_{xy})
\\&=&
w_a \sum_y j_{ay} + w_b \sum_y j_{by}
- w_a \sum_x j_{xa} - w_b \sum_x j_{xb}
\\&=&
w_a j_a + w_b j_b - w_a (-j_a) - w_b (-j_b)
\\&=&
2 (w_a - w_b) j_a
.
\end{eqnarray*}
Thus
\[
(w_a - w_b) j_a
=
\half
\sum_{x,y} (w_x - w_y) j_{xy}
\]
as was to be proven.
\proofend

If we now impose a voltage $v_a$ between $a$ and $b$
with $v_b=0$,
we obtain voltages $v_x$ and currents $i_{xy}$.
The currents $\mat{i}$ give a flow from $a$ to $b$
and so by the previous result,
we conclude that
\[
v_a i_a
=
\half
\sum_{x,y} (v_x - v_y) i_{xy} 
=
\half
\sum_{x,y} i_{xy}^2 R_{xy}
.
\]
Recall that $\Reff = v_a / i_a$.
Thus in terms of resistances we can write this as
\[
i_{xy}^2 \Reff = \half \sum_{x,y} i_{xy}^2 R_{xy}
.
\]

If we adjust $v_a$ so that $i_a = 1$,
we call the resulting flow \dt{the unit current flow}
from $a$ to $b$.
The unit current flow from $a$ to $b$
is a particular example of a \dt{unit flow} from $a$ to $b$,
which we define to be any flow $i_{xy}$ from $a$ to $b$
for which $i_a = -i_b = 1$.
The formula above shows that
the energy dissipated by the unit current flow
is just $\Reff$.
According to a basic result called Thomson's Principle,
this value is smaller than the energy dissipated by any other 
unit flow from $a$ to $b$.
Before proving this principle, let us watch it in action 
in the example worked out above.

Recall that, for this example, we found the true values and some
approximate values for the unit current flow; these were shown in 
Figure \figref{3.9}.
The energy dissipation for the true currents is
\[
E
=
\Reff
=
\frac{16}{19}
= .8421053
.
\]
Our approximate currents also form a unit flow and, for these, 
the energy dissipation is
\[
\bar{E} = (.4754)^2 \cdot 1 + (.5246)^2 \cdot 1 + (.3672)^2  \cdot 1
+ (.1082)^2 \cdot \half + (.6328)^2 \cdot \half
=
.8421177
.
\]

We note that $\bar{E}$ is greater than $E$, though just barely.

\tout{Thomson's Principle. (Thomson \cite{thomsonTait}).}
If $\mat{i}$ is the unit flow from $a$ to $b$
determined by Kirchhoff's Laws,
then the energy dissipation 
$\half \sum_{x,y} i_{xy}^2 R_{xy}$
minimizes the energy dissipation
$\half \sum_{x,y} j_{xy}^2 R_{xy}$
among all unit flows $\mat{j}$ from $a$ to $b$.
\toutend

\proofstart
Let $\mat{j}$ be any unit flow from $a$ to $b$
and let
$d_{xy} = j_{xy} - i_{xy}$.
Then $\mat{d}$ is a flow from $a$ to $b$ with
$d_a = \sum_x d_{ax} = 1-1 = 0$.
\begin{eqnarray*}
\sum_{x,y} j_{xy}^2 R_{xy}
&=&
\sum_{x,y} (i_{xy} + d_{xy})^2 R_{xy}
\\&=&
\sum_{x,y} i_{xy}^2 R_{xy}
+ 2 \sum_{x,y} i_{xy} R_{xy} d_{xy}
+ \sum_{x,y} d_{xy}^2 R_{xy}
\\&=&
\sum_{x,y} i_{xy}^2 R_{xy}
+ 2 \sum_{x,y} (v_x - v_y) d_{xy}
+ \sum_{x,y} d_{xy}^2 R_{xy}
.
\end{eqnarray*}
Setting $\mat{w} = \mat{v}$ and $\mat{j} = \mat{d}$
in the conservation of energy result above shows that 
the middle term is $4 (v_a - v_b) d_a = 0$.
Thus
\[
\sum_{x,y} j_{xy}^2 R_{xy}
=
\sum_{x,y} i_{xy}^2 R_{xy} + \sum_{x,y} d_{xy}^2 R_{xy}
\geq
\sum_{x,y} i_{xy}^2 R_{xy}
.
\]
This completes the proof.
\proofend

\exstart
\ex{3.5.1}{
The following is the so-called ``dual form'' of Thomson's Principle.
Let $u$ be any function on the points of the graph $G$ of a circuit
such that $u_a=1$ and $u_b=0$.
Then the energy dissipation
\[
\half \sum_{x,y} (u_x - u_y)^2 C_{xy}
\]
is minimized by the voltages $v_x$
that result when a unit voltage is established between $a$ and $b$,
i.e., $v_a=1$, $v_b=0$, and the other voltages are 
determined by Kirchhoff's Laws.
Prove this dual principle.
This second principle is known nowadays as the \dt{Dirichlet Principle},
though it too was discovered by Thomson.
}

\ex{3.5.2}{
In \sectref{2.4}\ we stated that,
to solve the Dirichlet problem 
by the method of relaxations,
we could start with an arbitrary initial guess.
Show that when we replace the value at a point
by the average of the neighboring points
the energy dissipation, as expressed in Exercise \exref{3.5.1},
can only decrease.
Use this to prove that the relaxation method converges
to a solution of the Dirichlet problem
for an arbitrary initial guess.
}
\exend

\chap{4}{Rayleigh's Monotonicity Law}

\sect{4.1}{Rayleigh's Monotonicity Law}

Next we will study Rayleigh's Monotonicity Law.
This law from electric network theory
will be an important tool in our future study of random walks.
In this section we will give an example of the use of this law.

Consider a random walk on streets of a city as in Figure \figref{4.1}. 
\prefig{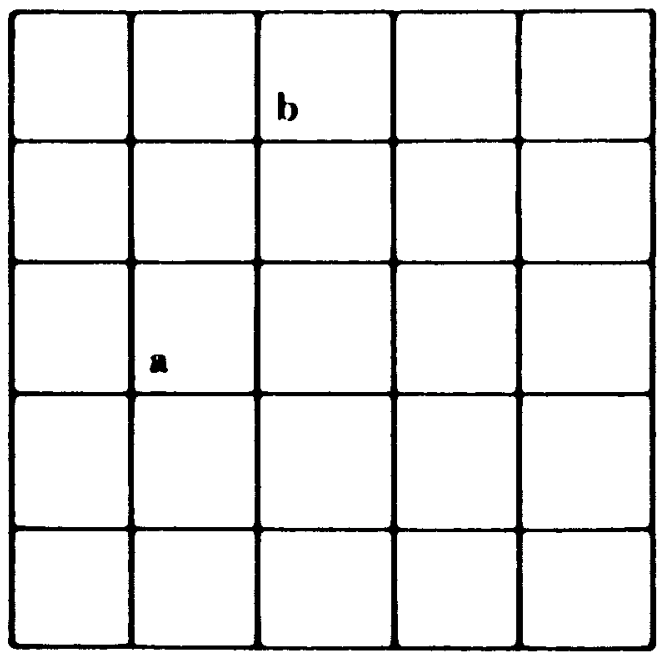}
Let $p_\esc$ be the probability that a walker starting from $a$ reaches $b$ 
before returning to $a$.
Assign to each edge a unit resistance
and maintain a voltage of one volt between $a$ and $b$;
then a current $i_a$  will flow into the circuit and
we showed in \sectref{3.4}\ that
\[
p_\esc
=
\frac{i_a}{C_a}
=
\frac{i_a}{4}
.
\]

Now suppose that one of the streets (not connected to $a$)
becomes blocked.
Our walker must choose from the remaining streets
if he reaches a corner of this street.
The modified graph will be as in Figure \figref{4.2}.
\prefig{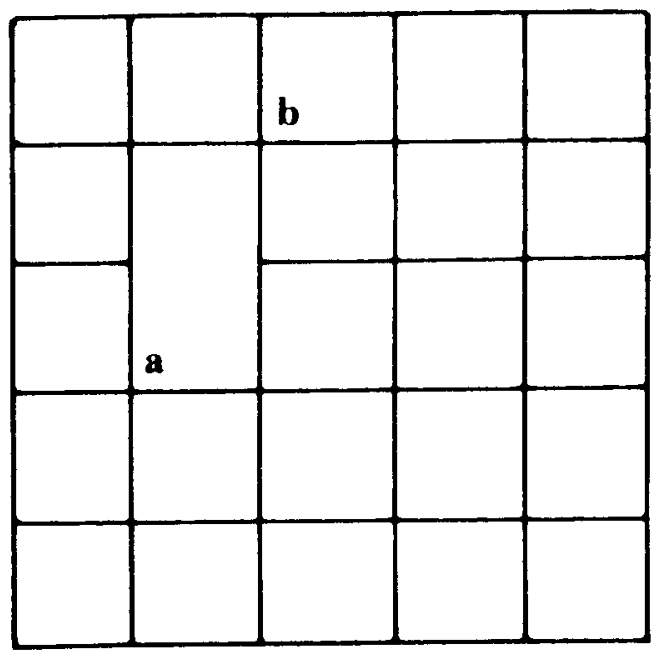}
We want to show that the probability of escaping to $b$ from $a$
is decreased.

Consider this problem in terms of our network.
Blocking a street corresponds to replacing a unit resistor
by an infinite resistor. 
This should have the effect of increasing the effective resistance
$\Reff$ of the circuit between $a$ and $b$.
If so, when we put a unit voltage between $a$ and $b$
less current will flow into the circuit and
\[
p_\esc
=
\frac{i_a}{4}
=
\frac{1}{4 \Reff}
\]
will decrease.

Thus we need only show that when we increase the resistance in one
part of a circuit,
the effective resistance increases.
This fact, known as
Rayleigh's Monotonicity Law, is almost self-evident.
Indeed, the father of electromagnetic theory,
James Clerk Maxwell, regarded this to be the case.
In his \emph{Treatise on Electricity and Magnetism}
(\cite{maxwell}, p. 427), he wrote
\begin{quote}
If the specific resistance of any portion of the conductor be changed,
that of the remainder being unchanged,
the resistance of the whole conductor will be increased
if that of the portion is increased, 
and diminished if that of the portion is diminished.
This principle may be regarded as self-evident \ldots .
\end{quote}

\tout{Rayleigh's Monotonicity Law.}
If the resistances of a circuit are increased,
the effective resistance $\Reff$ between any two points
can only increase.
If they are decreased, it can only decrease.
\toutend

\proofstart
Let $\mat{i}$ be the unit current flow from $a$ to $b$
with the resistors $R_{xy}$.
Let $\mat{j}$ be the unit current flow from $a$ to $b$
with the resistors $\bar{R}_{xy}$
with $\bar{R}_{xy} \geq R_{xy}$.
Then
\[
\bar{R}_\eff
=
\half \sum_{x,y} j_{xy}^2 \bar{R}_{xy}
\geq
\half \sum_{x,y} j_{xy}^2 R_{xy}
.
\]
But since $\mat{j}$ is a unit flow from $a$ to $b$,
Thomson's Principle tells us that the energy dissipation,
calculated with resistors $R_{xy}$,
is bigger than that for
the true currents determined by these resistors:
that is
\[
\half \sum_{x,y} j_{xy}^2 R_{xy}
\geq \half \sum_{x,y} i_{xy}^2 R_{xy}
= \Reff
.
\]
Thus $\bar{R}_\eff \geq \Reff$.
The proof for the case of decreasing resistances is the same.

\exstart

\ex{4.1.1}{
Consider a graph $G$ and let $R_{xy}$ and $\bar{R}_{xy}$
be two different assignments of resistances to the edges of $G$.
Let $\hat{R}_{xy} = \bar{R}_{xy} + R_{xy}$.
Let $\Reff$, $\bar{R}_\eff$, and $\hat{R}_\eff$
be the effective resistances when $R$, $\bar{R}$, and $\hat{R}$,
respectively, are used.
Prove that
\[
\hat{R}_\eff \geq \bar{R}_\eff + \Reff
.
\]
Conclude that the effective resistance of a network is a concave function
of any of its component resistances
(Shannon and Hagelbarger \cite{shannonhagelbarger}.)
}

\ex{4.1.2}{
Show that the effective resistance
of the circuit in Figure \figref{4.3}
is greater than or equal to the effective resistance
of the circuit in Figure \figref{4.4}.
\prefig{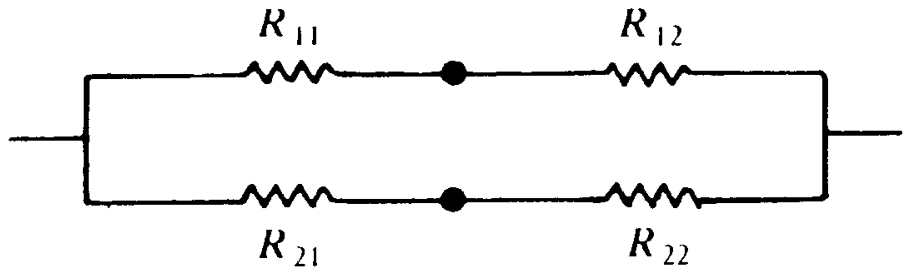}
\prefig{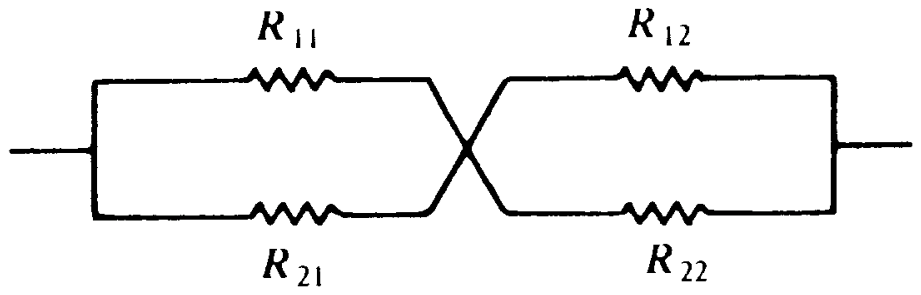}
Use this to show the following inequality for $R_{ij} \geq 0$:
\[
\recip{\recip{R_{11}+R_{12}}+\recip{R_{21}+R_{22}}}
\geq
\recip{\recip{R_{11}}+\recip{R_{21}}}
+
\recip{\recip{R_{12}}+\recip{R_{22}}}
.
\]
See the note by Lehman \cite{lehman}
for a proof of the general Minkowski inequality by this method.
}

\ex{4.1.3}{
Let $\matP$ be the transition matrix associated with an electric network
and let $a, b, r, s$ be four points on the network.
Let $\matPbar$ be a transition matrix defined on the state-space
$S = \{a, b, r, s\}$.
Let $\Pbar_{ii} = 0$
and for $i \neq j$
let $\Pbar_{ij}$ be the probability that,
if the chain $\matP$ is started in state $i$,
then the next time it is in the set $S - \{ i \}$
it is in the state $j$.
Show that $\matPbar$ is a reversible Markov chain and 
hence corresponds to an electric network of the form of a
Wheatstone Bridge,
shown in Figure \figref{4.5}.
\prefig{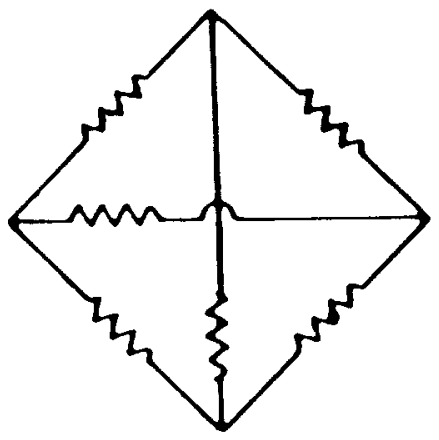}
Explain how this proves that,
in order to prove the Monotonicity Law,
it is sufficient to prove that $\Reff$ is a monotone function
of the component resistances for a Wheatstone Bridge.
Give a direct proof of the Monotonicity Law
for this special case.
}
\exend

\sect{4.2}{A probabilistic explanation of the Monotonicity Law}

We have quoted Maxwell's assertion that Rayleigh's Monotonicity Law
may be regarded as self-evident,
but one might feel that any argument in terms of electricity
is only self-evident if we know what electricity is.
In Cambridge, they tell the following story about Maxwell:
Maxwell was lecturing and, seeing a student dozing off,
awakened him, asking,
``Young man, what is electricity?''
``I'm terribly sorry, sir,''
the student replied,
~`I knew the answer but I have forgotten it.''
Maxwell's response to the class was,
``Gentlemen, you have just witnessed the greatest tragedy
in the history of science.
The one person who knew what electricity is has forgotten it.''

To say that our intuition about the Monotonicity Law is only as 
solid as our understanding of electricity is not really a valid argument, 
of course,
because in saying that this law is self-evident we are secretly 
depending on the analogy between electricity and the flow of water
(see Feynman \cite{feynman}, Vol. 2, Chapter 12).
We just can't believe that if a water main gets clogged
the total rate of flow out of the local reservoir
is going to increase.
But as soon as we admit this,
some pedant will ask if we're talking about flows with low Reynolds number,
or what, and we'll have to admit that we don't
understand water any better than we understand electricity.

Whatever our feelings about electricity or the flow of water,
it seems desirable to have an explanation of the Monotonicity Law
in terms of our random walker.
We now give such an explanation.

We begin by collecting some notation and results from previous sections.
As usual, we have a network of conductances (streets) 
and a walker who moves from point $x$ to point $y$
with probability
\[
P_{xy} = \frac{C_{xy}}{C_x}
\]
where $C_{xy}$ is the conductance from $x$ to $y$
and $C_x = \sum_y C_{xy}$.
We choose two preferred points $a$ and $b$.
The walker starts at $a$ and walks 
until he reaches $b$ or returns to $a$.
We denote by $v_x$ the probability that the walker,
starting at $a$,
reaches $a$ before $b$.
Then
$v_a=1$,
$v_b=0$,
and the function $v_x$ is harmonic at all points $x \neq a,b$.
We denote by $p_\esc$ the probability that the walker,
starting at $a$,
reaches $b$ before returning to $a$.
Then
\[
p_\esc = 1 - \sum_x p_{ax} v_x
.
\]

Now we have seen that the effective conductance between $a$ and $b$ is
\[
C_a p_\esc
.
\]
We wish to show that this increases whenever one of the conductances
$C_{rs}$ is increased.
If $a$ is different from $r$ or $s$,
we need only show that $p_\esc$ increases.
The case where $r$ or $s$ coincides with $a$ is easily 
disposed of (see Exercise \exref{4.2.1}).
The case where $r$ or $s$ coincides with $b$ is also easy 
(see Exercise \exref{4.2.2}).
Hence from now on we will assume that $r,s \neq a$ and $r,s \neq b$.

Instead of increasing $C_{rs}$,
we can think of adding a new edge of conductance $\epsilon$
between $r$ and $s$.
(See Figure \figref{4.6}.)
\prefig{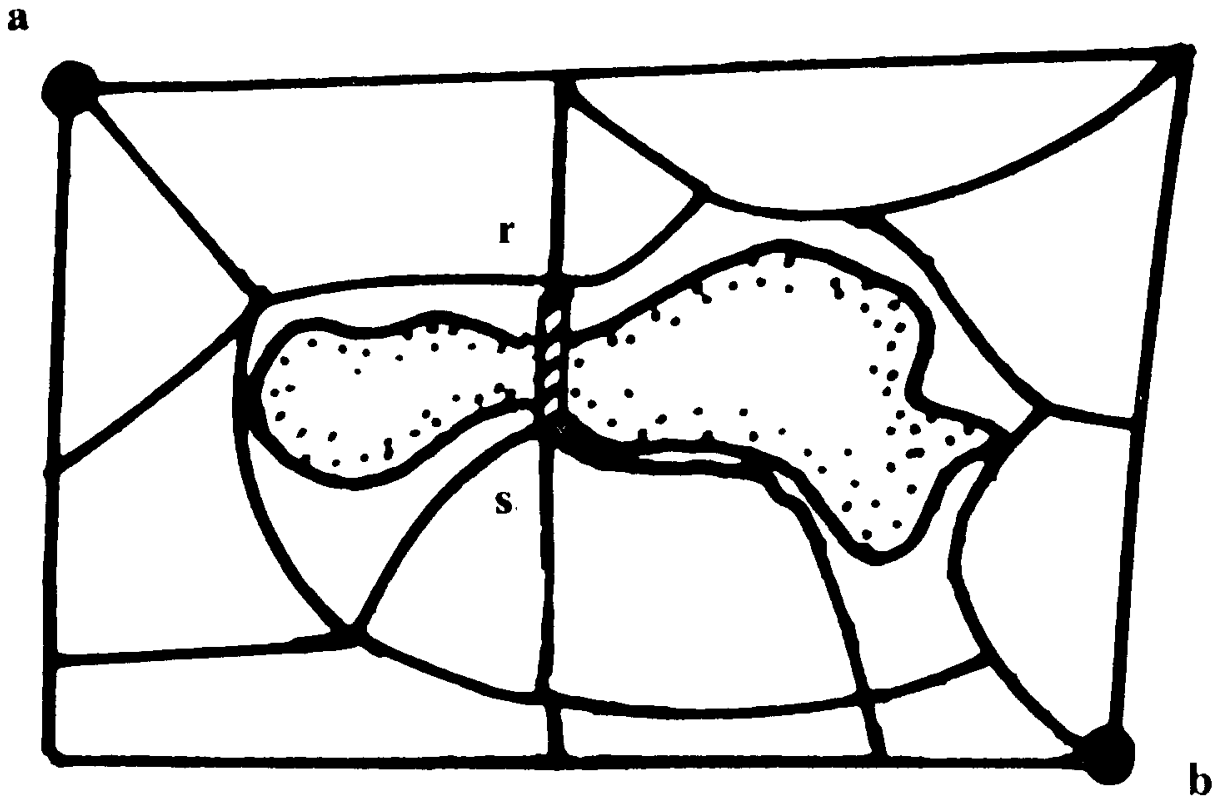}
We will call this new edge a ``bridge''
to distinguish it from the other edges.
Note that the graph with the bridge added
will have more than one edge between $r$ and $s$
(unless there was no edge between $r$ and $s$ in the original graph), 
and this will complicate any expression that involves summing over edges.
Everything we have said or will say holds for graphs 
with ``multiple edges'' as well as for graphs without them.
So far, we have chosen to keep our expressions simple
by assuming that an edge is determined by its
endpoints.
The trade-off is that in the manipulations below, 
whenever we write a sum over edges
we will have to add an extra term to account for the bridge.

Why should adding the bridge increase the escape probability?
The first thing you think is,
``Of course, it opens up new possibilities of escaping!''
The next instant you think,
``Wait a minute, it also opens up
new possibilities of returning to the starting point.
What ensures that the first effect will outweigh the second?''
As we shall see, the proper reply is,
``Because the walker will cross the bridge more often in the good 
direction than in the bad direction.''
To turn this blithe reply into a real explanation
will require a little work, however.

To begin to make sense of the notion that the bridge gets used 
more often in the good direction than the bad,
we will make a preliminary argument
that applies to any edge of any graph.
Let $G$ be any graph,
and let $rs$ be any edge with endpoints not $a$ or $b$.
$v_r > v_s$.
Since the walker has a better chance to escape from $s$ than from $r$,
this means that to cross this edge in the good direction is to go from 
$r$ to $s$
 We shall show that the walker will cross the edge from $r$ to $s$
more often on the average than from $s$ to $r$.

Let $u_x$ be the expected number of times the walker is at $x$
and $u_{xy}$ the expected number of times he crosses the edge $xy$
from $x$ to $y$ before he reaches $b$ or returns to $a$.
The calculation carried out in \sectref{3.3}\ shows that
$u_x/C_x$ is harmonic for $x \neq a,b$
with $u_a/C_a = 1/C_a$
and $b_b/C_b = 0$.
But the function $v_x/C_a$
also has these properties,
so by the Uniqueness Principle
\[
\frac{u_x}{C_x} = \frac{v_x}{C_a}
.
\]
Now
\[
u_{rs} 
=
u_r P_{rs}
=
u_r \frac{C_{rs}}{C_r}
=
v_r \frac{C_{rs}}{C_a}
\]
and
\[
u_{sr}
=
u_s P_{sr}
=
u_s \frac{C_{sr}}{C_s}
=
v_s \frac{C_{sr}}{C_a}
.
\]
Since $C_{rs}=C_{sr}$,
and since by assumption $v_r \geq v_s$,
this means that 
$u_{rs} \geq u_{sr}$.
Therefore,
we see that any edge leads the walker
more often to the more favorable of the points of the edge.

Now let's go back and think about the graph with the bridge.
The above discussion shows that the bridge helps in the sense that,
on the average,
the bridge is crossed more often in the direction that improves 
the chance of escaping.
While this fact is suggestive,
it doesn't tell us that we are more likely to escape
than if the bridge weren't there;
it only tells us what goes on
once the bridge is in place.
What we need is to make a ``before and after'' comparison.

Recall that we are denoting the conductance of the bridge by $\epsilon$.
To distinguish the quantities pertaining to the walks
with and without the bridge,
we will put $(\epsilon)$ superscripts
on quantities that refer to the walk with the bridge,
so that, e.g., $\peps_\esc$ denotes the escape probability 
with the bridge.

Now let $\eps{d}$
denote the expected net number of times the walker crosses
from $r$ to $s$.
As above, we have
\[
\eps{d}
=
\eps{u}_r \frac{\epsilon}{C_r+\epsilon}
- \eps{u}_s \frac{\epsilon}{C_s+\epsilon}
=
\braces{
\frac{\eps{u}_r}{C_r+\epsilon} - \frac{\eps{u}_s}{C_s+\epsilon}
} \epsilon
.
\]

\tout{Claim.}
\[
\peps_\esc = p_\esc + (v_r - v_s) \eps{d}
.
\]
\toutend

\whystart
Every time you use the bridge to go from $r$ to $s$,
you improve your chances of escaping by
\[
(1-v_s) - (1-v_r) = v_r - v_s
\]
assuming that you would continue your walk without using the bridge.
To get the probability of escaping with the bridge,
you take the probability of escaping without the bridge,
and correct it by adding in the change attributable to the bridge,
which is the difference in the original escape probabilities
at the ends of the bridge,
multiplied by the net number of times you expect to cross the bridge.
\whyend

\proofstart
Suppose you're playing a game
where you walk around the graph with the bridge,
and your fortune when you're at $x$ is $v_x$,
which is the probability
that you would return to $a$ before reaching $b$
if the bridge weren't there.
You start at $a$, and walk until you reach $b$ or return to $a$.

This is not a fair game.
Your initial fortune is 1 since you start at $a$ and $v_a = 1$.
Your expected final fortune is
\[
1 \cdot (1-\peps_\esc) + 0 \cdot \peps_\esc
=
1- \peps_\esc
.
\]
The amount you expect to lose by participating in the game is
\[
\peps_\esc
.
\]
(Note that escaping has suddenly become a bad thing!)

Let's analyze where it is that you expect to lose money.
First of all, you lose money when you take the first step away from $a$.
The amount you expect to lose is
\[
1 - \sum_x \Peps_{ax} v_x
=
p_\esc
.
\]

Now if your fortune were given by $\eps{v}_x$ instead of $v_x$,
the game would be fair after this first step.
However, the function $v_x$ is not harmonic
for the walk with the bridge;
it fails to be harmonic at $r$ and $s$.
Every time you step away from $r$,
you expect to lose an amount
\[
v_r
-
\braces{
\sum_x \frac{C_{rs}}{C_r+\epsilon} v_x + \frac{\epsilon}{C_r+\epsilon}v_s
}
=
(v_r-v_s) \frac{\epsilon}{C_r+\epsilon}
.
\]
Similarly, every time you step away from $s$ you expect to lose an amount
\[
(v_s - v_r) \frac{\epsilon}{C_s+\epsilon}
.
\]

The total amount you expect to lose by participating in the game is:
\begin{eqnarray*}
&&
\mbox{expected loss at first step}
+
\\&&
(\mbox{expected loss at $r$}) \cdot (\mbox{expected number of times at $r$})
+
\\&&
(\mbox{expected loss at $s$}) \cdot (\mbox{expected number of times at $s$})
\\&=&
p_\esc
+
\\&&
(v_r - v_s) \frac{\epsilon}{C_r+\epsilon} \eps{u}_r
+
\\&&
(v_s - v_r) \frac{\epsilon}{C_s+\epsilon} \eps{u}_s
\\&=&
p_\esc + (v_r - v_s) \eps{d}
.
\end{eqnarray*}
Equating this with our first expression for the expected cost of 
playing the game yields the formula we were trying to prove.

According to the formula just established,
\[
\peps_\esc - p_\esc
=
(v_r - v_s)
\braces{ \frac{\eps{u}_r}{C_r+\epsilon} - \frac{\eps{u}_s}{C_s+\epsilon} }
\epsilon
.
\]
For small $\epsilon$,
we will have
\[
 \frac{\eps{u}_r}{C_r+\epsilon} - \frac{\eps{u}_s}{C_s+\epsilon}
\approx
\frac{u_r}{C_r} - \frac{u_s}{C_s}   
=
\frac{v_r}{C_a} - \frac{v_s}{C_a}
,
\]
so for small $\epsilon$
\[
\peps_\esc - p_\esc
\approx
(v_r - v_s)^2 \frac{\epsilon}{C_a}
.
\]

This approximation allows us to conclude that
\[
\peps_\esc \geq p_\esc \geq 0
\]
for small $\epsilon$.
But this is enough to establish the monotonicity law,
since any finite change in $\epsilon$
can be realized by making an infinite chain of graphs
each of which is obtained from the one before
by adding a bridge of infinitesimal conductance.

To recapitulate,
the difference in the escape probabilities with and without the bridge
is obtained by taking the difference between
the original escape probabilities at the ends of the bridge,
and multiplying by the expected net number of crossings of the bridge.
This quantity is positive because the walker tends to cross the bridge
more often in the good direction than in the bad direction.

\exstart
\ex{4.2.1}{
Give a probabilistic argument to show that
$C_a p_\esc$ increases with $C_{ar}$ for any $r$.
Give an example to show that $p_\esc$ by itself may actually decrease.
}

\ex{4.2.2}{
Give a probabilistic argument to show that
$C_a p_\esc$ increases with $C_{rb}$ for any $r$.
}

\ex{4.2.3}{
Show that when $v_r = v_s$,
changing the value of $C_{rs}$
does not change $p_\esc$.
}

\ex{4.2.4}{
Show that
\[
\frac{\partial}{\partial R_{rs}} \Reff 
=
i_{rs}^2
.
\]
}

\ex{4.2.5}{
In this exercise we ask you to derive an exact formula
for the change in escape probability
\[
\peps_\esc - p_\esc
,
\]
in terms of quantities that refer only to the walk without the bridge.

(a) Let $N_{xy}$ denote the expected number of times in state $y$
for a walker who starts at $x$
and walks around the graph without the bridge
until he reaches $a$ or $b$.
It is a fact that
\[
\eps{u}_r
=
u_r
+ \eps{u}_r \frac{\epsilon}{C_r+\epsilon} (N_{sr} + 1 - N_{rr})
+ \eps{u}_s \frac{\epsilon}{C_s+\epsilon} (N_{rr} - N_{sr})
.
\]
Explain in words why this formula is true.

(b) This equation for $\eps{u}_r$ can be rewritten as follows:
\[
\frac{C_r}{C_r+\epsilon} \eps{u}_r
=
u_r + \eps{d} (N_{sr} - N_{rr})
.
\]
Prove this formula.
(Hint: Consider a game where your fortune at $x$ is $N_{xr}$,
and where you start from $a$ and walk on the graph with the bridge 
until you reach $b$ or return to $a$.)

(c) Write down the corresponding formula for $\eps{u}_s$,
and use this formula to get an expression for $\eps{d}$
in terms of quantities that refer to the walk without the bridge.

(d) Use the expression for $\eps{d}$ to express
$\peps_\esc - p_\esc$
in terms of quantities that refer to the walk without the bridge,
and verify that the value of your expression is
$\geq 0$ for $\epsilon \geq 0$.
}

\ex{4.2.6}{
Give a probabilistic interpretation of the energy dissipation rate.
}
\exend

\sect{4.3}{A Markov chain proof of the Monotonicity Law}

Let $\matP$ be the ergodic Markov chain associated with an electric network.
When we add an $\epsilon$ bridge from $r$ to $s$,
we obtain a new transition matrix $\mat{\Peps}$
that differs from $\matP$ only for transitions from $r$ and $s$.
We can minimize the differences between $\matP$ and $\mat{\Peps}$
by replacing $\matP$ by the matrix $\mat{\hat{P}}$
corresponding to the circuit without the bridge
but with an $\epsilon$ conductance added from $r$ to $r$ and from $s$ to $s$.
This allows the chain to stay in states $r$ and $s$
but does not change the escape probability from $a$ to $b$.
Thus, we can compare the escape probabilities
for the two matrices $\mat{\hat{P}}$ and $\mat{\Peps}$,
which differ only by
\[
\begin{array}{ll}
\hat{P}_{rr} = \frac{\epsilon}{C_r + \epsilon}
&
\Peps_{rr} = 0
\\
\hat{P}_{rs}=\frac{C_{rs}}{C_r+\epsilon}
&
\Peps_{rs} = \frac{C_{rs}+\epsilon}{C_r+\epsilon}
\\
\hat{P}_{ss} = \frac{\epsilon}{C_s+\epsilon}
&
\Peps_{ss} = 0
\\
\hat{P}_{sr} = \frac{C_{sr}}{C_S+\epsilon}
&
\Peps_{sr} = \frac{C_{sr}+\epsilon}{C_s+\epsilon}
\end{array}
.
\]

We make states $a$ and $b$ into absorbing states.
Let $\mat{\hat{N}}$ and $\mat{\eps{N}}$ be the fundamental matrices
for the absorbing chains
obtained from $\mat{\hat{P}}$ and $\mat{\Peps}$ respectively.
Then
$\mat{\hat{N}} = (\mat{I} - \mat{\hat{Q}})^{-1}$
and
$\mat{\eps{N}} = (\mat{I} - \mat{\eps{Q}})^{-1}$
where $\mat{\hat{Q}}$ and $\mat{\eps{Q}}$
differ only for the four components
involving only $r$ and $s$.
That is,
\[
\mat{\eps{Q}} =
\mat{\hat{Q}} +
\bordermatrix{
 &r&s \cr
r&  -\frac{\epsilon}{C_r+\epsilon} & \frac{\epsilon}{C_r+\epsilon} \cr
s&  \frac{\epsilon}{C_s+\epsilon} & -\frac{\epsilon}{C_s+\epsilon} }
= \mat{\hat{Q}} + \mat{h}\mat{k}
\]
where $\mat{h}$ is the column vector with only components
$r$ and $s$ non-zero
\[
\mat{h} =
\bordermatrix{
 &\ \cr
r&  \frac{\epsilon}{C_r+\epsilon} \cr
s& -\frac{\epsilon}{C_s+\epsilon} }
\]
and $\mat{k}$ is a row vector with only components $r$ and $s$ non-zero 
\[
\mat{k} =
\bordermatrix{
 &r&s \cr
 &-1&1 }
.
\]

J. G. Kemeny has pointed out to us
that if $\mat{A}$ is any matrix with inverse $\mat{N}$
and we add to $\mat{A}$ a matrix of the form $- \mat{h} \mat{k}$,
then $\mat{A} - \mat{h} \mat{k}$
has an inverse if and only if
$\mat{k} \mat{N} \mat{h} \neq 1$
and, if so,
$\mat{\bar{N}} = (\mat{A} - \mat{h} \mat{k})^{-1}$
is given by
\[
\mat{\bar{N}} = \mat{N} + c (\mat{N} \mat{h}) (\mat{k} \mat{N})
\]
where $c = 1/(1- \mat{k} \mat{N} \mat{h})$.
You are asked to prove this in Exercise \exref{4.3.1}.
Adding $- \mat{h} \mat{k}$
to $\mat{A} = \mat{I} - \mat{\hat{Q}}$ and using this result,
we obtain
\[
\mat{\eps{N}}
= \mat{\hat{N}} + c (\mat{\hat{N}} \mat{h}) (\mat{k} \mat{\hat{N}})
.
\]

Using the simple nature of $\mat{h}$ and $\mat{k}$ we obtain
\[
\eps{N}_{ij}
=
\hat{N}_{ij} +
\braces{
\frac{\hat{N}_{ir}\epsilon}{C_r+\epsilon}
- \frac{\hat{N}_{is}\epsilon}{C_s+\epsilon}
}
(\hat{N}_{sj} - \hat{N}_{rj})
\]
and
\[
c =
\recip{
1 +
\frac{\hat{N}_{rr}\epsilon}{C_r+\epsilon}
-
\frac{\hat{N}_{sr}\epsilon}{C_r+\epsilon}
+
\frac{\hat{N}_{ss}\epsilon}{C_s+\epsilon}
-
\frac{\hat{N}_{rs}\epsilon}{C_s+\epsilon}
}
.
\]

Since $\hat{N}_{rr}$
 is the expected number of times in $r$
starting in $r$ and
$\hat{N}_{sr}$ is the expected number of times in $r$ starting in $s$,
$\hat{N}_{rr} \geq \hat{N}_{sr}$.
Similarly
$\hat{N}_{ss} \geq \hat{N}_{rs}$
and so the denominator of $c$ is $\geq 1$.
In particular, it is positive.

Recall that the absorption probabilities for state $b$
are given by
\[
B_{xb}
=
\sum_y
N_{xy} P_{yb}
.
\]
Since
$\Peps_{xb} = \hat{P}_{xb}$,
\[
\eps{B}_{xb}
=
\hat{B}_{xb}
+ c
\braces{
\frac{\hat{N}_{xr}\epsilon}{C_r+\epsilon} 
-
\frac{\hat{N}_{xs}\epsilon}{C_s+\epsilon}
}
(\hat{B}_{sb} - \hat{B}_{rb})
.
\]
Since
$\Peps_{ax} = \hat{P}_{ax}$,
\[
\peps_\esc
=
\hat{p}_\esc
+ c
\braces{
\frac{\hat{u}_r \epsilon}{C_r+\epsilon}
-
\frac{\hat{u}_s \epsilon}{C_s+\epsilon}
}
(\hat{B}_{sb} - \hat{B}_{rb})
\]
where $\hat{u}_x$ is 
the expected number of times that the ergodic chain $\hat{P}$,
started at state $a$,
is in state $x$ before returning to $a$
reaching $b$ for the first time.
The absorption probability $B_{xa}$ is the quantity $v_x$
introduced in the previous section.
As shown there,
reversibility allows us to conclude that
\[
\frac{\hat{u}_x}{\hat{C}_x}
=
\frac{\hat{B}_{xa}}{\hat{C}_a}
=
\frac{\hat{B}_{xa}}{C_a}
\]
so that
\[
\peps_\esc
=
p_\esc
+
\frac{\epsilon c}{C_a}
(\hat{B}_{sb} - \hat{B}_{rb})^2
\]
and this shows that the Monotonicity Law is true.

The change from $\matP$ to $\mat{\hat{P}}$
was merely to make the calculations easier.
As we have remarked,
the escape probabilities are the same for the two chains
as are the absorption probabilities $B_{ib}$.
Thus we can remove the hats and write the same formula.
\[
\peps_\esc
=
p_\esc
+
\frac{\epsilon c}{C_a}
(B_{sb} - B_{rb})^2
.
\]

The only quantity in this final expression that seems to depend 
upon quantities from $\mat{\hat{P}}$ is $c$.
In Exercise \exref{4.3.2}\ you are asked
to show that $c$ can also
be expressed in terms of the fundamental matrix $\mat{N}$
obtained from the original $\matP$.

\exstart
\ex{4.3.1}{
Let $\mat{A}$ be a matrix with inverse $\mat{N} = \mat{A}^{-1}$
Let $\mat{h}$ be a column vector
and $\mat{k}$ a row vector.
Show that
\[
\mat{\bar{N}} = (\mat{A} - \mat{h} \mat{k})^{-1}
\]
exists if and only if 
$\mat{k} \mat{N} \mat{h} \neq 1$
and, if so,
\[
\mat{\bar{N}} =
\mat{N} +
\frac{(\mat{N} \mat{h})(\mat{k} \mat{N})}{1-\mat{k} \mat{N} \mat{h}}
.
\]
}

\ex{4.3.2}{
Show that $c$ can be expressed in terms of the fundamental matrix $\mat{N}$
of the original Markov chain $\matP$ by
\[
c =
\recip{
1 +
\frac{N_{rr}\epsilon}{C_r}
-
\frac{N_{sr}\epsilon}{C_r}
+
\frac{N_{ss}\epsilon}{C_s}
-
\frac{N_{rs}\epsilon}{C_s}
}
.
\]
}
\exend

\prt{II}{Random walks on infinite networks}

\chap{5}{Polya's recurrence problem}

\sect{5.1}{Random walks on lattices}

In 1921 George Polya \cite{polya:irrfahrt}
investigated random walks on certain infinite graphs,
or as he called them, ``street networks''.
The graphs he considered,
which we will refer to as lattices,
are illustrated in Figure \figref{5.1}.
\prefig{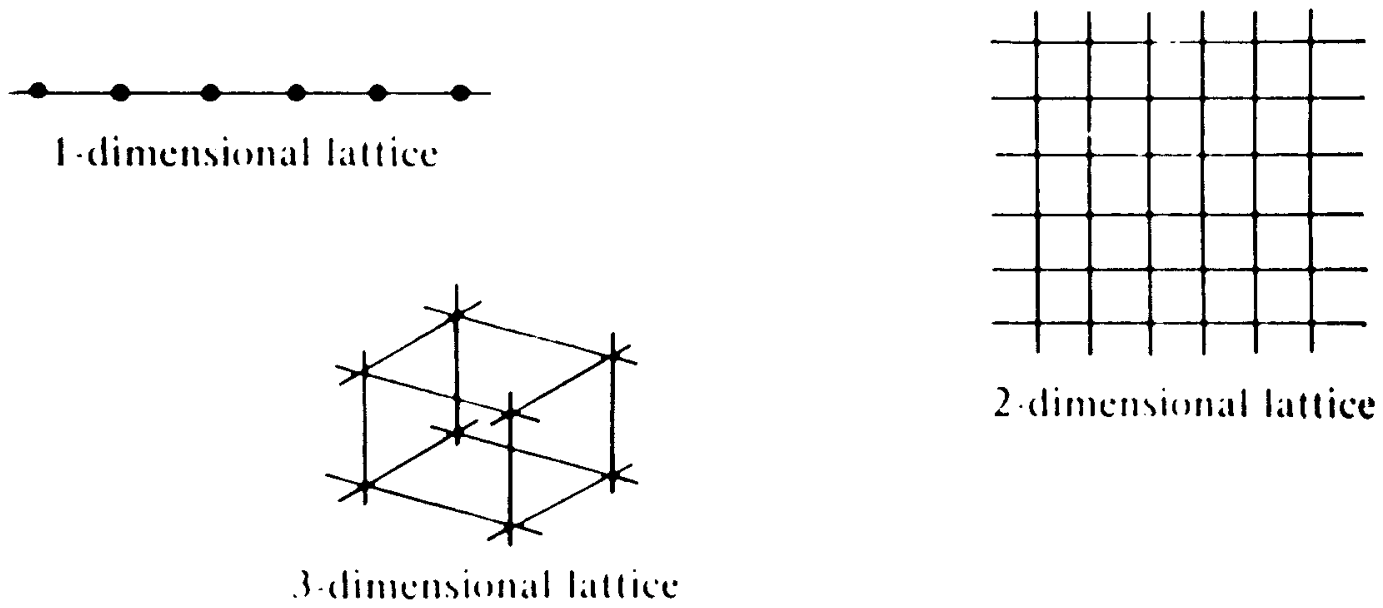}

To construct a $d$-dimensional lattice,
we take as vertices those points
$(x_1,\ldots,x_d)$ of $\R^d$
all of whose coordinates are integers,
and we join each vertex by an undirected line segment
to each of its $2d$ nearest neighbors.
These connecting segments,
which represent the edges of our graph,
each have unit length
and run parallel to one of the coordinate axes of $\R^d$.
We will denote this $d$-dimensional lattice by $\Z^d$.
We will denote the origin
$(0, 0, \ldots , 0)$ by $\orig$.

Now let a point walk around at random on this lattice.
As usual, by walking at random we mean that,
upon reaching any vertex of the graph,
the probability of choosing any one of the $2d$ edges leading out 
of that vertex is $\recip{2d}$.
We will call this random walk \dt{simple random walk}
in $d$ dimensions.

When $d = 1$,
our lattice is just an infinite line
divided into segments of length one.
We may think of the vertices of this graph
as representing the fortune of a gambler
betting on heads or tails in a fair coin tossing game.
Simple random walk in one dimension then represents the 
vicissitudes of his or her fortune,
either increasing or decreasing by one unit after each round of the game.

When $d =2$,
our lattice looks like an infinite network of streets and avenues,
which is why we describe the random motion of the 
wandering point as a ``walk''.

When $d = 3$,
the lattice looks like an infinite ``jungle gym'',
so perhaps in this case we ought to talk about a ``random climb'',
but we will not do so.
It is worth noting that when $d = 3$,
the wanderings of our point
could be regarded as an approximate representation of the random path 
of a molecule diffusing in an infinite cubical crystal.
Figure \figref{5.2} shows
a simulation of a simple random walk in three dimensions.
\prefig{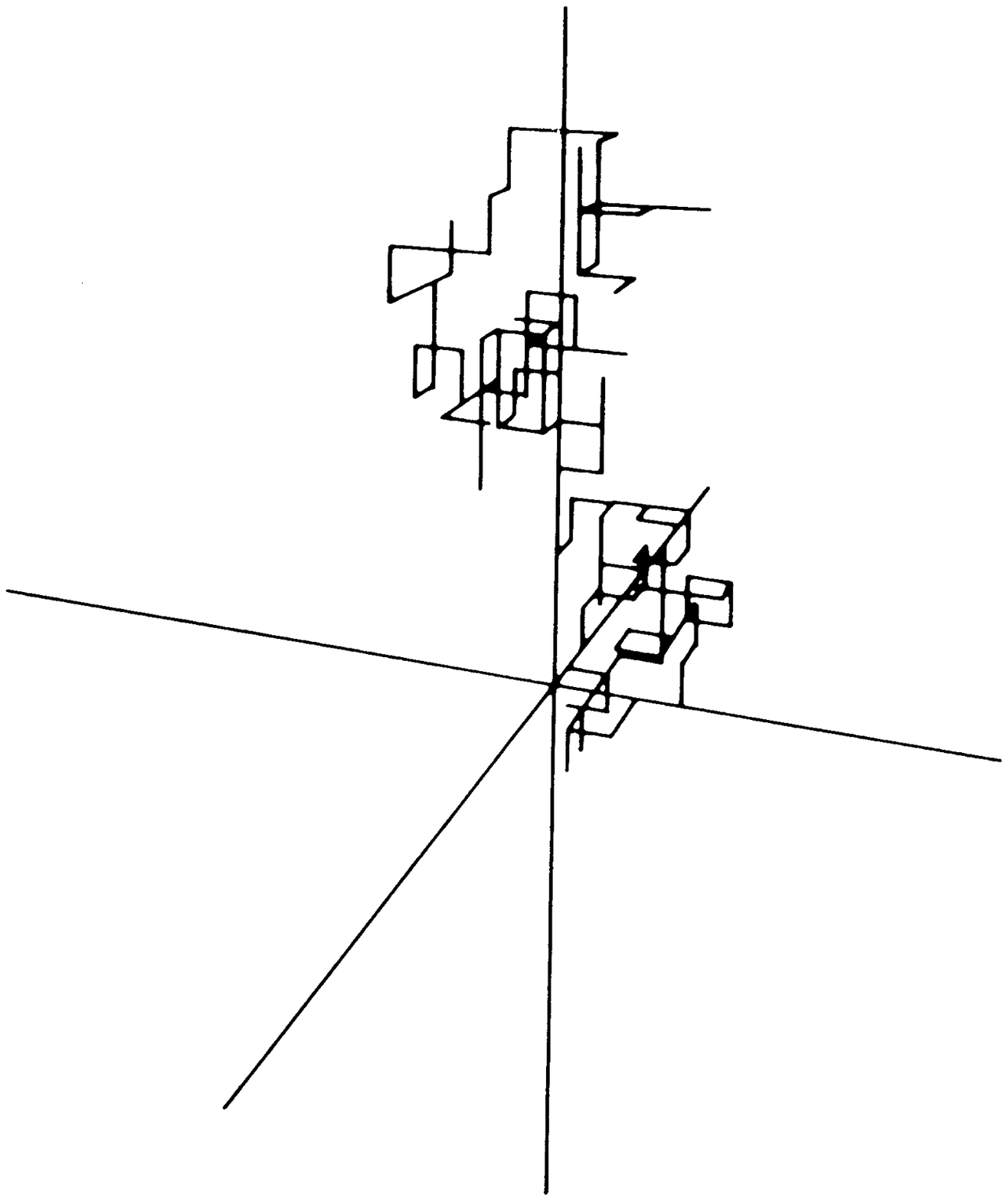}

\sect{5.2}{The question of recurrence}

The question that Polya posed amounts to this:
``Is the wandering point certain to return to its starting point
during the course of its wanderings?''
If so, we say that the walk is \dt{recurrent}.
If not, that is, if there is a positive probability 
that the point will never return to its starting point,
then we say that the walk is \dt{transient}.

If we denote the probability that the point never returns to its 
starting point by $p_\esc$,
then the chain is
recurrent if $p_\esc = 0$,
and
transient if $p_\esc > 0$.

We will call
the problem of determining recurrence or transience of a random walk
the \dt{type problem}.

\sect{5.3}{Polya's original question}

The definition of recurrence that we have given differs from Polya's
original definition.
Polya defined a walk to be recurrent if, 
with probability one,
it will pass through every single point of the lattice
in the course of its wanderings.
In our definition, we require only that the point
return to its starting point.
So we have to ask ourselves,
``Can the random walk be recurrent in our sense
and fail to be recurrent in Polya's sense?''

The answer to this question is,
``No, the two definitions of recurrence are equivalent.''
Why?
Because if the point must return once to its starting point,
then it must return there again and again,
and each time it starts away from the origin,
it has a certain non-zero probability of
hitting a specified target vertex
before returning to the origin.
And anyone can get a bull's-eye
if he or she is allowed an infinite number of darts,
so eventually the point will hit the target vertex.

\exstart
\ex{5.3.1}{
Write out a rigorous version of the argument just given.
}
\exend

\sect{5.4}{Polya's Theorem: recurrence in the plane, transience in space}

In \cite{polya:irrfahrt},
Polya proved the following theorem:

\tout{Polya's Theorem.}
Simple random walk on a $d$-dimensional lattice
is recurrent for $d = 1, 2$
and transient for $d > 2$.
\toutend

The rest of this \book\ will be devoted to
trying to understand this theorem.
Our approach will be to exploit once more
the connections between questions about a random walk on a graph
and questions about electric currents
in a corresponding network of resistors.
We hope that this approach,
by calling on methods that appeal to our physical intuition,
will leave us feeling that we know ``why'' the theorem is true.

\exstart
\ex{5.4.1}{
Show that Polya's theorem implies that
if two random walkers start at $\orig$
and wander independently,
then in one and two dimensions they will eventually meet again,
but in three dimensions there is positive probability that they won't.
}

\ex{5.4.2}{
Show that Polya's theorem implies that
a random walker in three dimensions
will eventually hit the line defined by $x=2,z=0$.
}
\exend

\sect{5.5}{The escape probability as a limit of escape probabilities for
finite graphs}

We can determine the type of an infinite lattice from properties 
of bigger and bigger finite graphs that sit inside it.
The simplest way to go about this
is to look at the lattice analog of balls (solid spheres) in space.
These are defined as follows:
Let $r$ be an integer---this will be the radius of the ball.
Let $\ball{G}{r}$ be the graph gotten from $\Z^d$ by throwing out vertices 
whose distance from the origin is $>r$.
By ``distance from the origin''
we mean here not the usual Euclidean distance,
but the distance ``in the lattice'';
that is, the length of the shortest path
along the edges of the lattice between the two points.
Let $\sphere{G}{r}$ be the ``sphere'' of radius $r$ about the origin,
i.e., those points that are exactly $r$ units from the origin.
In two dimensions, $\sphere{G}{r}$ looks like a square.
(See Figure \figref{5.3}.)
\prefig{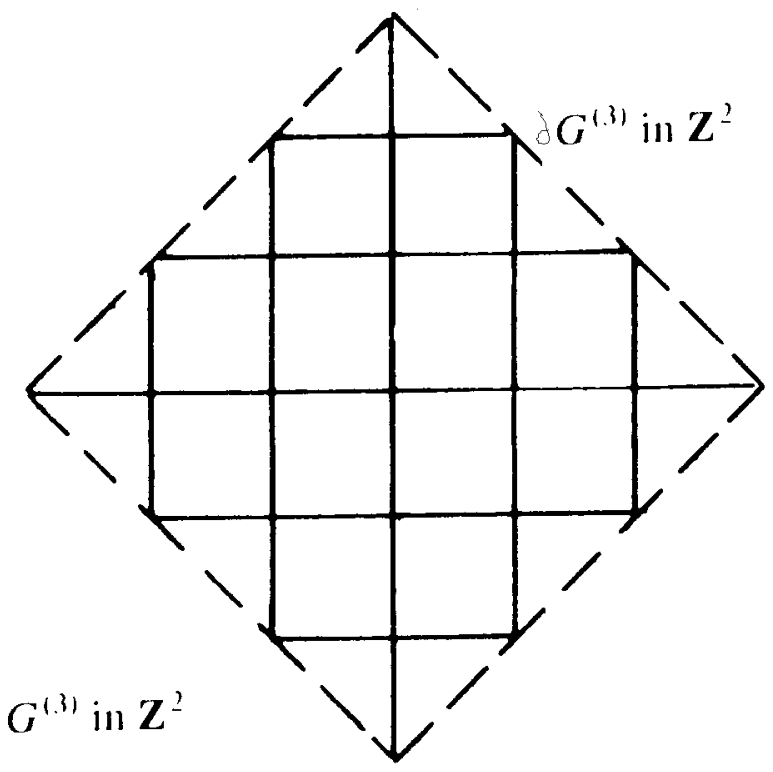}
In three dimensions, it looks like an octahedron.

We define a random walk on $\ball{G}{r}$ as follows:
The walk starts at $\orig$
and continues as it would on $\Z^d$
until it reaches a point on $\sphere{G}{r}$
and then it stays at this point.
Thus the walk on $\ball{G}{r}$
is an absorbing Markov chain
with every point of $\sphere{G}{r}$ an absorbing state.

Let $p^{(r)}_\esc$ be the probability that a random walk on $\ball{G}{r}$
starting at $\orig$,
reaches $\sphere{G}{r}$ before returning to $\orig$.
Then $p^{(r)}_\esc$ decreases as $r$ increases and
$p_\esc = \lim_{r \goto \infty} p^{(r)}_\esc$
is the \dt{escape probability}
for the infinite graph.
If this limit is 0, the infinite walk is recurrent.
If it is greater than 0, the walk is transient.

\sect{5.6}{Electrical formulation of the type problem}

Now that we have expressed things in terms of finite graphs,
we can make use of electrical methods.
To determine $p_\esc$ electrically,
we simply ground all the points of $\sphere{G}{r}$,
maintain $\orig$ at one volt,
and measure the current $i^{(r)}$ flowing into the circuit.
(See Figure \figref{5.4}.)
\prefig{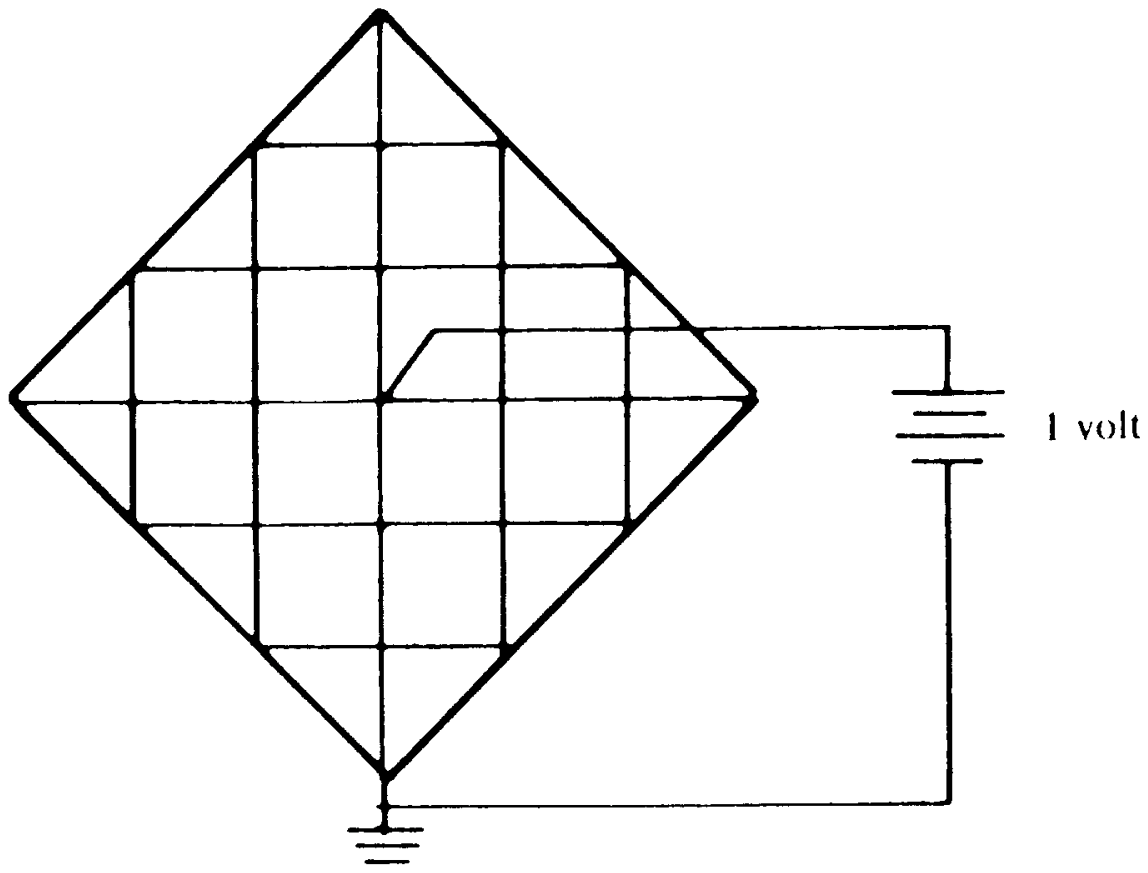}

From \sectref{3.4},
we have
\[
p^{(r)}_\esc = \frac{i^{(r)}}{2d}
,
\]
where $d$ is the dimension of the lattice.
(Remember that we have to divide
by the number of branches coming out of the starting point.) 
Since the voltage being applied is 1,
$i^{(r)}$ is just the effective conductance
between $\orig$ and $\sphere{G}{r}$,
i.e.,
\[
i^{(r)} = \recip{R^{(r)}_\eff}
.
\]
where $R^{(r)}_\eff$ is the effective resistance
from $\orig$ to $\sphere{G}{r}$.
Thus
\[
p^{(r)}_\esc = \recip{2 d R^{(r)}_\eff}
.
\]

Define $\Reff$,
the \dt{effective resistance from the origin to infinity},
to be
\[
\Reff = \lim_{r \goto \infty} R^{(r)}_\eff
.
\]
This limit exists since $R^{(r)}_\eff$ is an increasing function of $r$.
Then
\[
p_\esc = \recip{2 d \Reff}
.
\]
Of course $\Reff$ may be infinite;
in fact, this will be the case if and only if $p_\esc=0$.
Thus the walk is recurrent
if and only if the resistance to infinity is infinite,
which makes sense.
\prefig{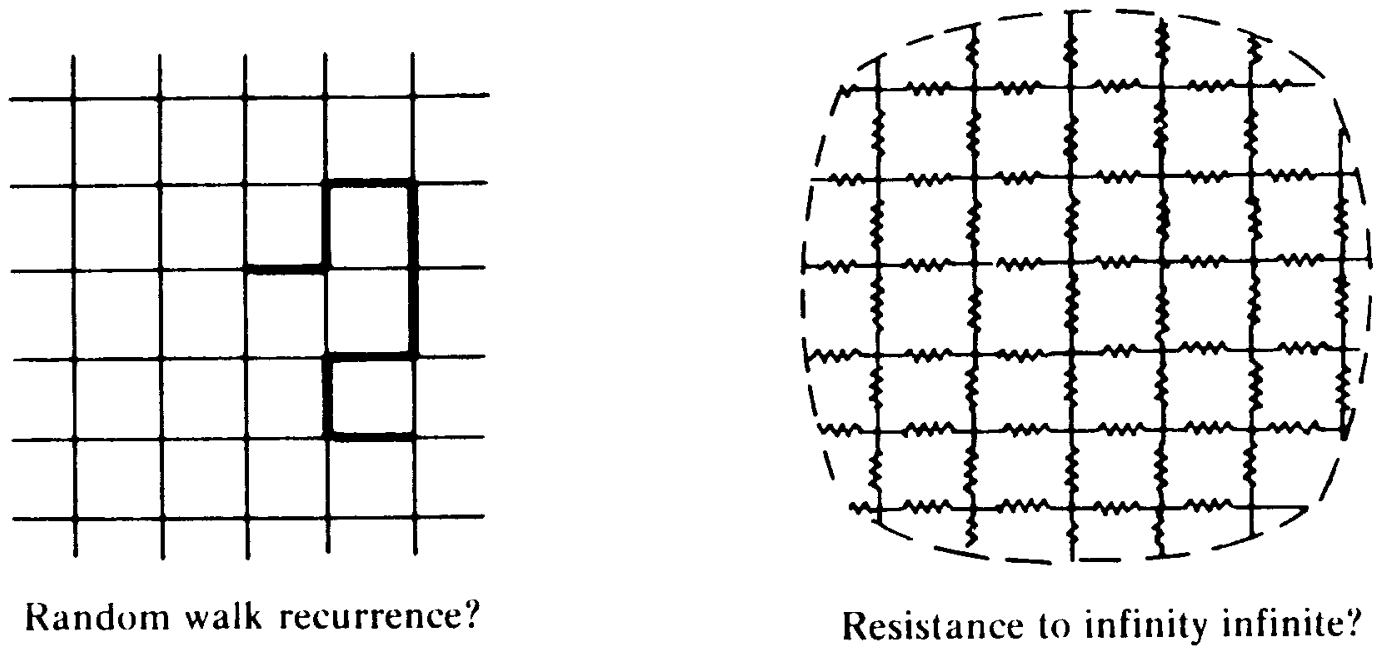}

The success of this electrical formulation of the type problem 
will come from the fact that
the resistance to infinity can be estimated
using classical methods of electrical theory.

\sect{5.7}{One Dimension is easy, but what about higher dimensions?}

We now know that simple random walk on a graph is recurrent
if and only if
a corresponding network of 1-ohm resistors
has infinite resistance ``out to infinity''.
Since an infinite line of resistors
obviously has infinite resistance,
it follows that simple random walk on the 1-dimensional lattice
is recurrent, as stated by Polya's theorem.
\prefig{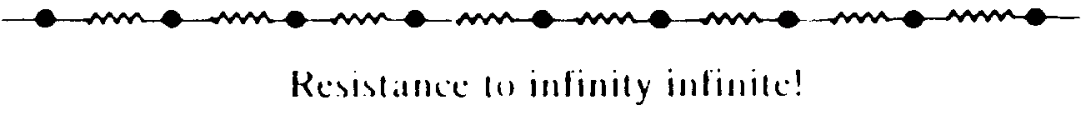}

What happens in higher dimensions?
We are asked to decide whether 
a $d$-dimensional lattice has infinite resistance to infinity.
The difficulty is that
the $d$-dimensional lattice $\Z^d$
lacks the rotational symmetry of the Euclidean space $\R^d$ in which it sits.

To see how this lack of symmetry complicates electrical problems, 
we determine,
by solving the appropriate discrete Dirichlet problem, 
the voltages for a one-volt battery
attached between $\orig$ and the points of $\sphere{G}{3}$ in $\Z^2$.
The resulting voltages are:
\[
\matrix{
 & & &0 \cr
 & &0&.091&0 \cr
 &0&.182&.364&.182&0 \cr
0&.091&.364&1&.364&.091&0 \cr
 &0&.182&.364&.182&0 \cr
 & &0&.091&0 \cr
 & & &0 }
.
\]
The voltages at points of $\sphere{G}{1}$ are equal,
but the voltages at points of $\sphere{G}{2}$ are not.
This means that the resistance from $\orig$ to $\sphere{G}{3}$
cannot be written
simply as the sum of the resistances
from $\orig$ to $\sphere{G}{1}$,
$\sphere{G}{1}$ to $\sphere{G}{2}$,
and $\sphere{G}{2}$ to $\sphere{G}{3}$.
This is in marked contrast to
the case of a continuous resistive medium
to be discussed in \sectref{5.8}.

\exstart
\ex{5.7.1}{
Using the voltages given for $\ball{G}{3}$,
find $\ball{R}{3}_\eff$ and $\ball{p}{3}_\esc$.
}

\ex{5.7.2}{
Consider a one-dimensional infinite network
with resistors $R_n = 1/2^n$ from $n$ to $n+1$
for $n= \ldots, -2, -1, 0, 1, 2, \ldots$.
Describe the associated random walk
and determine whether the walk is recurrent or transient.
}

\ex{5.7.3}{
A random walk moves on the non-negative integers;
when it is in state $n$, $n>0$,
it moves with probability $p_n$ to $n+1$
and with probability $1-p_n$, to $n-1$.
When at $\orig$, it moves to 1.
Determine a network that gives this random walk
and give a criterion in terms of the $p_n$
for recurrence of the random walk.
}
\exend

\sect{5.8}{Getting around the lack of rotational symmetry of the lattice}

Suppose we replace our $d$-dimensional resistor lattice
by a (homogeneous, isotropic) resistive medium
filling all of $\R^d$
and ask for the effective resistance to infinity.
Naturally we expect that the rotational symmetry
will make this continuous problem
easier to solve than the original discrete problem.
If we took this problem to a physicist,
he or she would probably produce something like the scribblings
illustrated in Figure \figref{5.7},
and conclude that the effective resistance is infinite
for $d = 1, 2$ and finite for $d>2$.
The analogy to Polya's theorem is obvious,
but is it possible to translate these calculations for continuous media
into information about what happens in the lattice?
\prefig{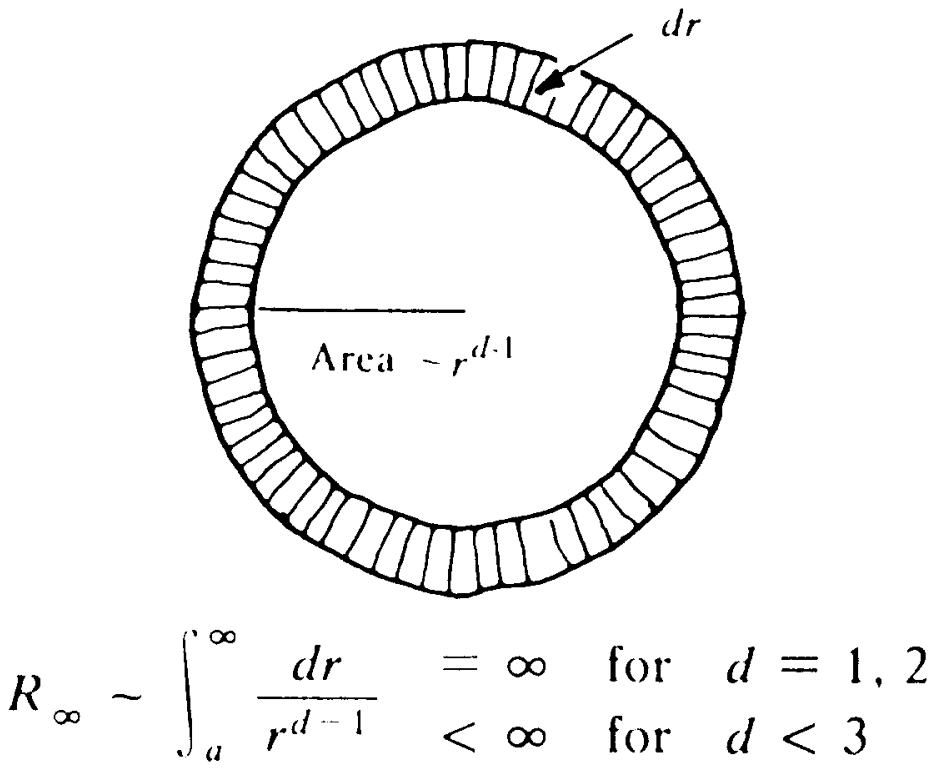}

This can indeed be done,
and this would certainly be the most natural approach to take.
We will come back to this approach at the end of the \book.
For now, we will take a different approach to
getting around the asymmetry of the lattice.
Our method will be to modify the lattice
in such a way as to obtain a graph
that is symmetrical enough
so that we can calculate its resistance out to infinity.
Of course, we will have to think carefully about 
what happens to that resistance when we make these modifications.

\sect{5.9}{Rayleigh: shorting shows recurrence in the plane,
cutting shows transience in space}

Here is a sketch of the method we will use to prove Polya's theorem.

To take care of the case $d = 2$,
we will modify the 2-dimensional resistor network
by shorting certain sets of nodes together
so as to get a new network
whose resistance is readily seen to be infinite.
As shorting can only decrease the effective resistance of the network,
the resistance of the original network must also be infinite.
Thus the walk is recurrent when $d=2$.

To take care of the case $d = 3$,
we will modify the 3-dimensional network
by cutting out certain of the resistors
so as to get a new network
whose resistance is readily seen to be finite.
As cutting can only increase the resistance of the network,
the resistance of the original network must also be finite.
Thus the walk is transient when $d = 3$.

The method of applying shorting and cutting
to get lower and upper bounds for the resistance of a resistive medium
was introduced by Lord Rayleigh.
(Rayleigh \cite{rayleigh};
see also Maxwell \cite{maxwell},
Jeans \cite{jeans}, 
PoIya and Szego \cite{polyaszego}).
We will refer to Rayleigh's techniques collectively as
\dt{Rayleigh's short-cut method}.
This does not do Rayleigh justice,
for Rayleigh's method is a whole bag of tricks
that goes beyond mere shorting and cutting---but who can resist a pun?

Rayleigh's method was apparently first applied to random walks
by C. St. J. A. Nash-Williams \cite{nashwilliams},
who used the shorting method to establish recurrence
for random walks on the 2-dimensional lattice.

\chap{6}{Rayleigh's short-cut method}

\sect{6.1}{Shorting and cutting}

In its simplest form,
Rayleigh's method involves modifying
the network whose resistance we are interested in
so as to get a simpler network.
We consider two kinds of modifications,
shorting and cutting. 
Cutting involves nothing more than
clipping some of the branches of the network, 
or what is the same,
simply deleting them from the network.
Shorting involves
connecting a given set of nodes together
with perfectly conducting wires,
so that current can pass freely between them.
In the resulting network, the nodes that were shorted 
together behave as if they were a single node.
\prefig{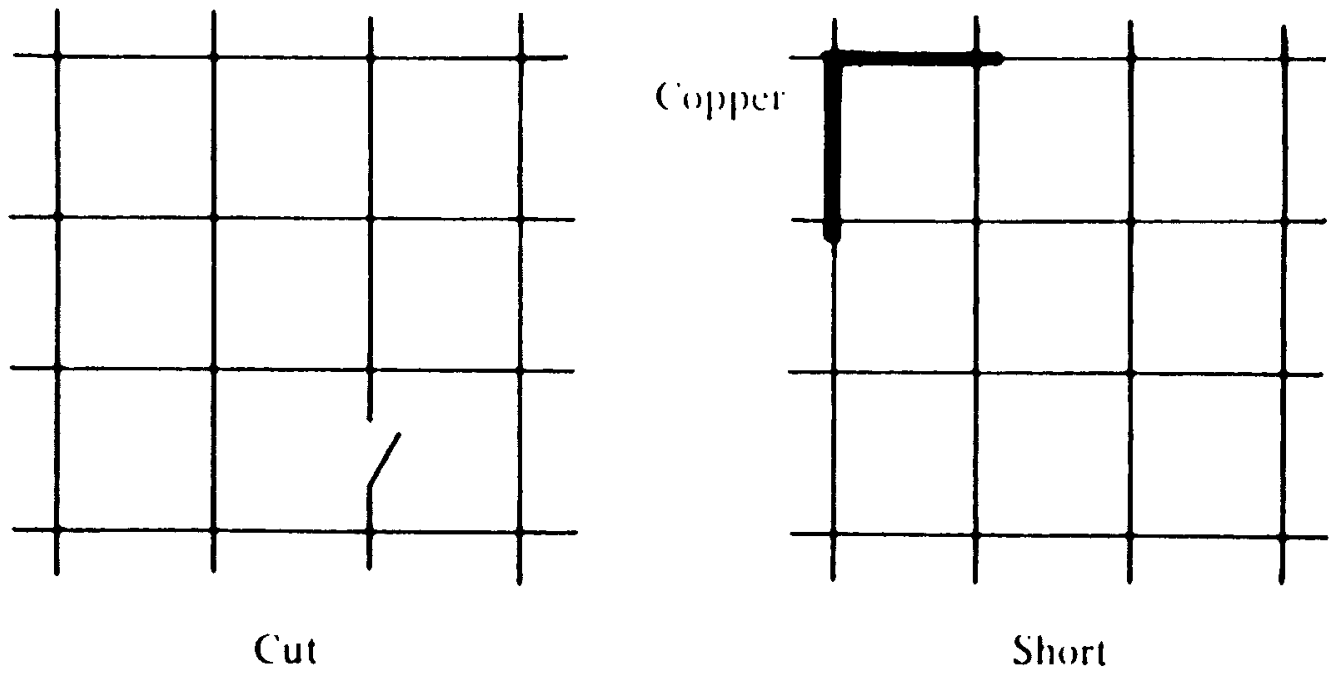}

\sect{6.2}{The Shorting Law and the Cutting Law; Rayleigh's idea}

The usefulness of these two procedures (shorting and cutting) 
stems from the following observations:

\tout{Shorting Law.}
Shorting certain sets of nodes together
can only decrease the effective resistance of the network
between two given nodes.
\toutend

\tout{Cutting Law.}
Cutting certain branches
can only increase the effective resistance
between two given nodes.
\toutend

These laws are both equivalent to Rayleigh's Monotonicity Law, 
which was introduced in \sectref{4.1}\ (see Exercise \exref{6.1.1}):

\tout{Monotonicity Law.}
The effective resistance between two given nodes
is monotonic in the branch resistances.
\toutend

Rayleigh's idea was to use the Shorting Law and the Cutting Law above
to get lower and upper bounds for the resistance of a network.
In \sectref{6.3}\ we apply this method
to solve the recurrence problem 
for simple random walk in dimensions 2 and 3.

\exstart
\ex{6.1.1}{
Show that the Shorting Law and the Cutting Law are both 
equivalent to the Monotonicity Law.
}
\exend

\sect{6.3}{The plane is easy}

When $d =2$,
we apply the Shorting Law as follows:
Short together nodes on squares about the origin,
as shown in Figure \figref{6.2}.
\prefig{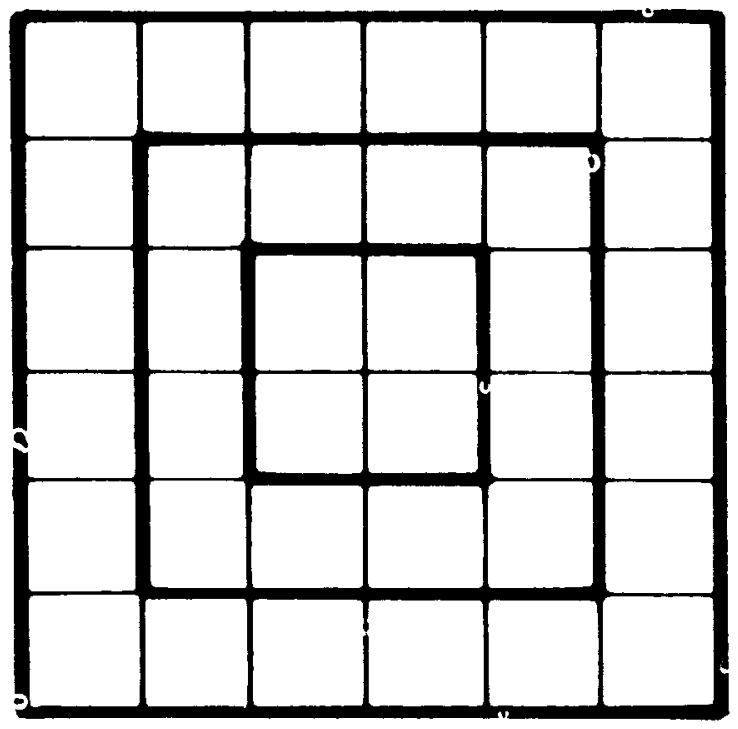}
The network we obtain
is equivalent to the network shown in Figure \figref{6.3}.
\prefig{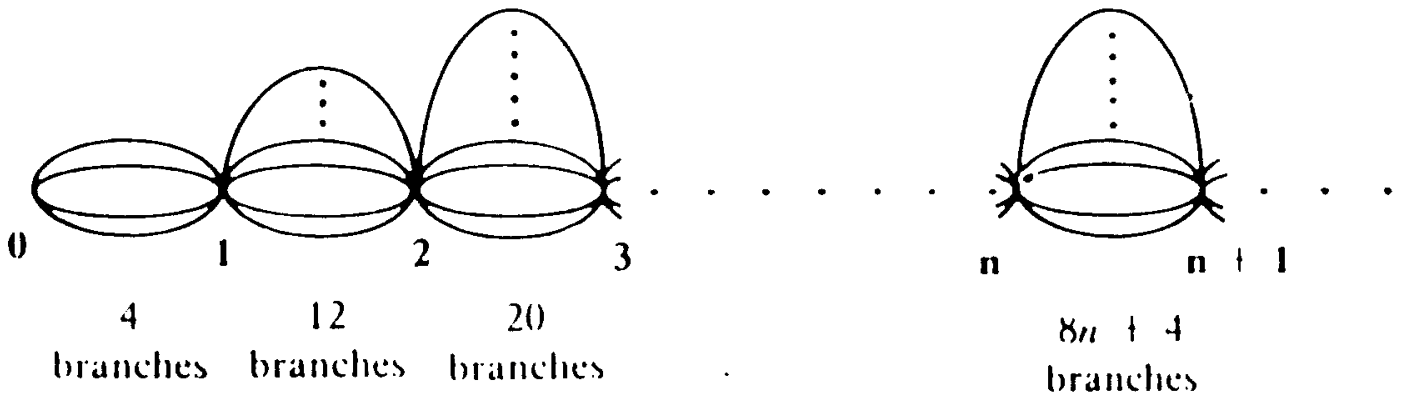}

Now as $n$ 1-ohm resistors in parallel
are equivalent to a single resistor of resistance $\recip{n}$  ohms,
the modified network is equivalent to
the network shown in Figure \figref{6.4}.
\prefig{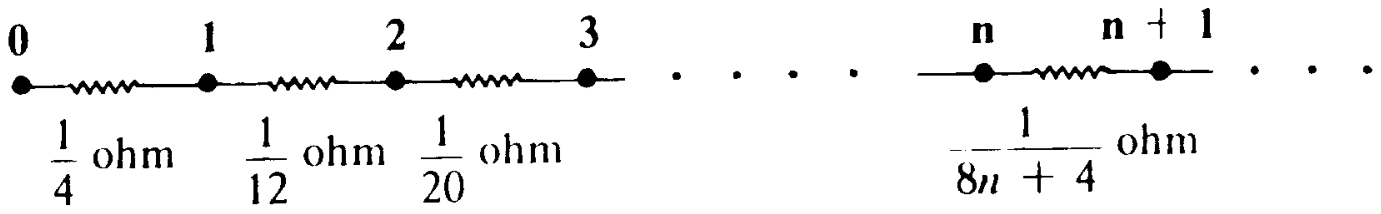}
The resistance of this network out to infinity is
\[
\sum_{n=0}^{\infty} \recip{8n+4} = \infty
.
\]
As the resistance of the old network can only be bigger,
we conclude that it too must be infinite,
so that the walk is recurrent when $d = 2$.

\exstart
\ex{6.3.1}{
Using the shorting technique,
give an upper bound for $\ball{p}{3}_\esc$,
and compare this with the exact value obtained in 
Exercise \exref{5.7.1}.
}
\exend

\sect{6.4}{Space: searching for a residual network}

When $d = 3$,
what we want to do is
delete certain of the branches of the network
so as to leave behind a residual network
having manifestly finite resistance.
The problem is to reconcile the ``manifestly''
with the ``finite''.
We want to cut out enough edges so that
the effective resistance of what is left is easy to calculate,
while leaving behind enough edges so that
the result of the calculation is finite.

\sect{6.5}{Trees are easy to analyze}

Trees---that is,
graphs without circuits---are undoubtedly the easiest to work with.
For instance,
consider the full binary tree,
shown in Figure \figref{6.5}.
\prefig{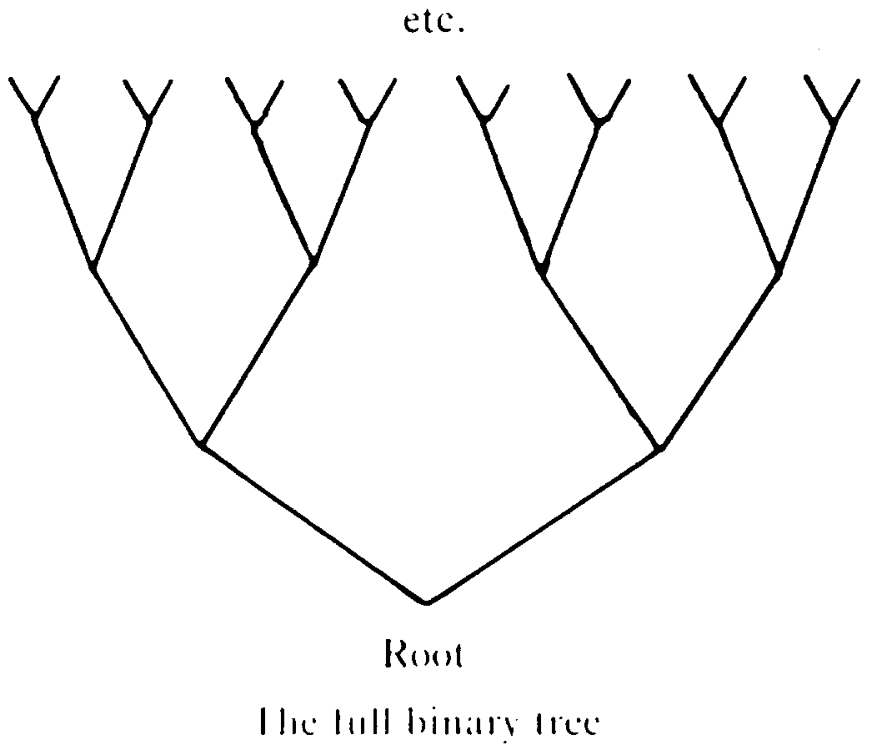}
\prefig{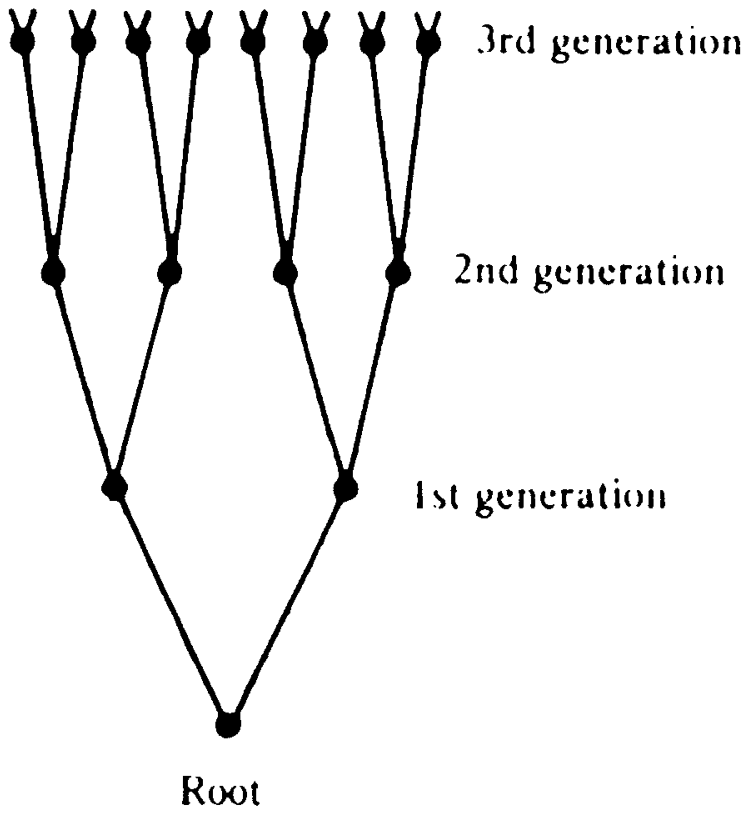}
Notice that sitting inside this tree just above the root
are two copies of the tree itself.
This self-similarity property can be used to compute
the effective resistance $\Rinfty$ from the root out to infinity.
(See Exercise \exref{6.5.1}.)
It turns out that $\Rinfty=1$.
We will demonstrate this below by a more direct method.

To begin with,
let us determine the effective resistance $R_n$
between the root and the set of $n$th generation branch points.
To do this,
we should ground the set of branch points,
hook the root up to a 1-volt battery,
and compute 
\[
R_n = \recip{\mbox{current through battery}}
.
\]
For $n = 3$,
the circuit that we would obtain is shown in Figure \figref{6.7}.
\prefig{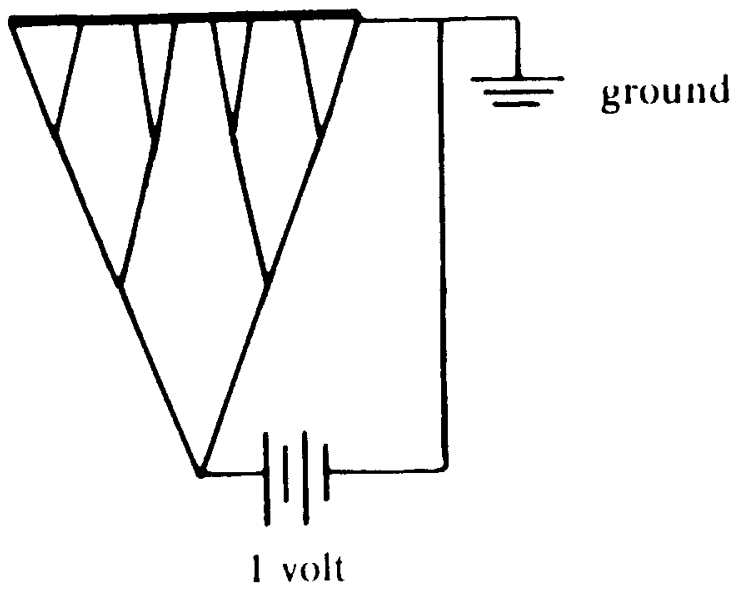}

In the resulting circuit,
all branch points of the same generation
are at the same voltage (by symmetry).
Nothing happens when you short together
nodes that are at the same potential.
Thus shorting together branch points of the same generation
will not affect the
distribution of currents in the branches.
In particular,
this modification will not affect
the current through the battery,
and we conclude that 
\[
R_n
=
\recip{\mbox{current in original circuit}}
=
\recip{\mbox{current in modified circuit}}
.
\]
For $n = 3$,
the modified circuit is shown in Figure \figref{6.8}.
\prefig{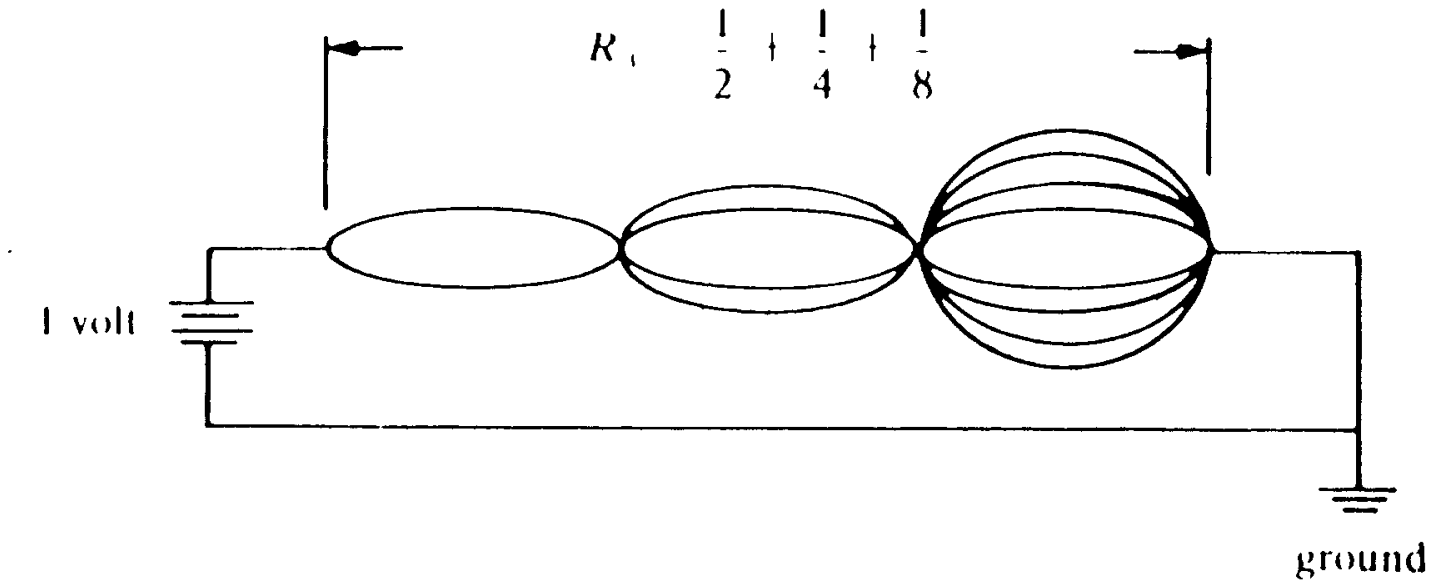}
This picture shows that
\[
R_3 = \recip{2} + \recip{4} + \recip{8} = 1-\recip{2^3}
.
\]
More generally,
\[
R_n =
\recip{2} + \recip{4} + \ldots + \recip{2^n}
=
1-\recip{2^n}
.
\]
Letting $n \goto \infty$,
we get
\[
\Rinfty
=
\lim_{n \goto \infty} R_n
=
\lim_{n \goto \infty} 1-\recip{2^n}
=
1
.
\]

Figure \figref{6.9}\ shows another closely related tree,
the \dt{tree homogeneous of degree three}:
Note that all nodes of this tree are similar---there is no
intrinsic way to distinguish one from another.
\prefig{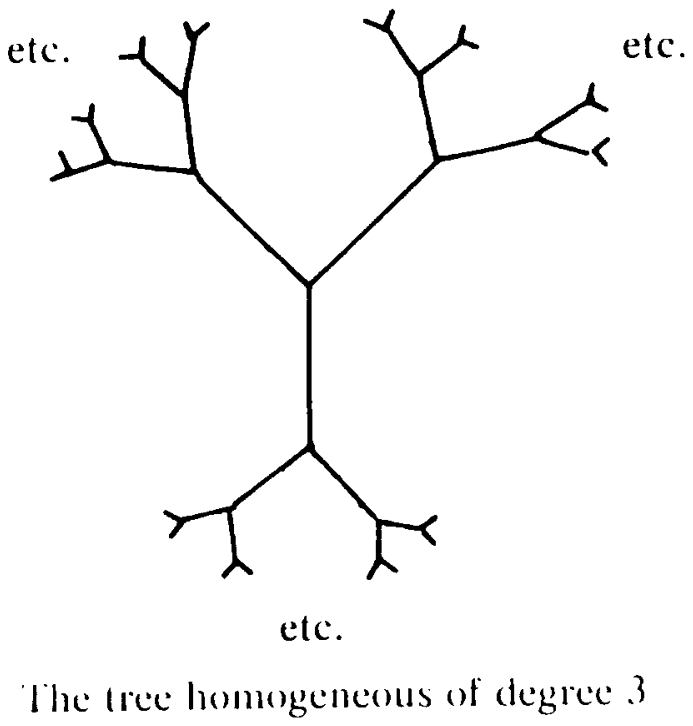}
This tree is obviously a
close relative of the full binary tree.
Its resistance to infinity is $2/3$.

\exstart
\ex{6.5.1}{
(a) Show, using the self-similarity of the full binary tree, 
that the resistance $\Rinfty$ to infinity satisfies the equation
\[
\Rinfty = \frac{\Rinfty+1}{2}
\]
and conclude that either $\Rinfty =1$ or $\Rinfty = \infty$.

(b) Again using the self-similarity of the tree,
show that 
\[
R_{n+1} = \frac{R_n+1}{2}
\]
where $R_n$ denotes the resistance out to
the set of the $n$th-generation branch points.
Conclude that
\[
\Rinfty = \lim_{n \goto \infty} R_n = 1
.
\]
}
\exend

\sect{6.6}{The full binary tree is too big}

Nothing could be nicer than the two trees we have just described. 
They are the prototypes of networks having
manifestly finite resistance to infinity.
Unfortunately,
we can't even come close to finding either of these trees
as a subgraph of the three-dimensional lattice.
For in these trees,
the number of nodes in a ``ball'' of radius $r$
grows exponentially with $r$,
whereas in a $d$-dimensional lattice,
it grows like $r$, i.e., much slower.
(See Figure \figref{6.10}.) 
\prefig{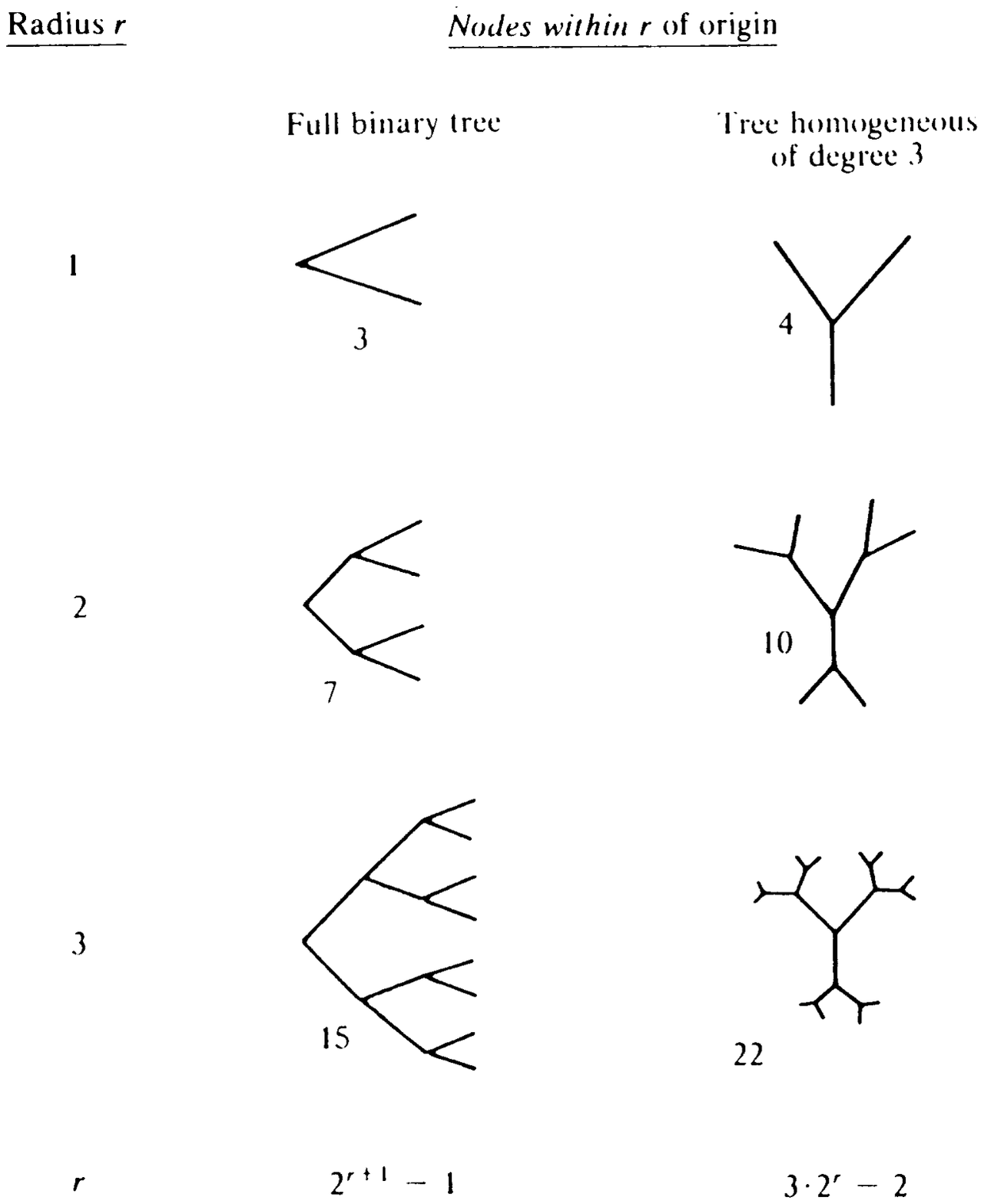}
There is simply no room for these trees
in any finite-dimensional lattice.

\sect{6.7}{$\NTthree$: a ``three-dimensional'' tree}

These observations suggest that we would do well to look for
a nice tree $\NTthree$
where the number of nodes within a radius $r$ of the root
is on the order of $r^3$.
For we might hope to find something resembling $\NTthree$
in the 3-dimensional lattice,
and if there is any justice in the world, 
this tree would have finite resistance to infinity,
and we would be done.

Before introducing $\NTthree$,
let's have a look at $\NTtwo$,
our choice for the tree most likely to succeed
in the 2-dimensional lattice
(see Figure \figref{6.11}).
\prefig{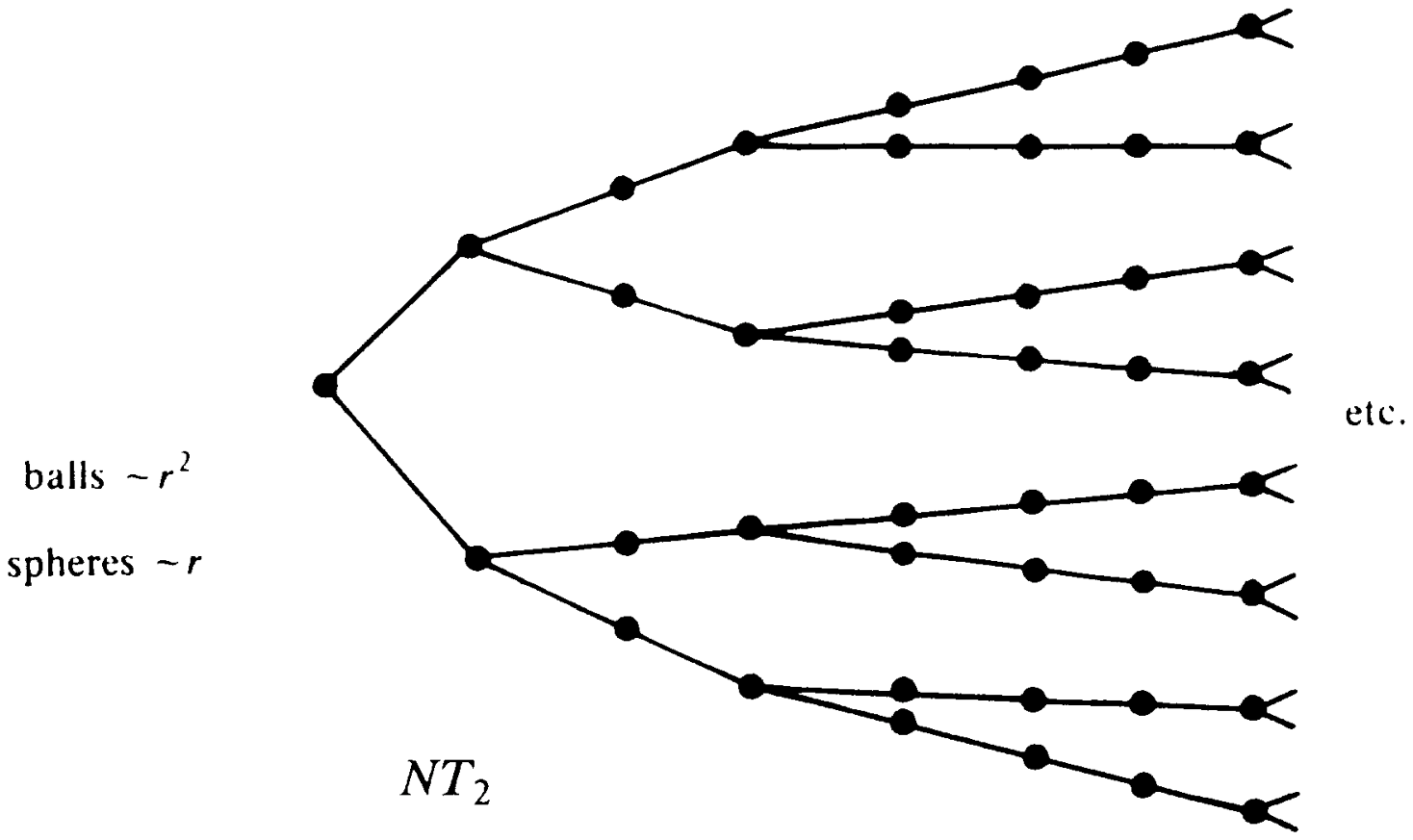}
The idea behind $\NTtwo$ is that,
since a ball of radius $r$ in the graph
ought to contain something like $r^2$ points,
a sphere of radius $r$ ought to contain something like $r$ points,
so the number of points in a sphere 
should roughly double
when the radius of the sphere is doubled.
For this reason,
we make the branches of our tree split in two
every time the distance from the origin is (roughly) doubled.

Similarly, in a 3-dimensional tree,
when we double the radius, 
the size of a sphere should roughly quadruple.
Thus in $\NTthree$,
we make the branches of our tree split in four
where the branches of $\NTtwo$ would have split in two.
$\NTthree$ is shown in Figure \figref{6.12}.
\prefig{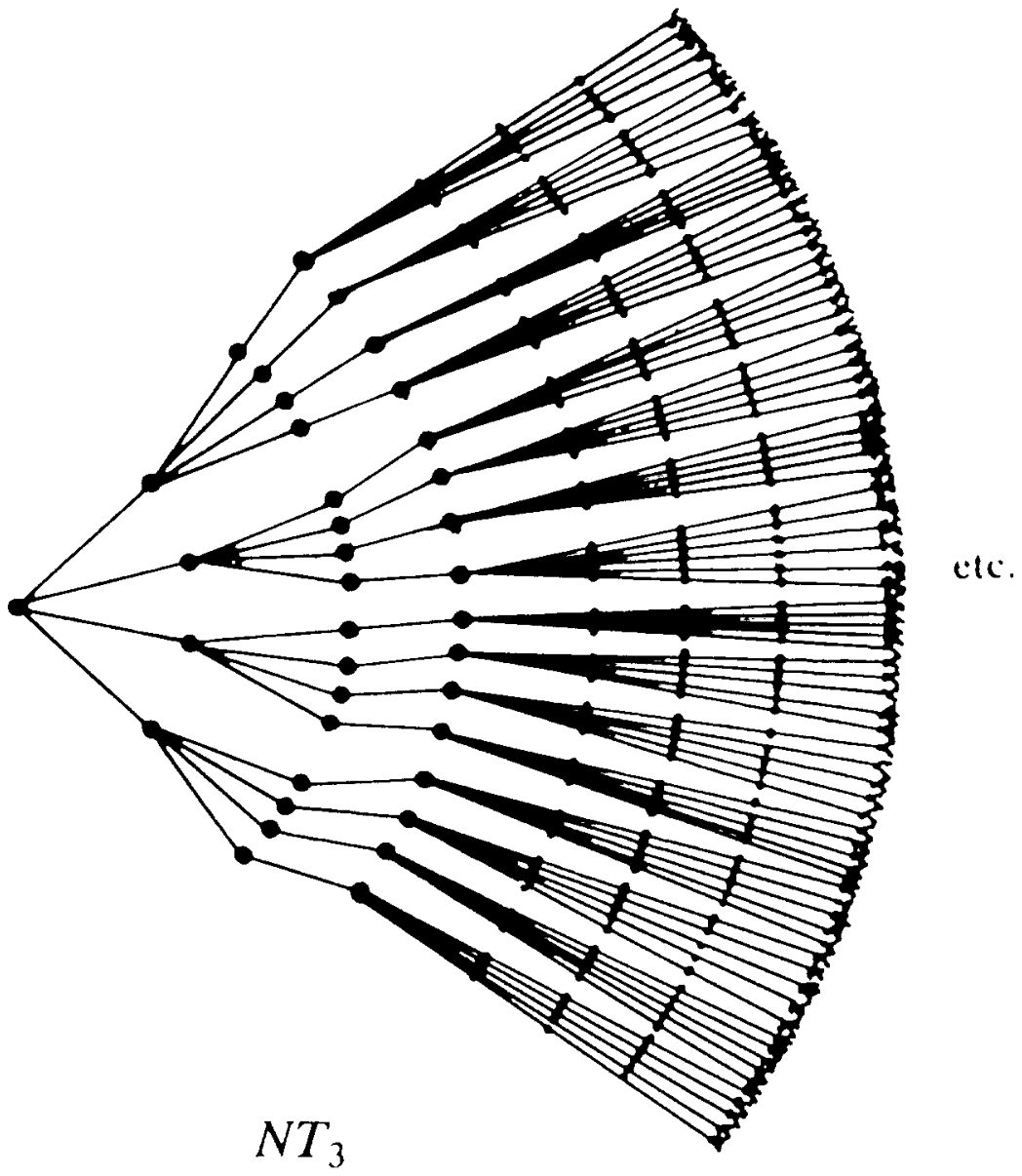}
Obviously, $\NTthree$ is none too happy about being drawn in the plane.

\sect{6.8}{$\NTthree$ has finite resistance}

To see if we're on the right track,
let's work out the resistance of our new trees.
These calculations are shown
in Figures \figref{6.13} and \figref{6.14}.
\prefig{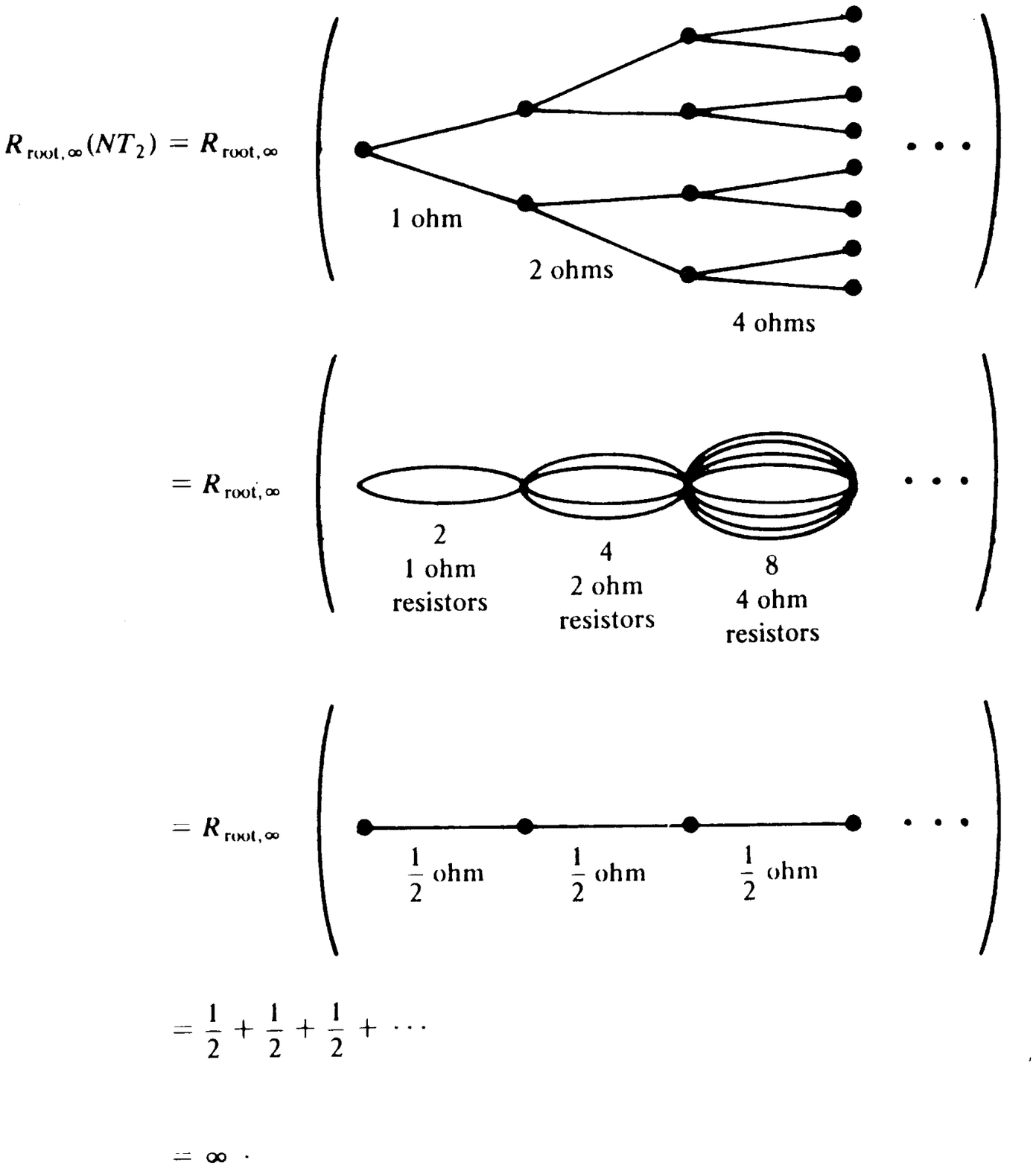}
\prefig{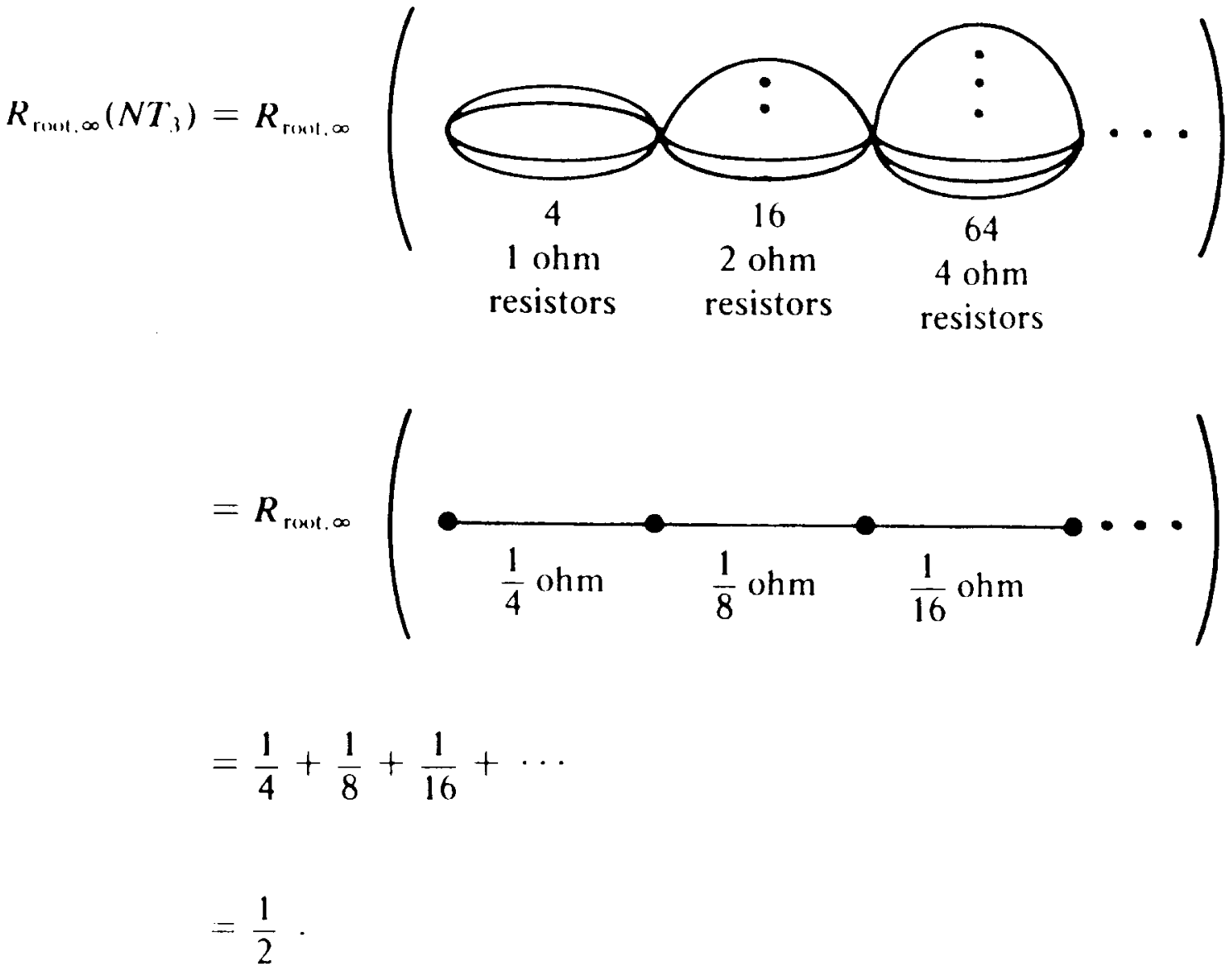}

As we would hope, the resistance of $\NTtwo$ is infinite,
but the resistance of $\NTthree$ is not.

\exstart
\ex{6.8.1}{
Use self-similarity arguments
to compute the resistance of $\NTtwo$ and $\NTthree$.
}
\exend

\sect{6.9}{But does $\NTthree$ fit in the three-dimensional lattice?}

We would like to embed $\NTthree$ in $\Z^3$.
We start by trying to embed $\NTtwo$ in $\Z^2$.
The result is shown in Figure \figref{6.15}.
\prefig{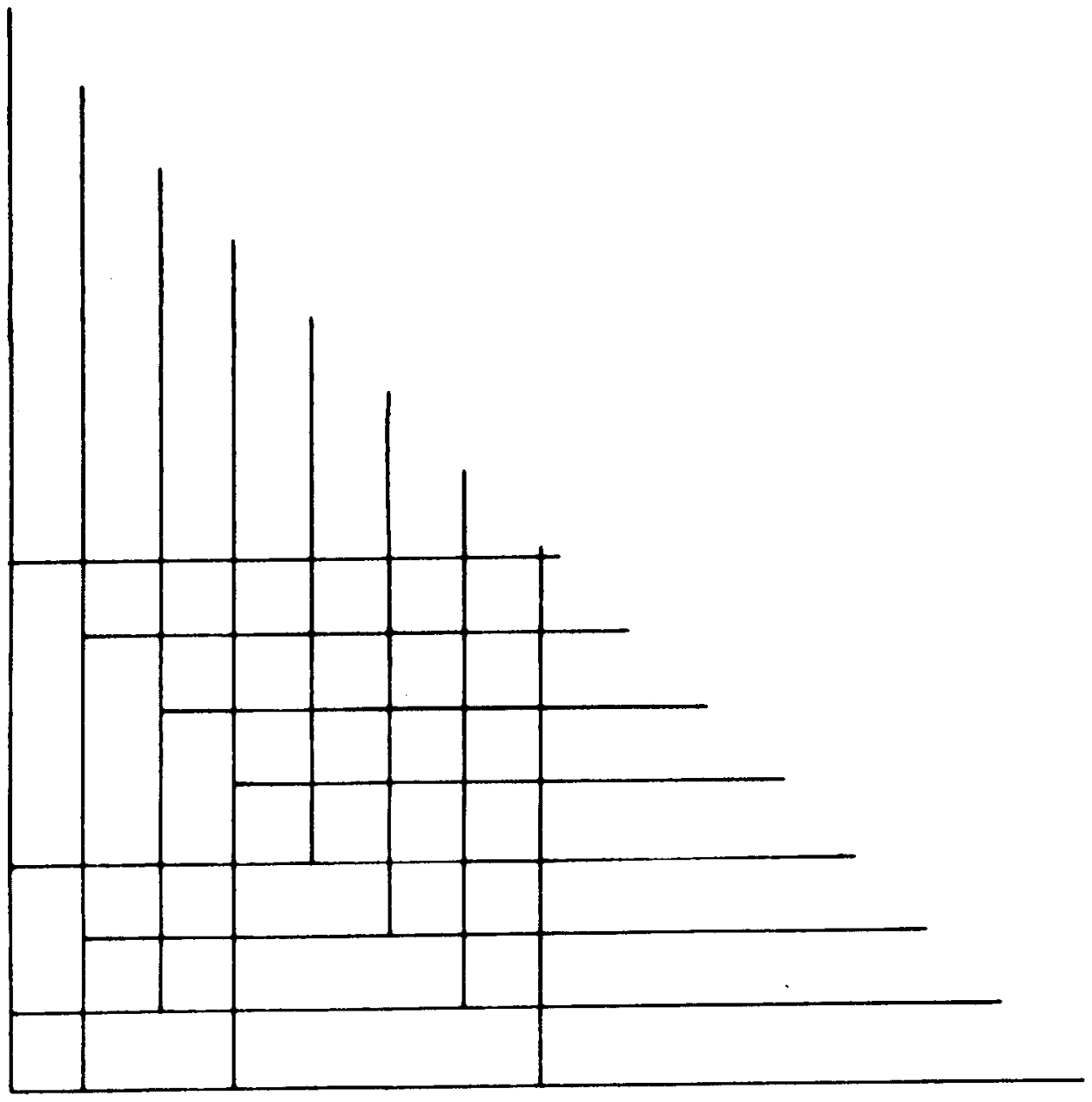}

To construct this picture,
we start from the origin and draw 2 rays,
one going north, one going east.
Whenever a ray intersects the line
$x + y = 2^n-1$ for some $n$,
it splits into 2 rays, one going north, and one going east.
The sequence of pictures in Figure \figref{6.16}
shows successively larger portions of the graph,
along with the corresponding portions of $\NTtwo$.
\prefig{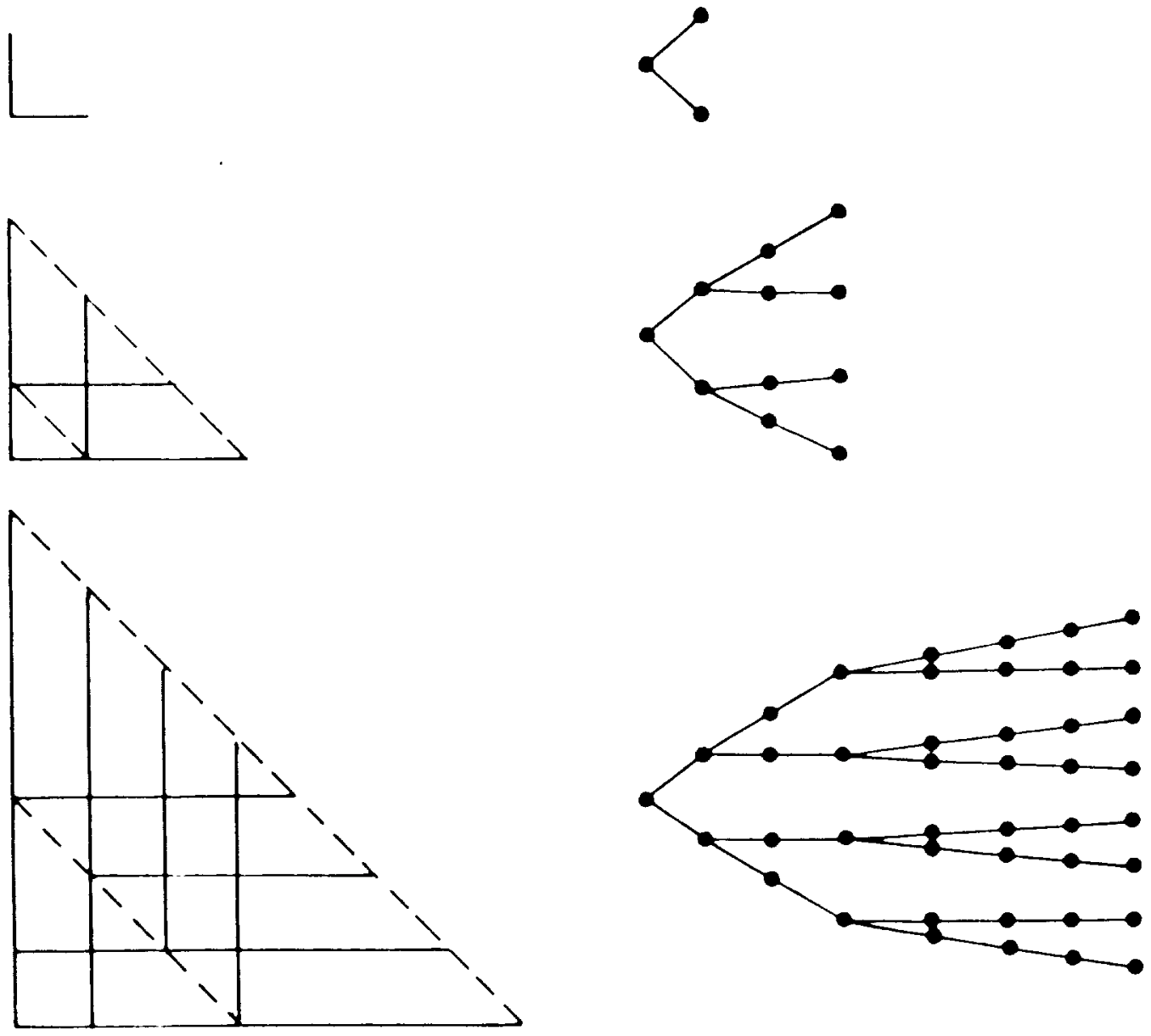}

Of course this isn't really an embedding,
since certain pairs of points that were distinct in $\NTtwo$
get identified, that is,
they are made to correspond to a single point in the picture.
In terms of our description, 
sometimes a ray going north and a ray going east
pass through each other.
This could have been avoided by allowing the rays to ``bounce''
instead of passing through each other,
at the expense of embedding
not $\NTtwo$ but a close relative---see Exercise \exref{6.9.3}.
However,
because the points of each identified pair
are at the same distance from the root of $\NTtwo$, 
when we put a battery between the root and the $n$th level
they will be at the same potential.
Hence, the current flow is not affected by these identifications,
so the identifications have no effect on $\Reff$.
For our purposes, then,
we have done just as well as if we had actually embedded $\NTtwo$.

To construct the analogous picture in three dimensions,
we start three rays off from the origin going north, east, and up.
Whenever a ray intersects the plane $x + y + z = 2^n - 1$ for some $n$,
it splits into three rays, going north, east, and up.
This process is illustrated in Figure \figref{6.17}.
\prefig{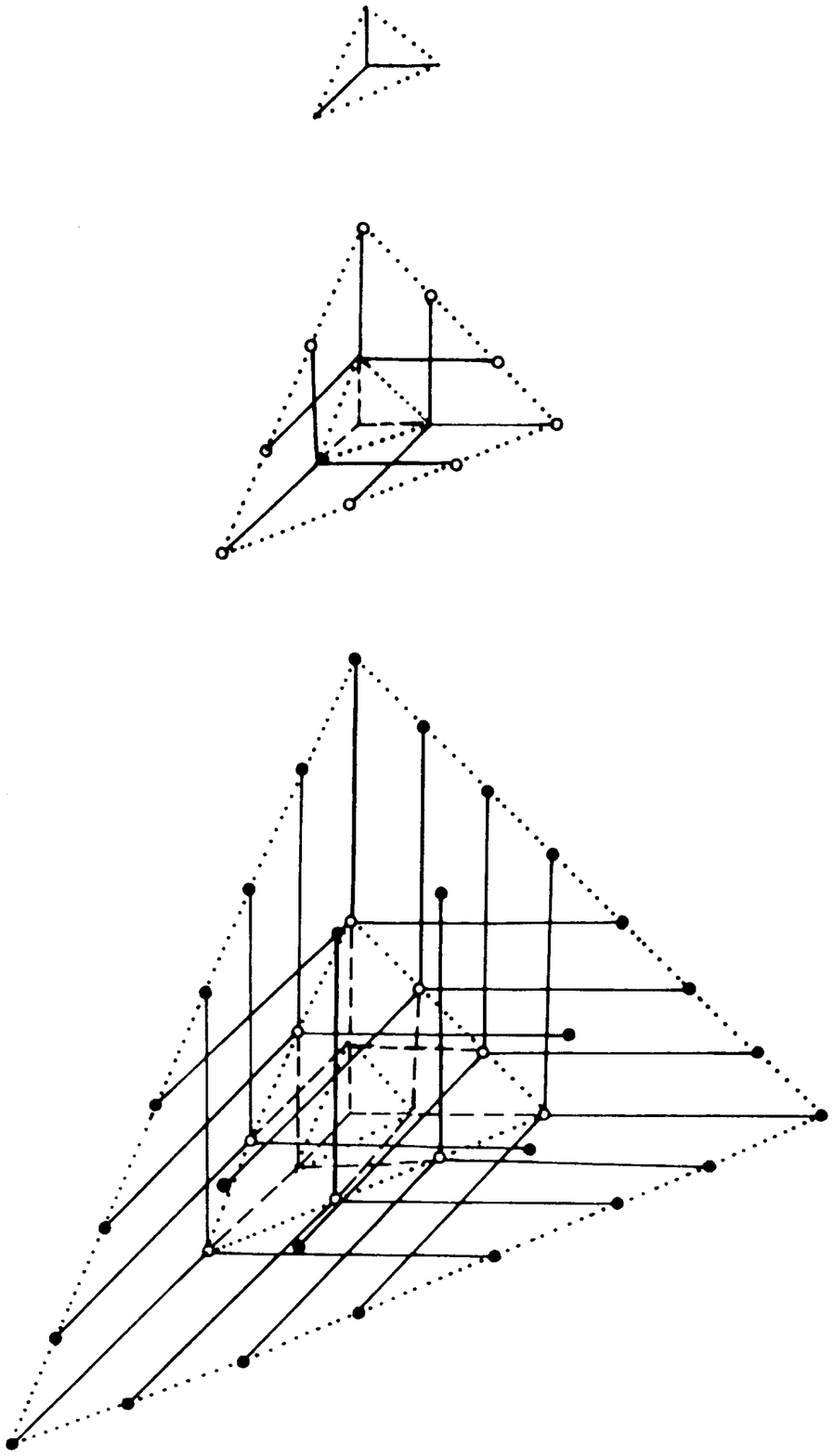}
Surprisingly, the subgraph of the 3-dimensional lattice obtained 
in this way is not $\NTthree$!
Rather, it represents an attempt to embed
the tree shown in Figure \figref{6.18}.
\prefig{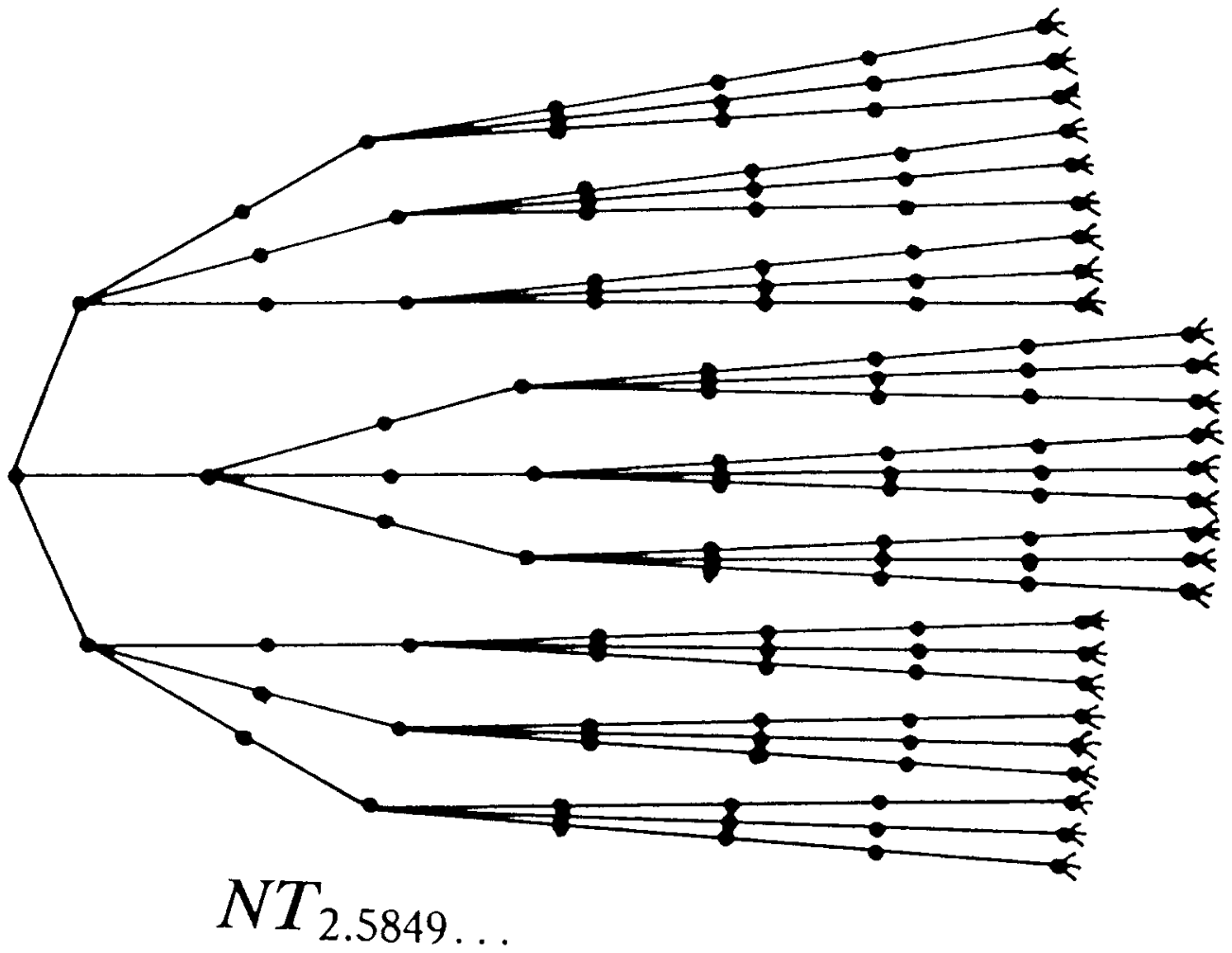}
We call this tree $\NTsemi$
because it is $2.5849\ldots$-dimensional
in the sense that when you double the radius of a ball,
the number of points in the ball gets multiplied roughly by 6
and
\[
6 = 2^{\log_2 6} = 2^{2.5849\ldots}
.
\]
Again, certain pairs of points of $\NTsemi$
have been allowed to correspond to the same point in the lattice,
but once again the intersections have no effect on $\Reff$.

So we haven't come up with our embedded $\NTthree$ yet.
But why bother?
The resistance of $\NTsemi$ out to infinity is
\[
\frac{1}{3} + \frac{2}{9} + \frac{4}{27} + \ldots
=
\frac{1}{3} \braces{ 1 + \frac{2}{3} + (\frac{2}{3})^2 + \ldots }
=
\frac{1}{3} \recip{1-\frac{2}{3}} = 1
.
\]

Thus we have found an infinite subgraph of the 3-dimensional lattice
having finite resistance out to infinity,
and we are done.

\exstart
\ex{6.9.1}{
This exercise deals with the escape probability $p_\esc$ for 
simple random walk in 3 dimensions.
The idea is to turn upper and lower bounds for the resistance of the lattice
into bounds for $p_\esc$.
Bounds are the best we can ask for using our method.
The determination of the exact value will be discussed in
\sectref{7.5}.
It is roughly .66.

(a) Use a shorting argument to get an upper bound for $p_\esc$.

(b) We have seen that the resistance of the 3-dimensional lattice 
is at most one ohm.
Show that the corresponding lower bound for $p_\esc$ is $1/6$
Show that this lower bound can be increased to $1/3$
with no extra effort.
}

\ex{6.9.2}{
Prove that simple random walk in any dimension $d > 3$
is transient.
}

\ex{6.9.3}{
Show how the not-quite embeddings of $\NTtwo$ and $\NTsemi$
can be altered to yield honest-to-goodness embeddings of
``stretched-out'' versions of these trees,
obtained by replacing each edge of the tree by three edges in series.
(Hint: ``bounce''.)
}
\exend

\sect{6.10}{What we have done; what we will do}

We have finally finished our electrical proof of Polya's theorem. 
The proof of recurrence for $d = 1, 2$ was straight-forward,
but this could hardly be said of the proof for $d = 3$.
After all,
we started out trying to embed $\NTthree$
and ended up by not quite embedding
something that was not quite $\NTthree$!

This is not bad in itself,
for one frequently sets out to do something
and in the process of trying to do it gets a better idea.
The real problem is that this explicit construction
is just too clever, too artificial.
We seem to be saying that a simple random walk in 3 dimensions is transient 
because it happens to contain
a beautifully symmetrical subgraph
that is in some sense $2.5849\ldots$-dimensional!
Fine, but what if we hadn't stumbled upon this subgraph?
Isn't there some other, more natural way?

We will see that indeed there are more natural approaches to showing
transience for $d = 3$.
One such approach uses the same idea of embedding trees,
but depends on the observation that
one doesn't need to be too careful about sending edges to edges.
Another approach,
based on relating the lattice not to a tree but to Euclidean space,
was already hinted at in \sectref{5.8}.
The main goal for the rest of this \book\ will be to explore
these more natural electrical approaches to Polya's theorem.

Before jumping into this, however,
we are going to go back and take a look at
a classical---i.e., probabilistic---approach to Polya's theorem.
This will give us something to compare our electrical proofs with.

\chap{7}{The classical proofs of Polya's Theorem}

\sect{7.1}{Recurrence is equivalent to an infinite expected number of
returns}

For the time being, all of our random walks will be simple.
Let $u$ be the probability that a random walker,
starting at the origin,
will return to the origin.
The probability that the walker will be there exactly $k$ times
(counting the initial time) is
$u^k(1 - u)$.
Thus, if $m$ is the expected number of times at the origin,
\[
m = \sum_{k=1}{\infty} k u^{k-1} (1-u)
=
\recip{1-u}
.
\]
If $m = \infty$ then $u=1$,
and hence the walk is recurrent. If $m < \infty$
then $u < 1$,
so the walk is transient.
Thus $m$ determines the type of the walk.

We shall use an alternate expression for $m$.
Let $u_n$ be the probability that the walk, starting at $\orig$,
is at $\orig$ on the nth step.
Since the walker starts at $\orig$, $u_0=1$.
Let $e_n$ be a random variable that takes on the value 1
if, at time $n$ the walker is at $\orig$ and 0 otherwise.
Then
\[
T = \sum_{n = 0}^\infty e_n
\]
is the total number of times at $\orig$ and
\[
m = \Expect(T) = \sum_{n=0}^\infty \Expect(e_n)
.
\]
But $\Expect(e_n)=1 \cdot u_n + 0 \cdot (1-u_n) = u_n$.
Thus
\[
m = \sum_{n=0}^\infty u_n
.
\]
Therefore,
the walk will be recurrent if the series $\sum_{n=0}^\infty$
diverges and transient if it converges.

\exstart
\ex{7.1.1}{
Let $N_{\mat{x} \mat{y}}$
be the expected number of visits to $\mat{y}$ for a random 
walker starting in $\mat{x}$.
Show that $N_{\mat{x} \mat{y}}$is finite
if and only if the walk is transient.
}
\exend

\sect{7.2}{Simple random walk in one dimension}

Consider a random walker in one dimension, started at $\orig$.
To return to $\orig$,
the walker must take the same number of steps to the right
as to the left;
hence, only even times are possible.
Let us compute $u_{2n}$.
Any path that returns in $2n$ steps
has probability $1/2^n$.
The number of possible paths
equals the number of ways that we can choose the $n$ times to go 
right from the $2n$ possible times.
Thus
\[
u_{2n} = \binom{2n}{n} \recip{2^{2n}}
.
\]

We shall show that
$\sum_n u_{2n} = \infty$
by using Stirling's approximation:
\[
n! \sim \sqrt{2 \pi n} e^{-n} n^n
.
\]
Thus
\[
u_{2n}
=
\frac{(2n!)}{n!n!} \recip{2^{2n}}
\sim
\frac{\sqrt{2 \pi \cdot 2n} e^{-2n} (2n)^{2n}}
{(\sqrt{2 \pi n} e^{-n} n^n)^2 2^{2n}}
=
\recip{\sqrt{\pi n}}
.
\]
Therefore,
\[
\sum_n u_{2n}
\sim
\sum_n \recip{\sqrt{\pi n}}
=
\infty
\]
and a simple random walk in one dimension is recurrent.

Recall that this case was trivial by the resistor argument.

\exstart
\ex{7.2.1}{
We showed in \sectref{1.5}\ that
a random walker starting at $x$ with $0 < x < N$
has probability $x/N$ of reaching $N$ before 0.
Use this to show that a simple random walk in one dimension is recurrent.
}

\ex{7.2.2}{
Consider a random walk in one dimension that moves from $n$ to $n+1$
with probability $p$
and to $n - 1$
with probability $q = 1 - p$.
Assume that $p > 1/2$.
Let $h_x$
be the probability, starting at $x$,
that the walker ever reaches 0.
Use Exercise \exref{1.6.1}\ to show that
$h_x = (q/p)^x$ for $x \geq 0$ and $h_x = 1$ for $x<0$.
Show that this walk is transient.
}

\ex{7.2.3}{
For a simple random walk in one dimension,
it follows from Exercise \exref{1.5.4}\ that the expected time,
for a walker starting at $x$ with $0 < x < N$,
to reach 0 or $n$ is $x(N - x)$.
Prove that for the infinite walk,
the expected time to return to 0 is infinite.
}

\ex{7.2.4}{
Let us regard a simple random walk in one dimension as
the fortune of a player in a penny matching game
where the players have unlimited credit.
Show that the result is a martingale (see \sectref{1.6}).
Show that you can describe a stopping system
that guarantees that you make money.
}
\exend

\sect{7.3}{Simple random walk in two dimensions}

For a random walker in two dimensions to return to the origin, 
the walker must have gone the same number of times
north and south and the same number of times east and west.
Hence, again,
only even times for return are possible.
Every path that returns in $2n$ steps
has probability $1/4^{2n}$.
The number of paths that do this by taking $k$ steps to the north, 
$k$ south, $n - k$ east and $n - k$ west is
\[
\multinom{2n}{k,k,n-k,n-k}
=
\frac{2n!}{k!k!(n-k)!(n-k)!}
.
\]
Thus
\begin{eqnarray*}
u_{2n}
&=&
\recip{4^{2n}}
\sum_{k=0}^n
\frac{(2n)!}{k!k!(n-k)!(n-k)!}
\\&=&
\recip{4^{2n}}
\sum_{k=0}^n
\frac{(2n)!}{n!n!} \frac{n!n!}{k!k!(n-k)!(n-k)!}
\\&=&
\recip{4^{2n}}
\binom{2n}{n} \sum_{k=0}^n \binom{n}{k}^2
.
\end{eqnarray*}
But
$\sum_{k=0}^n \binom{n}{k}^2 = \binom{2n}{n}$
(see Exercise \exref{7.3.1}).
Hence
\[
u_{2n}
=
\braces{\recip{2^{2n}} \binom{2n}{n}}^2
.
\]
This is just the square of the one dimension result
(not by accident, either---see \sectref{7.6}).
Thus we have
\[
m =
\sum_n u_{2n}
\approx \sum_n \recip{\pi n}
= \infty
.
\]

Recall that the resistor argument was very simple in this case also.

\exstart
\ex{7.3.1}{
Show that 
$\sum_{k=0}^n \binom{n}{k}^2 = \binom{2n}{n}$
(Hint: Think of choosing $n$ balls from a box that has $2n$ balls,
$n$ black and $n$ white.)
}
\exend

\sect{7.4}{Simple random walk in three dimensions}

For a walker in three dimensions to return,
the walker must take an equal number of steps back and forth
in each of the three different directions.
Thus we have
\[
u_{2n} =
\recip{6^{2n}}
\sum_{j,k} \frac{(2n)!}{j!j!k!k!(n-j-k)!(n-j-k)!}
\]
where the sum is taken over all $j,k$ with $j+k \leq n$.
Following Feller \cite{feller1},
we rewrite this sum as
\[
u_{2n} =
\recip{2^{2n}} \binom{2n}{n}
\sum_{j,k}
\braces{\recip{3^n} \frac{n!}{j!k!(n-j-k)!}}^2
.
\]
Now consider placing $n$ balls randomly into three boxes, $A, B, C$. 
The probability that
$j$ balls fall in $A$, $k$ in $B$, and $n -j - k$ in $C$
is
\[
\recip{3^n}\multinom{n}{j,k,n-j-k}
=
\recip{3^n} \frac{n!}{j!k!(n-j-k)!}
.
\]
Intuitively, this probability is largest when
$j$, $k$, and $n - j - k$
are as near as possible to $n/3$,
and this can be proved (see Exercise \exref{7.4.1}).
Hence, replacing one of the factors in the square by this larger value,
we have:
\[
u_{2n} \leq
\recip{2^{2n}}
\binom{2n}{n}
\braces{\recip{3^n}
  \frac{n!}{\floor{\frac{n}{3}}!\floor{\frac{n}{3}}!\floor{\frac{n}{3}}!}}
\braces{\sum_{j,k}\recip{3^n} \frac{n!}{j!k!(n-j-k)!}}
,
\]
where $\floor{n/3}$ denotes the greatest integer $\leq n/3$.
The last sum is 1
since it is the sum of all probabilities
for the outcomes of putting $n$ balls into three boxes.
Thus
\[
u_{2n} \leq
\recip{2^{2n}}
\binom{2n}{n}
\braces{\recip{3^n} \frac{n!}{\floor{\frac{n}{3}}!^3}}
.
\]
Applying Stirling's approximation yields
\[
u_{2n} \leq
\frac{K}{n^{3/2}}
\]
for suitable constant $K$.
Therefore
\[
m = \sum_n u_{2n}
\leq K \sum_n \recip{n^{3/2}}
< \infty
,
\]
and a simple random walk in three dimensions is recurrent.

While this is a complex calculation,
the resistor argument was also complicated in this case.
We will try to make amends for this presently.

\exstart
\ex{7.4.1}{
Prove that $\binom{n}{j,k,n-j-k}$ is largest
when $j$, $k$, and $n -j - k$ are as close as possible to $n/3$.
}

\ex{7.4.2}{
Find an appropriate value for the "suitable constant" $K$
that was mentioned above,
and derive an upper bound for $m$.
Use this to get a lower bound
for the probability of escape
for simple random walk in three dimensions.
}
\exend

\sect{7.5}{The probability of return in three dimensions: exact calculations}

Since the probability of return in three dimensions
is less than one,
it is natural to ask,
``What is this probability?''
For this we need an exact calculation.
The first such calculation was carried out
in a pioneering paper of McCrea and Whipple \cite{mccreawhipple}.
The solution outlined here follows
Feller \cite{feller1}, Exercise 28, Chapter 9,
and Montroll and West \cite{montrollwest}.

Let $p(a, b, c; n)$ be the probability that a random walker, 
starting at $\orig$,
is at $(a, b, c)$ after $n$ steps.
Then $p(a, b, c; n)$ is completely determined by the fact that
\[
p(0, 0, 0; 0) =1
\]
and
\begin{eqnarray*}
p (a, b, c; n)
&=&
\recip{6} p(a - 1, b, c; n-1)
+
\recip{6} p(a + 1, b, c; n-1)
+
\\&&
\recip{6} p(a, b-1, c; n-1)
+
\recip{6} p(a, b+1, c; n-1)
+
\\&&
\recip{6} p(a, b, c-1; n-1)
+
\recip{6} p(a, b, c+1; n-1)
.
\end{eqnarray*}

Using the technique of generating functions,
it is possible to derive a solution of these equations as
\begin{eqnarray*}
&&p (a, b, c; n)
\\&=&
\frac{1}{(2 \pi)^3}
\cdot
\\&&
\int_{-\pi}^\pi
\int_{-\pi}^\pi
\int_{-\pi}^\pi
\braces{\frac{\cos x + \cos y + \cos z}{3}}^n
\cos(xa) \cos(yb) \cos(zc) dx dy dz
.
\end{eqnarray*}
Of course,
we can just verify that this formula satisfies our equations once
someone has suggested it.
Having this formula, we put $a = b = c = 0$
and sum over $n$ to obtain the expected number of returns as
\[
m =
\frac{3}{(2 \pi)^3}
\int_{-\pi}^\pi
\int_{-\pi}^\pi
\int_{-\pi}^\pi
\recip{3-(\cos x + \cos y + \cos z)}
dx dy dz
.
\]
This integral was first evaluated by Watson \cite{watson}
in terms of elliptic integrals, which are tabulated.
A simpler result was obtained by Glasser and Zucker \cite{glasserzucker}
who evaluated this integral in terms of gamma functions.
Using this result, we get
\[
m = \frac{\sqrt{6}}{32\pi^3}
\Gamma \braces{\frac{1}{24}}
\Gamma \braces{\frac{5}{24}}
\Gamma \braces{\frac{7}{24}}
\Gamma \braces{\frac{11}{24}}
=
1.516386059137\ldots
,
\]
where
\[
\Gamma(x) = \int_0^\infty e^{-t} t^{x-1} dt
\]
is Euler's gamma function.
(Incidentally, the value given by Glasser and Zucker \cite{glasserzucker}
for the integral above needs to be corrected
by multiplying by a factor of $1/(384\pi)$.)

Recall that $m = 1/(1 - u)$ so that $u = (m - 1)/m$.
This gives
\[
u = .340537329544 \ldots
.
\]

\sect{7.6}{Simple random walk in two dimensions is the same as two
independent one-dimensional random walks}

We observed that the probability of return at time $2n$ in two dimensions
is the square of the corresponding probability in one dimension. 
Thus it is the same as the probability that two independent walkers,
one walking in the $x$ direction and the other in the $y$ direction,
will, at time $2n$, both be at $0$.
Can we see that this should be the case?
The answer is yes.
Just change our axes by 45 degrees to new axes $\bar{x}$ and $\bar{y}$
as in Figure \figref{7.1}.
\prefig{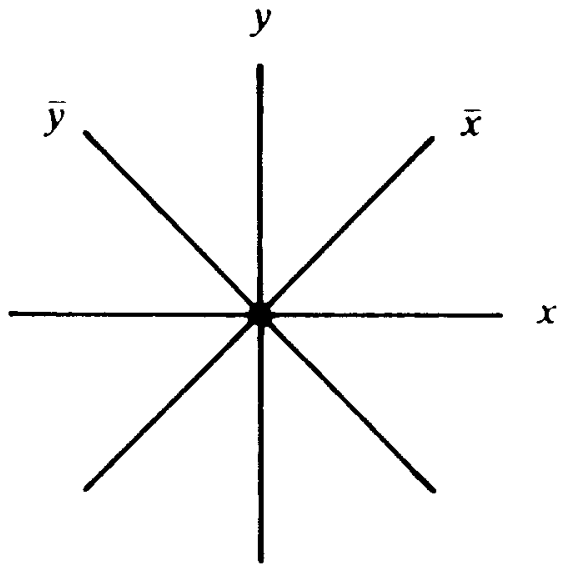}

Look at the possible outcomes for the first step using the $x, y$
coordinates and the $\bar{x},\bar{y}$ coordinates.
We have
\[
\begin{array}{ll}
\mbox{$x,y$ coordinates} & \mbox{$\bar{x},\bar{y}$ coordinates} \\
(0,1) & (1/\sqrt{2},1/\sqrt{2}) \\
(0,-1) & (-1/\sqrt{2},-1/\sqrt{2}) \\
(1,0) & (1/\sqrt{2},-1/\sqrt{2}) \\ 
(-1,0) & (-1/\sqrt{2},1/\sqrt{2})
\end{array}
.
\]
Assume that we have two independent walkers,
one moving with step size $\recip{\sqrt{2}}$ randomly
along the $\bar{x}$ axis and the other moving 
with the same step size along the $\bar{y}$ axis.
Then, if we plot their positions using the $x,y$ axes,
the four possible outcomes for the first step
would agree with those given in the second column of the
table above.
The probabilities for each of the four pairs of outcomes
would also be $(1/2) \cdot (1/2) = 1/4$.
Thus, we cannot distinguish a simple random walk in two dimensions
from two independent walkers along the $\bar{x}$ and $\bar{y}$ axes
making steps of magnitude $1/\sqrt{2}$.
Since the probability of return does not depend upon
the magnitude of the steps,
the probability that our two independent walkers
are at $(0, 0)$ at time $2n$
is equal to the product of the probabilities
that each separately is at 0 at time $2n$,
namely $(1/2^{2n}) \binom{2n}{n}$.
Therefore, the probability that the standard walk
will be at $(0, 0)$ at time $2n$ is
$((1/2^{2n}) \binom{2n}{n})^2$
as observed earlier.

\sect{7.7}{Simple random walk in three dimensions is not the same as
three independent random walks}

In three dimensions,
the probability that three independent walkers
are each back to 0 after time $2n$ is
\[
u_{2n} =
\braces{\binom{2n}{n} \recip{2^{2n}}}^3
.
\]
This does not agree with
our result for a simple random walk in three dimensions.
Hence, the same trick cannot work.
However, it is interesting to consider a random walk
which is the result of three independent walkers.
Let $(i, j, k)$ be the position of three independent random walkers.
The next position is one of the eight possibilities
$(i \pm 1, j \pm 1, k \pm 1)$
Thus we may regard their progress as a random walk
on the lattice points in three dimensions.
If we center a cube of side 2 at $(i, j, k)$,
then the walk moves with equal probabilities
to each of the eight corners of the cube.
It is easier to show that this random walk is transient
(using classical methods)
than it is for simple random walk.
This is because we can again use the one-dimension calculation.
The probability $u_{2n}$ for return at time $2n$ is
\[
u_{2n} =
\braces{\binom{2n}{n} \recip{2^{2n}}}^3
\sim \braces{\recip{\sqrt{\pi n}}}^3
.
\]
Thus
\[
m =
\sum_n u_{2n}
\sim \sum_n \braces{\recip{\pi n}}^{3/2}
< \infty
\]
and the walk is transient.

The fact that this three independent walkers model
and simple random walk are of the same type
suggests that when two random walks are ``really about the same'',
they should either both be transient or both be recurrent.
As we will soon see, this is indeed the case.
Thus we may infer the transience of simple random walk in 3 dimensions
from the transience of the three independent walkers model
without going through the involved calculation of \sectref{7.4}.

\chap{8}{Random walks on more general infinite networks}

\sect{8.1}{Random walks on infinite networks}

From now on we assume that $G$ is an infinite connected graph.
We assume that it is of \dt{bounded degree},
by which we mean that there is some integer $E$
such that the number of edges from any point is at most $E$.
We assign to each edge $xy$ of $G$ a conductance $C_{xy} > 0$
with $R_{xy} = \recip{C_{xy}}$.
The graph $G$ together with the conductances $\mat{C} = (C_{xy})$
is called a \dt{network}
and denoted by $(G, \mat{C})$.
Given a network $(G, \mat{C})$,
we define a random walk by
\[
P_{xy} = \frac{C_{xy}}{C_x}
\]
where
$C_x = \sum_y C_{xy}$.
When all the conductances are equal,
we obtain a random walk
that moves along each edge with the same probability: 
In agreement with our previous terminology,
we call this walk \dt{simple random walk} on $G$.

We have now a quite general class of infinite-state Markov chains.
As in the case of finite networks,
the chains are reversible Markov chains:
That is, there is a positive vector $\mat{w}$
such that $w_x P_{xy} = w_y P_{yx}$.
As in the finite case,
we can take $w_x = C_x$,
since $C_x P_{xy} = C_{xy} = C_{yx} = C_y P_{yx}$.

\sect{8.2}{The type problem}

Let $(G, \mat{C})$ be an infinite network with random walk $\matP$.
Let $\orig$ be a reference point.
Let $p_\esc$ be the probability that a walk starting at $\orig$
will never return to $\orig$.
If $p_\esc=0$ we say that $\matP$ is \dt{recurrent},
and if $\pesc > 0$ we say that it is \dt{transient}.
You are asked to show in Exercise \exref{8.2.1} that
the question of recurrence or transience of $\matP$
does not depend upon the choice of the reference point.
The \dt{type problem}
is the problem of determining if a random walk (network)
is recurrent or transient.

In \sectref{5.5}\ we showed how to rephrase the type problem
for a lattice in terms of finite graphs sitting inside it.
In \sectref{5.6}\ we showed that
the type problem is equivalent to an electrical network problem
by showing that simple random walk on a lattice is recurrent
if and only if the lattice has infinite resistance to infinity.
The same arguments apply with only minor modifications to
the more general infinite networks as well.
This means that we can use Rayleigh's short-cut method
to determine the type of these more general networks.

\exstart
\ex{8.2.1}{
Show that the question of recurrence or transience of $\matP$ does 
not depend upon the choice of the reference point.
}
\exend

\sect{8.3}{Comparing two networks}

Given two sets of conductances $\mat{C}$ and $\mat{\Cbar}$ on $G$,
we say that $(G, \mat{\Cbar}) < (G, \mat{C})$
if $\Cbar_{xy} < C_{xy}$ for all $xy$,
or equivalently,
if $\bar{R}_{xy} > R_{xy}$ for all $xy$.
Assume that $(G, \mat{\Cbar}) < (G, \mat{C})$.
Then by the Monotonicity Law,
$\bar{R}_\eff \geq \Reff$.
Thus if random walk on $(G,\mat{\Cbar})$ is transient,
i.e., if $\bar{R}_\eff < \infty$,
then random walk on $(G, \mat{C})$ is also transient.
If random walk on $(G, \mat{C})$ is recurrent,
i.e., if $\Reff = \infty$,
then random walk on $(G,\mat{\Cbar})$ is also recurrent.

\tout{Theorem.}
If $(G, \mat{C})$ and $(G, \mat{\Cbar})$ are networks,
and if there exist constants $u, v$ with $0 < u \leq  v < \infty$
such that
\[
u C_{xy} \leq \Cbar_{xy} \leq v C_{xy}
\]
for all $x$ and $y$,
then random walk on $(G, \mat{\Cbar})$
is of the same type as
random walk on $(G, \mat{C})$.
\toutend

\proofstart
Let $U_{xy} = u C_{xy}$ and $V_{xy} = v C_{xy}$.
Then
$(G,\mat{U}) \leq (G, \mat{\Cbar}) \leq (G,\mat{V})$.
But the random walks for $(G,\mat{U})$ and $(G,\mat{V})$
are the same as random walk on $(G, \mat{C})$.
Thus random walk for $(G, \mat{\Cbar})$
is of the same type as random walk on $(G,\mat{C})$.
\proofend

\tout{Corollary.}
Let $(G, \mat{C})$ be a network.
If for every edge $xy$ of $G$
we have $0 < u < C_{xy} < v < \infty$ for some constants $u$ and $v$,
then the random walk on $(G, \mat{C})$
has the same type as simple random walk on $G$.
\toutend

\exstart
\ex{8.3.1}{
Consider the two-dimensional lattice.
For each edge,
we toss a coin
to decide what kind of resistor to put across this edge.
If heads turns up, 
we put a two-ohm resistor across this edge;
if tails turns up,
we put a one-ohm resistor across the edge.
Show that the random walk on the resulting network is recurrent.
}

\ex{8.3.2}{
Consider the analogous problem to Exercise \exref{8.3.1} in 3 dimensions.
}
\exend

\sect{8.4}{The $k$-fuzz of a graph}

For any integer $k$, the $k$-\dt{fuzz} of a graph $G$
is the graph $\fuzz{G}{k}$ obtained from $G$
by adding an edge $xy$ if it is possible to go from $x$ to $y$
in at most $k$ steps.
For example, the 2-fuzz of the two-dimensional lattice is shown 
in Figure \figref{8.1};
please note that horizontal and vertical edges of length 2,
such as those joining $(0, 0)$ to $(0, 2)$, have not been indicated.
\prefig{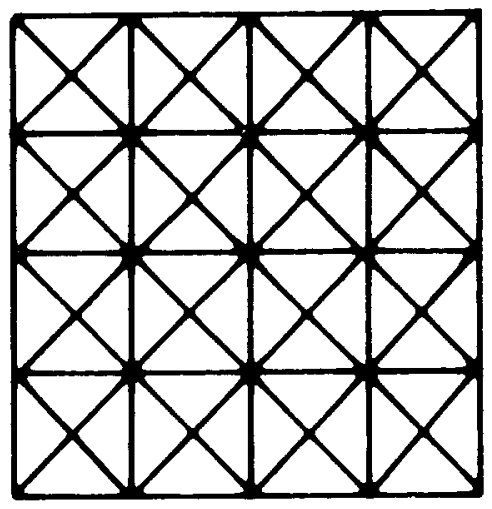}

\tout{Theorem.}
Simple random walk on $G$ and on the $k$-fuzz $\fuzz{G}{k}$ of $G$
have the same type.
\toutend

\proofstart
Let $\matP$ be simple random walk on $G$.
Define $\matPbar = (\matP + \matP^2 + \ldots + \matP^k)/k$.
Then $\matPbar$ may be considered to be $\matP$,
watched at one of the first $k$ steps chosen at random,
then at a time chosen at random from the next $k$ steps
after this time, etc.
Thinking of $\matPbar$ in this way,
we see that $\matP$ is in state $\orig$
at least once for every time $\matPbar$ in state $\orig$.
Hence, if $\matPbar$ is recurrent so is $\matP$.
Assume now that $\matPbar$ is transient.
Choose a finite set $S$
so that $\orig$ cannot be reached in $k$ steps
from a point outside of $S$.
Then, since the walk $\matPbar$ will be outside $S$ from some time on,
the walk $\matP$ cannot be at $\orig$ after this time,
and $\matP$ is also transient.
Therefore, $\matP$ and $\matPbar$ are of the same type.

Finally, we show that $\matPbar$ has the same type as
simple random walk on $\fuzz{G}{k}$.
Here it is important to remember our restriction that
$G$ is of bounded degree,
so that for some $E$ no vertex has degree $>E$.
We know that $\matP$ is reversible with
$\mat{w} \matP = \mat{w}$,
where $w_x$ is the number of edges coming out of $x$.
From its construction,
$\matPbar$ is also reversible and 
$\mat{w} \matPbar = \mat{w}$.
$\matPbar$ is the random walk on a network
$(G_k,\mat{\Cbar})$
with $\Cbar_{xy} = w_x \Pbar_{xy}$.
If $\Pbar_{xy} > 0$,
there is a path
$x,x_1,x_2,\ldots,x_{m-1},y$ in $G$ from $x$ to $y$
of length $m \leq k$.
Then
\[
\Pbar_{xy}
\geq \recip{k} (\recip{E})^m
\geq \recip{k} (\recip{E})^k
.
\]
Thus
\[
0 < \recip{k} (\recip{E})^k \leq \Pbar_{xy} \leq 1
\]
and
\[
0 < \recip{k} (\recip{E})^k \leq \Cbar_{xy} \leq E
.
\]
Therefore, by the theorem on the irrelevance of bounded twiddling
proven in \sectref{8.3},
$\matPbar$ and simple random walk on $\fuzz{G}{k}$ are of the same type.
So $G$ and $\fuzz{G}{k}$ are of the same type.

NOTE: This is the only place where we use probabilistic methods of
proof. For the purist who wishes to avoid probabilistic methods, 
Exercise \exref{8.6.2}\ indicates an alternative electrical proof.

We show how this theorem can be used.
We say that a graph $G$ can be
\dt{embedded} in a graph $\Gbar$
if the points $x$ of $G$
can be made to correspond in a one-to-one fashion
to points $\bar{x}$ of $\Gbar$
in such a way that if $xy$
is an edge in $G$, then $\bar{x}\bar{y}$ is an edge in $\Gbar$.

\tout{Theorem.}
If simple random walk on $G$ is transient,
and if $\Gbar$ can be embedded in a $k$-fuzz $\fuzz{\Gbar}{k}$ of $\Gbar$
then simple random walk on $\Gbar$ is also transient.
Simple random walk on $G$ and $\Gbar$
are of the same type if each graph can be embedded
in a $k$-fuzz of the other graph.
\toutend

\proofstart
Assume that simple random walk on $G$ is transient
and that $G$ can be embedded in a $k$-fuzz $\fuzz{\Gbar}{k}$ of $\Gbar$.
Since $\Reff$ for $G$ is finite 
and $G$ can be embedded in $\fuzz{\Gbar}{k}$,
$\Reff$ for $\fuzz{\Gbar}{k}$ is finite.
By our previous theorem,
the same is true for $\Gbar$ and simple random walk on $\Gbar$ is transient.

If we can embed $G$ in $\fuzz{\Gbar}{k}$ and $\Gbar$ in $\fuzz{G}{k}$,
then the random walk on $G$ is transient
if and only if the random walk on $\Gbar$ is.
\proofend

\exstart
\ex{8.4.1}{
We have assumed that there is a bound $E$
for the number of edges coming out of any point.
Show that if we do not assume this,
it is not necessarily true that $G$ and $\fuzz{G}{k}$ are of the same type.
(Hint: Consider a network
something like that shown in Figure \figref{8.2}.)
\prefig{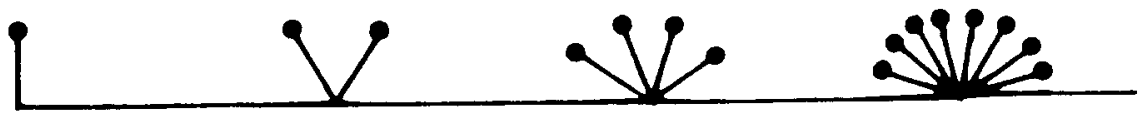}
}
\exend

\sect{8.5}{Comparing general graphs with lattice graphs}

We know the type of simple random walk on a lattice $\Z^d$.
Thus to determine the type of simple random walk on an arbitrary graph $G$, 
it is natural to try to compare G with $\Z^d$.
This is feasible for graphs that can be
drawn in some Euclidean space $\R^d$ in a civilized manner.

\tout{Definition.}
A graph $G$ can be drawn in a Euclidean space $\R^d$ in a
\dt{civilized manner}
if its vertices can be embedded in $\R^d$ so that 
for some $r < \infty$, $s > 0$
\propstart
\prop{(a)}
The length of each edge is $\leq r$.
\prop{(b)}
The distance between any two points is $> s$.
\propend
\toutend

Note that we make no requirement about
being able to draw the edges of $G$ so they don't intersect.

\tout{Theorem.}
If a graph can be drawn in $\R^d$ in a civilized manner, 
then it can be embedded in a $k$-fuzz of the lattice $\Z^d$.
\toutend

\proofstart
We carry out the proof for the case $d=2$.
Assume that $G$ can be drawn in a civilized manner in $\R^2$.
We want to show that $G$ can be embedded in a $k$-fuzz of $\Z^2$.
We have been thinking of $\Z^2$ as being drawn in $\R^2$
with perpendicular lines and adjacent points a unit distance 
apart on these lines,
but this embedding is only one particular way of representing $\Z^2$.
To emphasize this, let's talk about $L^2$ instead of $\Z^2$. 
Figure \figref{8.3} shows another way of drawing $L^2$ in $\R^2$.
\prefig{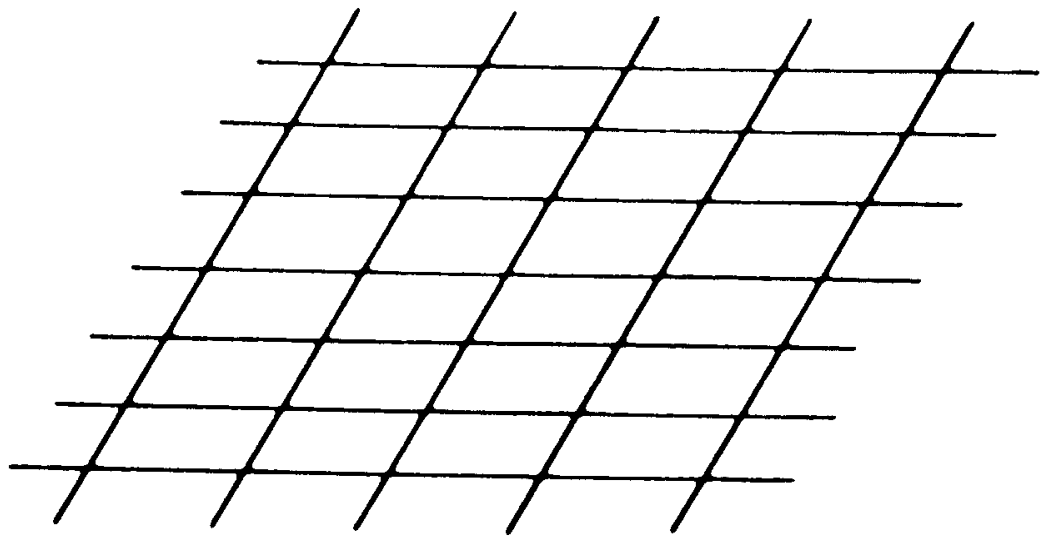}
From a graph-theoretical point of view,
this is the same as $\Z^2$.
In trying to compare $G$ to $L$,
we take advantage of this 
flexibility by drawing $L^2$
so small that points of $G$
can be moved onto points of $L^2$
without bumping into each other.

Specifically, let $L^2$ be a two-dimensional rectangular lattice 
with lines a distance $s/2$ apart.
In any square of $L^2$,
there is at most one point of $G$.
Move each point $x$ of $G$
to the southwest corner $\bar{x}$ of the square that it is in,
as illustrated in Figure \figref{8.4}.
\prefig{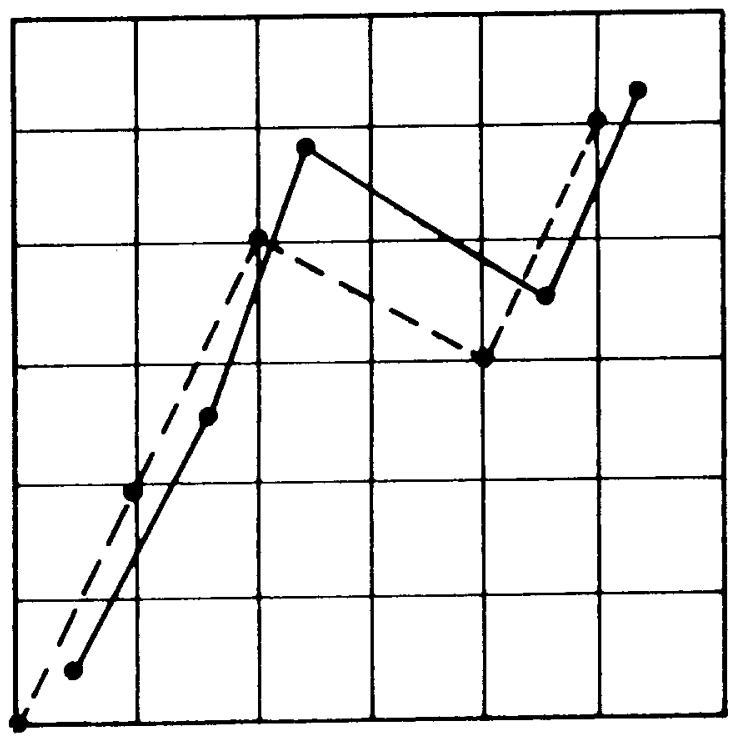}

Now since any two adjacent points $x$, $y$ in $G$ were within $r$ of each
other in $\R^2$,
the corresponding points $\bar{x}$, $\bar{y}$ in $L^2$
will have Euclidean distance $<r + 2s$.
Choose $k$ so that any two points of $L^2$
whose Euclidean distance is $< r + 2s$
can be connected by a path in $L^2$ of at most $k$ steps.
Then $\bar{x}$ and $\bar{y}$
will be adjacent in $\fuzz{L^2}{k}$
and---since the prescription for $k$ does not depend on $x$ and $y$---we
have embedded $G$ in the $k$-fuzz $\fuzz{L^2}{k}$.

\tout{Corollary.}
If $G$ can be drawn in a civilized manner in $\R^1$ or $\R^2$,
then simple random walk on $G$ is recurrent.
\toutend

\proofstart
Assume, for example,
that $G$ can be drawn in a civilized manner in $\R^2$.
Then $G$ can be embedded in a $k$-fuzz $\fuzz{\Z^2}{k}$ of $\Z^2$.
If simple random walk on $G$ were transient,
then the same would be true for $\fuzz{\Z^2}{k}$ and $\Z^2$.
But we know that simple random walk on $\Z^2$ is recurrent.
Thus simple random walk on $G$ is recurrent.
\proofend

Our first proof that random walk in three dimensions is transient
consisted in showing that we could embed a transient tree in $\Z^3$.
We now know that it would have been sufficient to show how to draw a 
transient tree in $\R^3$ in a civilized manner:
This is easier (see Exercise \exref{8.5.1}).

The corollary implies that simple random walk on any sufficiently
symmetrical graph in $\R^2$ is recurrent.
For example, simple random walk on the regular graph made up of hexagons
shown in Figure \figref{8.5}
is recurrent.
\prefig{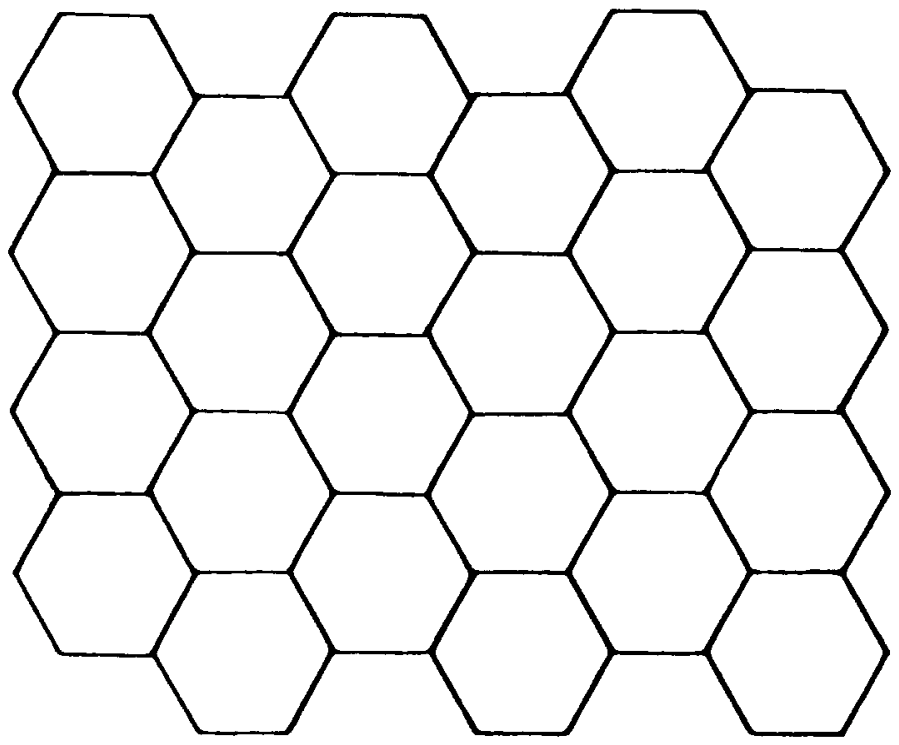}

We can even consider very irregular graphs.
For example, on the cover of the January 1977 Scientific American,
there is an example due to Conway of an infinite non-periodic tiling
using Penrose tiles of the form shown in Figure \figref{8.6}.
\prefig{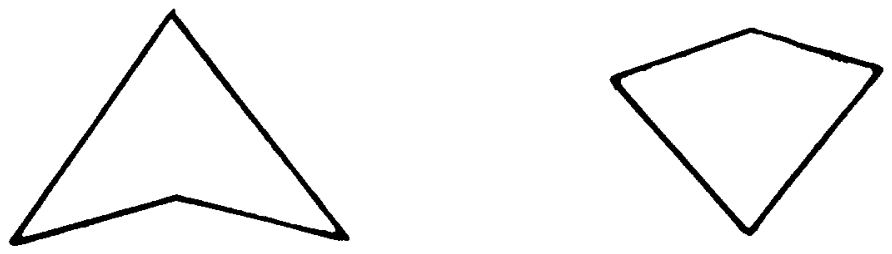}
It is called the cartwheel pattern;
part of it is shown in Figure \figref{8.7}.
\prefig{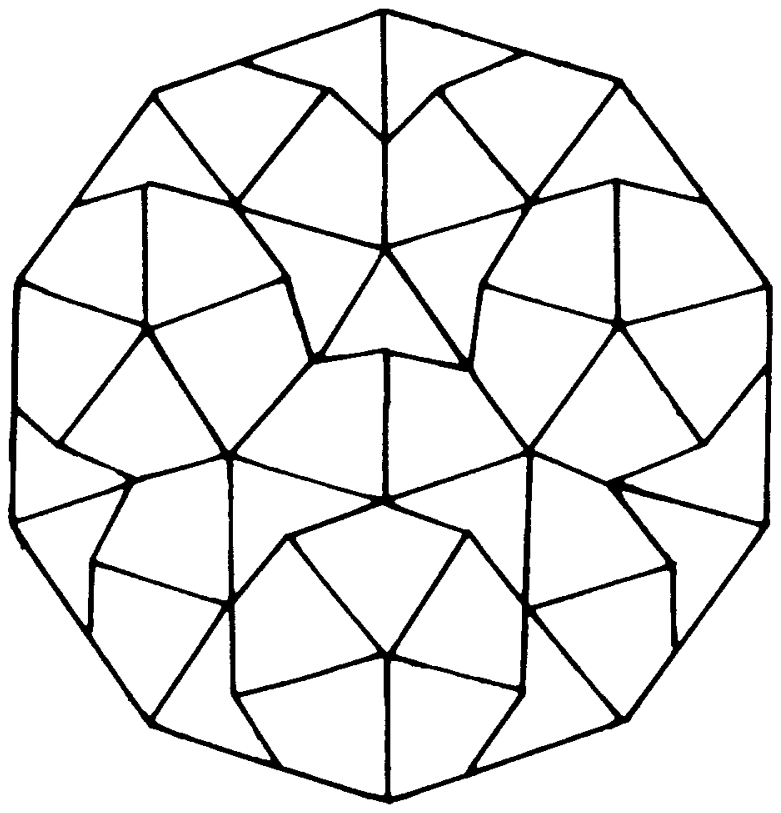}
A walker walking randomly on the edges
of this very irregular infinite tiling will still 
return to his or her starting point.

Assume now that $G$ can be drawn in a civilized manner in $\R^3$.
Then to show that simple random walk on $G$ is of the same type as $\Z^3$, 
namely transient, it is sufficient to show that we can embed $\Z^3$
in a $k$-fuzz of $G$.
This is clearly possible for any regular lattice in $\R^3$.
The three lattices that have been most studied
and for which exact probabilities for return have been found
are called the SC, BCC, and FCC lattices.
The SC (simple cubic) lattice is just $\Z^3$.
The walker moves each time to a new point
by adding a random choice from the six vectors
\[
(\pm 1, 0, 0), (0, \pm 1,0), (0,0,\pm 1)
.
\]
For the BCC (body-centered cubic) lattice,
the choice is one of the eight vectors
\[
(\pm 1, \pm 1, \pm 1)
.
\]
This was the walk that resulted from three independent
one-dimensional walkers.
For the FCC (face-centered cubic) lattice,
the random choice is made from the twelve vectors 
\[
(\pm 1, \pm 1, 0), (\pm 1,0,\pm 1), (0,\pm 1,\pm 1)
.
\]
For a discussion of exact calculations for these three lattices,
see Montroll and West \cite{montrollwest}

As we have seen,
once the transience of any one of these three walks is established,
no calculations are necessary to determine that the 
other walks are transient also.
Thus we have yet another way of establishing
Polya's theorem in three dimensions:
Simply verify transience of the walk on the
BCC lattice via the simple three-independent-walkers computation, 
and infer that walk on the SC lattice is also transient since the BCC 
lattice can be embedded in a $k$-fuzz of it.

\exstart
\ex{8.5.1}{
When we first set out to prove Polya's theorem for  $d=3$,
our idea was to embed $\NTthree$ in $\Z^3$.
As it turned out, what we ended up embedding
was not $\NTthree$ but $\NTsemi$,
and we didn't quite embed it at that.
We tried to improve the situation by finding
(in Exercise \exref{6.9.3}) an honest-to-goodness embedding
of a relative of $\NTsemi$, but
$\NTthree$ was still left completely out in the cold.
Now, however, we are in a position to embed $\NTthree$,
if not in $\Z^3$ then at least in a $k$-fuzz of it.
All we need to do is to draw $\NTthree$ in $\R^3$ in a civilized manner.
Describe how to do this, 
and thereby give one more proof of Polya's theorem for $d = 3$.
}

\ex{8.5.2}{
Find a graph that can be embedded in a civilized manner in $\R^3$
but not in $\R^2$, but is nonetheless recurrent.
}

\ex{8.5.3}{
Assume that $G$ is drawn in a civilized manner in $\R^3$.
To show that simple random walk on $G$ is transient,
it is enough to know that $\Z^3$ can be
embedded in a $k$-fuzz of G.
Try to come up with a nice condition that will guarantee
that this is possible.
Can you make this condition simple,
yet general enough so that it will settle
all reasonably interesting cases?
In other words, can you make the condition nice enough to allow us to
remember only the condition, and forget about the general method 
lying behind it?
}
\exend

\sect{8.6}{Solving the type problem by flows:  a variant of the cutting
method}

In this section we will introduce a variant of the cutting method 
whereby we use Thomson's Principle directly to estimate the effective 
resistance of a conductor.

Thomson's Principle says that,
given any unit flow through a resistive medium,
the dissipation rate of that flow gives an upper bound 
for the effective resistance of the medium.
This suggests that
to show that a given infinite network is transient,
it should be enough
to produce a unit flow out to infinity having finite energy dissipation.

In analogy with the finite case,
we say that $\mat{j}$ is a
\dt{flow from $\orig$ to infinity} if

\propstart
\prop{(a)}
$j_{xy} = j_{yx}$.
\prop{(b)}
$\sum_y j_{xy} = 0$ if $x \neq \orig$.
\propend

We define $j_\orig = \sum_y j_{\orig y}$.
If $j_\orig = 1$,
we say that $\mat{j}$ is a
\dt{unit flow to infinity}.
Again in analogy with the finite case,
we call $\half \sum_{x,y} j_{xy}^2 R_{xy}$
the \dt{energy dissipation}
of the flow $\mat{j}$.

\tout{Theorem.}
The effective resistance $\Reff$ from $\orig$ to $\infty$
is less than or equal to
the energy dissipation of any unit flow from $\orig$ to infinity.
\toutend

\proofstart
Assume that we have a unit flow $\mat{j}$
from $\orig$ to infinity
with energy dissipation
\[
E = \half \sum_{x,y} j_{xy}^2 R_{xy}
.
\]
We claim that $\Reff \leq E$.
Restricting $j_{xy}$ to the edges of the finite graph $\ball{G}{r}$,
we have a unit flow from $\orig$ to $\sphere{G}{r}$ in $\ball{G}{r}$.
Let $i^{(r)}$ be the unit current flow in $\ball{G}{r}$
from 0 to $\sphere{G}{r}$.
By the results of \sectref{3.5},
\[
\R^{(r)}_\eff
= 
\half \sum_{\ball{G}{r}} (i^{(r)}_{xy})^2 R_{xy}
\leq
\half \sum_{\ball{G}{r}} j_{xy}^2 R_{xy}
\leq
\half \sum_{x,y} j_{xy}^2 R_{xy}
=
E
,
\]
where $\sum_{\ball{G}{r}}$ indicates the sum over all pairs $x,y$
such that $xy$ is an edge of $\ball{G}{r}$.

\exstart
\ex{8.6.1}{
We have billed the method of
using Thomson's Principle directly
to estimate the effective resistances of a network
as a variant of the cutting method.
Since the cutting method was derived from Thomson's Principle,
and not vice versa,
it would seem that we have got the cart before the horse.
Set this straight by giving an informal (``heuristic'')
derivation of Thomson's Principle from the cutting method.
(Hint: see Maxwell \cite{maxwell}, Chapter VIII, Paragraph 307.) 
For more on this question, see Onsager \cite{onsager}.
}

\ex{8.6.2}{
Let $G$ be an infinite graph of bounded degree
and $\fuzz{G}{k}$ the $k$-fuzz of $G$.
Using electric network arguments,
show that $\Reff < \infty$ for $G$
if and only if $\Reff < \infty$ for $\fuzz{G}{k}$.
}
\exend

\sect{8.7}{A proof, using flows, that simple random walk in three
dimensions is transient}

We now apply this form of the cutting method to give another proof that
simple random walk on the threedimensional lattice is transient. 
All we need is a flow to infinity with finite dissipation.
The flow we are going to describe is not the first flow one would think of.
In case you are curious,
the flow described here was constructed as a side effect of an 
unsuccessful attempt to derive the isoperimetric inequality
(see Polya \cite{polya:solve}) 
from the ``max-flow min-cut'' theorem
(Ford and Fulkerson \cite{fordFulkerson}).
The idea is to find a flow in the positive orthant
having the property that the same amount
flows through all points at the same distance from $\orig$.

Again,
it is easiest to show the construction for the two-dimensional case.
Let $G$ denote the part of $\Z^2$ lying in the first quadrant. 
The graph $\ball{G}{4}$ is shown in Figure \figref{8.8}.
\prefig{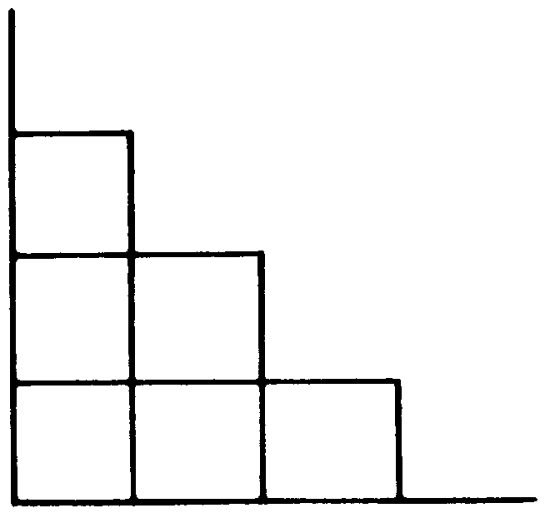}

We choose our flow so that it always goes away from $\orig$.
Into each point that is not on either axis there are two flows,
one vertical and one horizontal.
We want the sum of the corresponding values of $j_{xy}$
to be the same for all points the same distance from $\orig$.
These conditions completely determine the flow.
The flow out of the point $(x,y)$
with $x + y = n$
is as shown in Figure \figref{8.9}.
\prefig{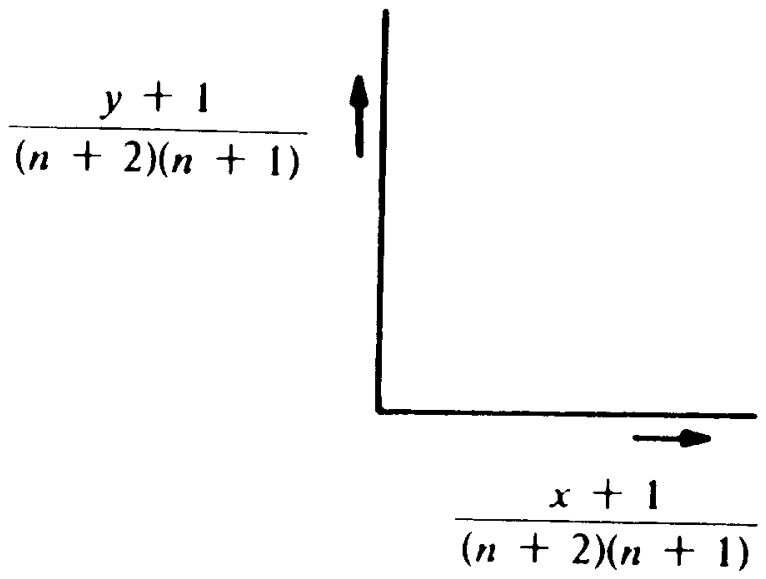}
The values for the currents out to the fourth
level are shown in Figure \figref{8.10}.
\prefig{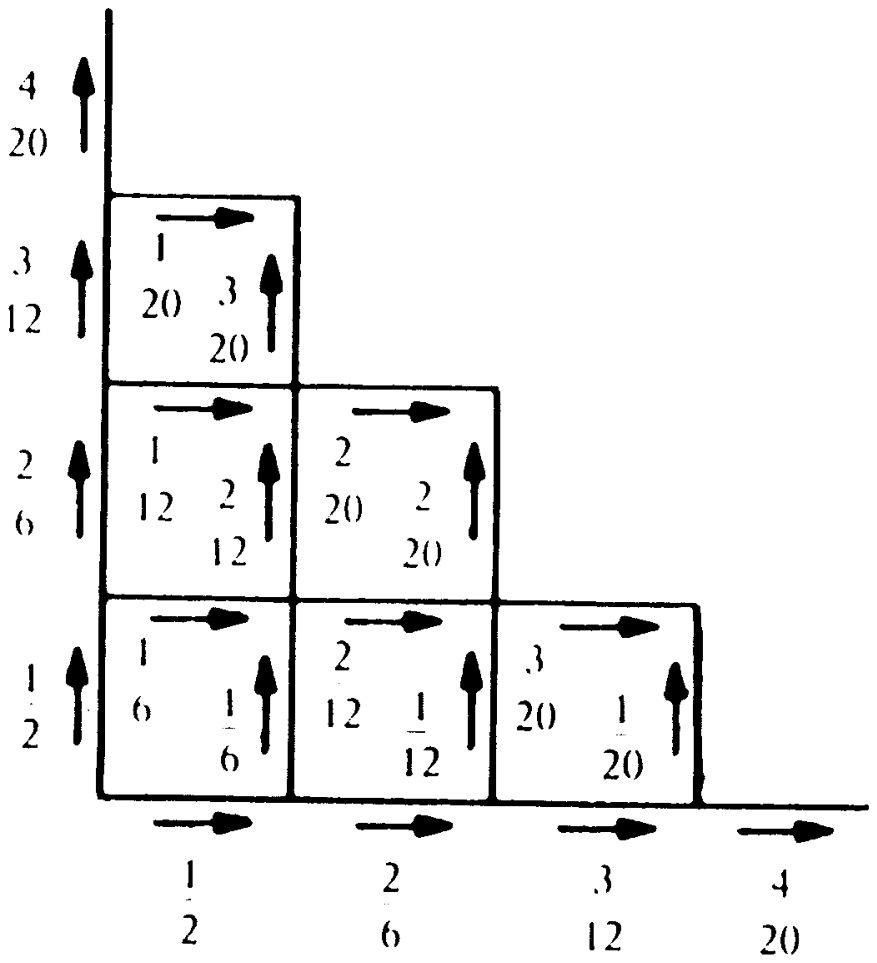}
In general, the flow out of a point $(x, y)$
with $x + y = n$ is
\[
\frac{x+1}{(n+2)(n+1)}
+ \frac{y+1}{(n+2)(n+1)}
=
\recip{n+1}
\]
and the flow into this point is
\[
\frac{x}{n(n+1)}
+ \frac{y}{n(n+1)}
=
\recip{n+1}
\]
Thus the net flow at $(x, y)$ is 0.
The flow out of $\orig$ is
$(1/2) + (1/2) = 1$.
For this two-dimensional flow,
the energy dissipation is infinite,
as it would have to be.
For three dimensions,
the uniform flow is defined as follows:
Out of $(x, y, z)$ with $x + y + z = n$
we have the flow indicated in Figure \figref{8.11}.
\prefig{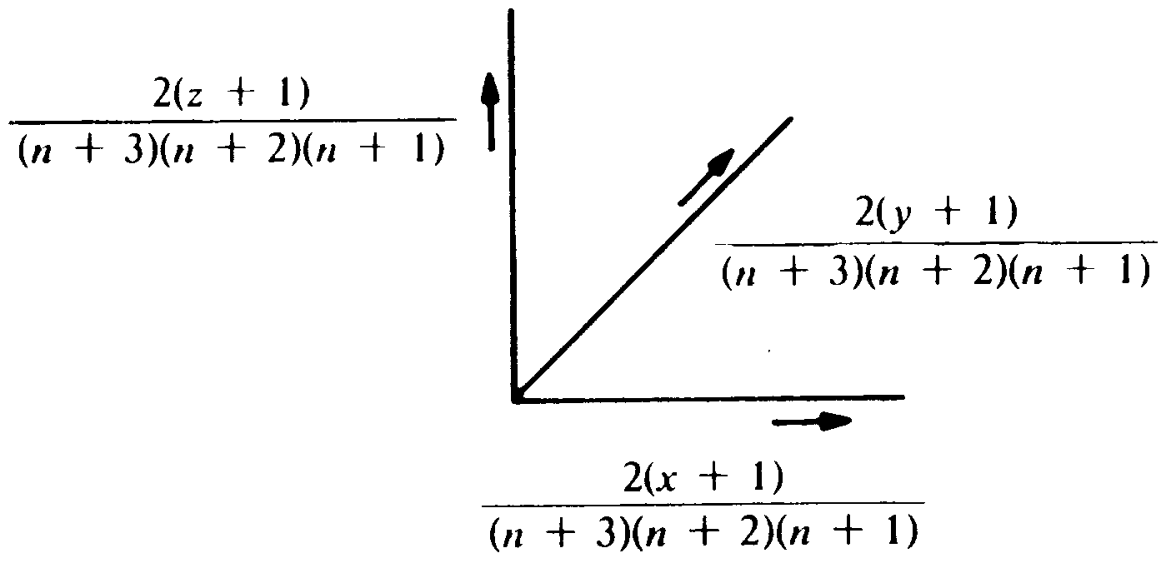}
The total flow out of $(x, y, z)$ is then
\begin{eqnarray*}
&&
\frac{2(x+1)}{(n+3)(n+2)(n+1)}
+
\frac{2(y+1)}{(n+3)(n+2)(n+1)}
+
\frac{2(z+1)}{(n+3)(n+2)(n+1)}
\\&=&
\frac{2}{(n+2)(n+1)}
.
\end{eqnarray*}
The flow into $(x, y, z)$
comes from the points $(x - 1, y, z)$, $(x, y - 1, z)$, $(x, y, z - 1)$
and, hence, the total flow into $(x, y, z)$ is
\[
\frac{2x}{(n+2)(n+1)n}
+
\frac{2y}{(n+2)(n+1)n}
+
\frac{2z}{(n+2)(n+1)n}
=
\frac{2}{(n+2)(n+1)}
.
\]
Thus the net flow for $(x, y, z)$ is 0.
The flow out of $\orig$ is $(1/3)+(1/3)+(1/3)=1$
We have now to check finiteness of energy dissipation. 
The flows coming out of the edges at the $n$th level
are all $\leq 2/(n+1)^2$.
There are $(n + 1)(n + 2)/2$ points a distance $n$ from $\orig$,
and thus there are 
$(3/2)(n+1)(n+2) \leq 3 (n+1)^2$ edges
coming out of the $n$th level.
Thus the energy dissipation $E$ has
\[
E \leq \sum_n 3 (n+1)^2 \braces{\frac{2}{(n+1)^2}}^2
= 12 \sum_n \recip{(n+1)^2}
< \infty
,
\]
and the random walk is transient.

\sect{8.8}{The end}

We have come to the end of our labors,
and it seems fitting to look back
and try to say what it is we have learned.

To begin with,
we have seen how phrasing certain mathematical questions in physical terms
allows us to draw on a large body of physical lore,
in the form of established methods and ways of thought,
and thereby often leads us to the answers to those questions.

In particular, we have seen the utility of considerations involving energy.
In took hundreds of years for the concept of energy to emerge
and take its rightful place in physical theory,
but it is now recognized as perhaps the most fundamental concept
in all of physics.
By phrasing our probabilistic problems in physical terms,
we were naturally led to considerations of energy,
and these considerations showed us the way
through the difficulties of our problems.

As for Polya's theorem and the type problem in general,
we have picked up a bag of tricks,
known collectively as ``Rayleigh's short-cut method'',
which we may expect will allow us to determine the type 
of almost any random walk we are likely to embark on.
In the process,
we have gotten some feeling for the connection between
the dimensionality of a random walk and its type.
Furthermore,
we have settled one of the main questions
likely to occur to someone encountering Polya's theorem,
namely:
``If two walks look essentially the same,
and if one has been shown to be transient,
must not the other also be transient?''

Another question likely to occur to someone contemplating Polya's theorem
is the question raised in \sectref{5.8}:
``Since the lattice $\Z^d$ is in some sense
a discrete analog of a resistive medium
filling all of $\R^d$,
should it not be possible to go quickly and naturally
from the trivial computation of
the resistance to infinity of the continuous medium
to a proof of Polya's theorem?''
Our shorting argument allowed us to do this in the two-dimensional case;
that leaves the case of three (or more) dimensions.
Again, it is considerations of energy
that allow us to make this connection.
The trick is to start with the flow field
that one gets by solving the continuous problem,
and adapt it to the lattice,
so as to get a lattice flow to infinity
having finite dissipation rate.
We leave the working out of this as an exercise,
so as not to rob readers of the fun of doing it for themselves.

\exstart
\ex{8.8.1}{
Give one final proof of Polya's theorem in 3 dimensions
by showing how to adapt the $1/r^2$ radial flow field
to the lattice.
(Hint: ``cubes''.)
}
\exend

\section*{Acknowledgements}
This work is derived from the book
\emph{Random Walks and Electric Networks},
originally published in 1984
by the Mathematical Association of American
in their Carus Monographs series.
We are grateful to the MAA for permitting this work to be
freely redistributed under the terms of the GNU General Public
License.
(See Figure \figref{walksgpl}.)
\prefig{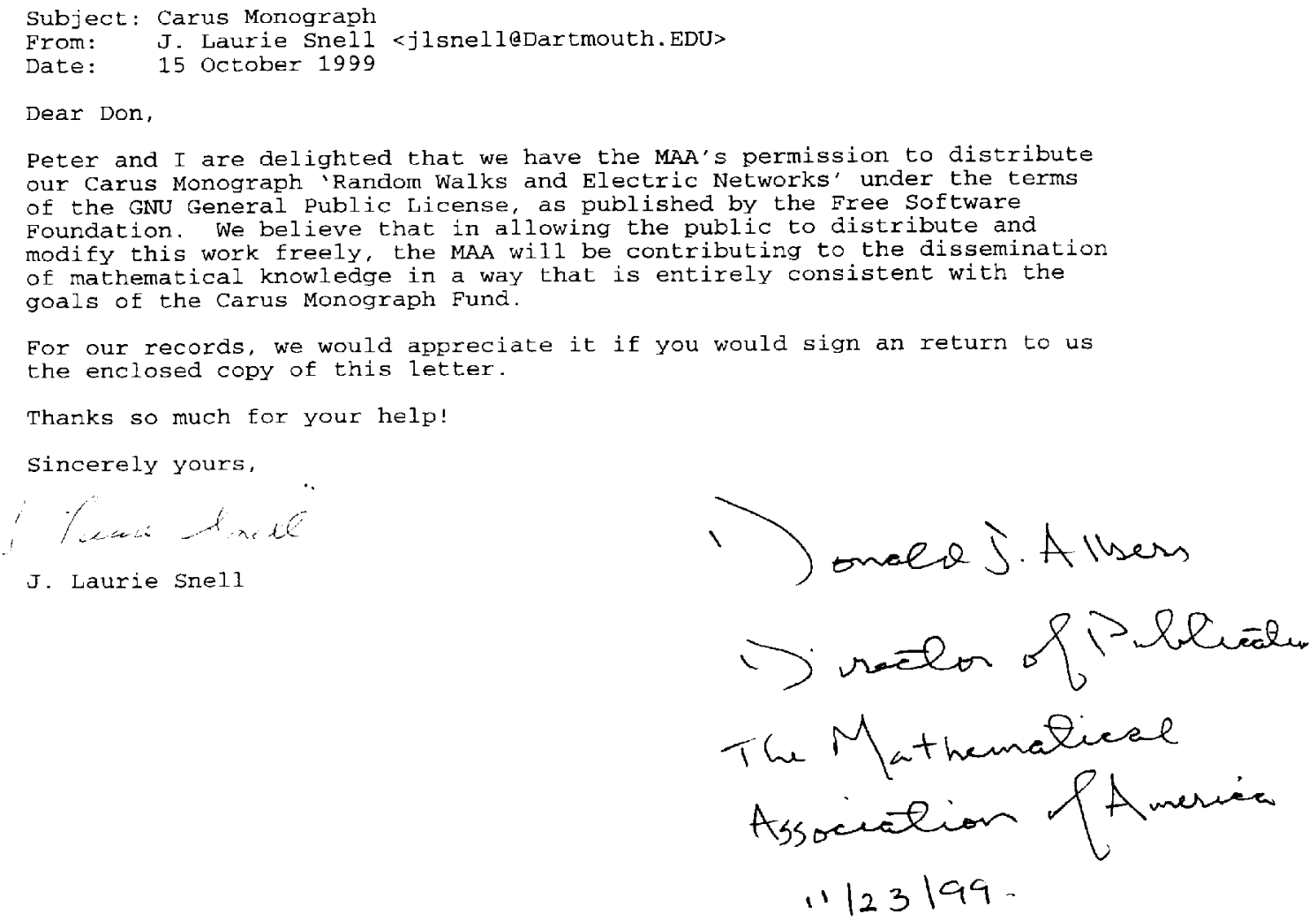}

\newpage
\bibliography{walks}

\begin{thebibliography}{10}

\bibitem{abbott}
Edwin Abbott.
\newblock {\em Flatland}.
\newblock 1899.

\bibitem{bollobas}
B.~Bollob\'{a}s.
\newblock {\em Graph Theory}.
\newblock 1979.

\bibitem{courantfriedrichslewy}
R.~Courant, K.~Friedrichs, and H.~Lewy.
\newblock {\"U}ber die pertiellen {D}ifferenzengleichungen der mathematischen
  {P}hysik.
\newblock {\em Math. Ann.}, 100:32--74, 1928.

\bibitem{doob}
J.~L. Doob.
\newblock {\em Stochastic Processes}.
\newblock 1953.

\bibitem{feller1}
W.~Feller.
\newblock {\em An Introduction to Probability Theory and Its Applications},
  volume~I.
\newblock 1968.

\bibitem{feynman}
R.~P. Feynman.
\newblock {\em The {F}eynmann Lectures on Physics}.
\newblock 1964.

\bibitem{fordFulkerson}
L.~R. Ford and D.~R. Fulkerson.
\newblock Maximal flow through a network.
\newblock {\em Can. J. Math.}, 8:399--404, 1956.

\bibitem{glasserzucker}
M.~L. Glasser and I.~J. Zucker.
\newblock Extended {W}atson integrals for the cubic lattice.
\newblock {\em Proc. Natl. Acad. Sci., USA}, 74:1800--1801, 1977.

\bibitem{griffeathliggett}
D.~Griffeath and T.~M. Liggett.
\newblock Critical phenomena for {S}pitzer's reversible nearest particle
  systems.
\newblock {\em Ann. Probab.}, 10:881--895, 1982.

\bibitem{hershgriego}
R.~Hersh and R.~J. Griego.
\newblock Brownian motion and potential theory.
\newblock {\em Scientific American}, pages 66--74, March 1969.

\bibitem{jeans}
J.~Jeans.
\newblock {\em The Mathematical Theory of Electricity and Magnetism, 5th
  Edition}.
\newblock 1966.

\bibitem{kakutani}
S.~Kakutani.
\newblock Markov processes and the {D}irichlet problem.
\newblock {\em Proc. Jap. Acad.}, 21:227--233, 1945.

\bibitem{kelly}
F.~Kelly.
\newblock {\em Reversibility and Stochastic Networks}.
\newblock 1979.

\bibitem{kemenysnellknapp}
J.~G. Kemeny, J.~L. Snell, and A.~W. Knapp.
\newblock {\em Denumerable {Markov} Chains}.
\newblock 1966.

\bibitem{kemenysnellthompson}
J.~G. Kemeny, J.~L. Snell, and G.~L. Thompson.
\newblock {\em Finite {Markov} Chains, 3rd Edition}.
\newblock 1974.

\bibitem{kesten}
Harry Kesten.
\newblock {\em Percolation Theory for Mathematicians}.
\newblock 1982.

\bibitem{kingman}
J.~F.~C. Kingman.
\newblock Markov population processes.
\newblock {\em J. Appl. Prob.}, 6:1--18, 1969.

\bibitem{lehman}
A.~Lehman.
\newblock A resistor network inequality, {P}roblem 60-5.
\newblock {\em SIAM Review}, 4:150--154, 1965.

\bibitem{levy}
Paul L\'evy.
\newblock {\em Th\'eorie de l'addition des variable al\'eatoires}.
\newblock 1937.

\bibitem{lyons}
T.~J. Lyons.
\newblock A simple criterion for transience of a reversible {M}arkov chain.
\newblock {\em Ann. Probab.}, 11:393--402, 1983.

\bibitem{maxwell}
J.~C. Maxwell.
\newblock {\em Treatise on Electricity and Magnetism, 3rd Edition}.
\newblock 1891.

\bibitem{mccreawhipple}
W.~H. McCrea and F.~J.~W. Whipple.
\newblock Random paths in two and three dimensions.
\newblock {\em Proc. of the Royal Soc. of Edinburgh}, 60:281--298, 1940.

\bibitem{montrollwest}
E.~W. Montroll and B.~J. West.
\newblock Fluctuation phenomena.
\newblock In {\em Studies in Statistical Mathematics}, volume~7, pages 61--175.
  1979.

\bibitem{nashwilliams}
C.~St. J.~A. Nash-Williams.
\newblock Random walk and electric currents in networks.
\newblock {\em Proc. Camb. Phil. Soc.}, 55:181--194, 1959.

\bibitem{onsager}
L.~Onsager.
\newblock Reciprocal relations in irreversible processes {I}.
\newblock {\em Phys. Rev.}, 37:405--426, 1931.

\bibitem{polya:irrfahrt}
G.~Polya.
\newblock {\"{U}}ber eine {A}ufgabe betreffend die {I}rrfahrt im
  {S}trassennetz.
\newblock {\em Math. Ann.}, 84:149--160, 1921.

\bibitem{polya:solve}
G.~Polya.
\newblock {\em How to Solve It, 2nd Edition}.
\newblock 1957.

\bibitem{polyaszego}
G.~Polya and G.~Szego.
\newblock {\em Isoperimetric Inequalities of Mathematical Physics}.
\newblock 1951.

\bibitem{rayleigh}
J.~W.~S. Rayleigh.
\newblock On the theory of resonance.
\newblock In {\em Collected scientific papers}, volume~1, pages 33--75. 1899.

\bibitem{royden}
H.~L. Royden.
\newblock Harmonic functions on open {R}iemann surfaces.
\newblock {\em Trans. Amer. Math. Soc.}, 73:40--94, 1952.

\bibitem{shannonhagelbarger}
C.~E. Shannon and D.~W. Hagelbarger.
\newblock Concavity of resistance function.
\newblock {\em J. of Appl. Phys.}, 27:42--43, 1956.

\bibitem{snell}
J.~L. Snell.
\newblock Probability and martingales.
\newblock {\em The Mathematical Intelligencer}, 4, 1982.

\bibitem{thomsonTait}
W.~Thomson and P.~G. Tait.
\newblock {\em Treatise on Natural Philosophy}.
\newblock 1879.

\bibitem{ville}
J.~Ville.
\newblock {\em \'Etude critique de la notion de collectif}.
\newblock 1939.

\bibitem{watson}
G.~N. Watson.
\newblock Three triple integrals.
\newblock {\em Quarterly J. Math.}, 10:266--276, 1939.

\end{thebibliography}
\bibliographystyle{plain}

\end{document}